\documentclass[reqno]{amsart}

\usepackage{etex}
\usepackage{amsmath,amssymb,amsthm,amsfonts,mathrsfs}
\usepackage[frame,cmtip,arrow,matrix,line,graph,curve]{xy}
\usepackage{graphpap,color,paralist,pstricks}
\usepackage[mathscr]{eucal}
\usepackage[matha,mathb]{mathabx}
\usepackage[pdftex,colorlinks,backref=page,citecolor=blue]{hyperref}
\usepackage{tikz}
\usetikzlibrary{decorations,calc,matrix,cd}
\usepackage{epic,eepic}
\usepackage{yfonts}
\usepackage{enumerate}
\usepackage{blkarray}
\usepackage{mathtools}
\usepackage{setspace}

\usepackage{xcolor}

\definecolor{darkgreen}{rgb}{0,.5,0}
\definecolor{brown}{rgb}{0.5,0.3,0}

\usepackage[T1]{fontenc}
\usepackage[utf8]{inputenc}

\theoremstyle{definition}
\newtheorem{theorem}{Theorem}[section]
\newtheorem{definition}[theorem]{Definition}
\newtheorem{conjecture}[theorem]{Conjecture}
\newtheorem{lemma}[theorem]{Lemma}
\newtheorem{proposition}[theorem]{Proposition}
\newtheorem{corollary}[theorem]{Corollary}

\theoremstyle{remark}
\newtheorem{remark}[theorem]{Remark}
\newtheorem{example}[theorem]{Example}

\makeatletter
\newtheorem*{rep@theorem}{\rep@title} \newcommand{\newreptheorem}[2]{%
\newenvironment{rep#1}[1]{%
\def\rep@title{\bf #2 \ref{##1}}%
\begin{rep@theorem} }%
{\end{rep@theorem} } }
\makeatother
\newreptheorem{theorem}{Theorem}
\newreptheorem{conjecture}{Conjecture}

\numberwithin{equation}{subsection}

\newcommand{\set}[1]{\left\{#1\right\}}
\newcommand{\abs}[1]{\left|#1\right|}
\newcommand{\angl}[1]{\left<#1\right>}

\def\PD{\text{PD}}
\def\PA{P\!A}       
\def\HL{\text{HL}}
\def\HR{\text{HR}}
\def\lef{L}
\def\bl{{p}}             

\DeclareMathOperator{\D}{{\mathscr D}}
\DeclareMathOperator{\bc}{BC}
\DeclareMathOperator{\BC}{\overline{BC}}
\DeclareMathOperator{\IN}{{IN}}
\DeclareMathOperator{\Sym}{{Sym}}
\DeclareMathOperator{\M}{\mathrm M}         
\DeclareMathOperator{\cl}{cl}               
\DeclareMathOperator{\N}{\mathrm N}         

\DeclareMathOperator{\nbc}{{\mathit{nbc}}}  
\DeclareMathOperator{\bnbc}{{\beta-\mathit{nbc}}}  
\DeclareMathOperator{\bcone}{{\beta-\mathrm{cone}}}  
\DeclareMathOperator{\e}{{\mathbf e}}  
\DeclareMathOperator{\f}{{\mathbf f}}  
\DeclareMathOperator{\st}{{st}}  
\DeclareMathOperator{\clst}{{\overline{st}}}  
\DeclareMathOperator{\MW}{{MW}}  
\DeclareMathOperator{\spn}{{span}}  
\DeclareMathOperator{\Eff}{{Eff}} 
\DeclareMathOperator{\TSigma}{{\widetilde{\Sigma}}}  
\def\tallSigma{{\vphantom{\TSigma}\smash{\Sigma}}}  
\DeclareMathOperator{\TDelta}{{\widetilde{\Delta}}}  
\DeclareMathOperator{\hr}{{hr}}  

\DeclareMathOperator{\stellar}{{stellar}}
\DeclareMathOperator{\link}{{link}}

\DeclareMathOperator{\Hom}{{Hom}}
\DeclareMathOperator{\rank}{{rk}}

\DeclareMathOperator{\conv}{{conv}}
\DeclareMathOperator{\trop}{{trop}}
\DeclareMathOperator{\mult}{{mult}}
\DeclareMathOperator{\cone}{{cone}}
\DeclareMathOperator{\closure}{{closure}}

\DeclareMathOperator{\csm}{{csm}}
\DeclareMathOperator{\Kone}{{\mathscr K}}  
\DeclareMathOperator{\ebar}{{\overline{\e}}}  
\newcommand{\xbar}{{\overline{x}}}  

\renewcommand{\P}{{\mathscr P}}
\newcommand{\Q}{{\mathscr Q}}
\newcommand{\PP}{{\mathbb P}}

\newcommand{\B}{{\mathscr B}}
\newcommand{\F}{{\mathscr F}}
\newcommand{\G}{{\mathscr G}}
\renewcommand{\SS}{{\mathscr S}}
\newcommand{\T}{{\mathscr T}}

\newcommand{\Z}{{\mathbb Z}}

\newcommand{\R}{{\mathbb R}}

\newcommand{\C}{{\mathscr C}}
\newcommand{\BB}{{\mathscr B}}

\renewcommand{\k}{{\Bbbk}}
\newcommand{\E}{\overline{E}}

\renewcommand{\i}{\overline{i}}
\renewcommand{\k}{\overline{k}}
\newcommand{\n}{\overline{n}}

\newcommand{\0}{\overline{0}}
\renewcommand{\1}{\overline{1}}
\renewcommand{\2}{\overline{2}}
\newcommand{\3}{\overline{3}}

\newcommand{\PAA}{{A_{\M,\M^\perp}}}  
\newcommand{\SSS}{{S_{\M,\M^\perp}}}  
\newcommand{\III}{{I_{\M,\M^\perp}}}  
\newcommand{\JJJ}{{J_{\M,\M^\perp}}}  
\newcommand{\og}{\overline{\gamma}}

\begin{document}

\title{Lagrangian geometry of matroids}

\author{Federico Ardila}
\address{San Francisco State University and Universidad de Los Andes}
\email{federico@sfsu.edu}

\author{Graham Denham}
\address{University of Western Ontario}
\email{gdenham@uwo.ca}

\author{June Huh}
\address{Princeton University and Korea Institute for Advanced Study}
\email{huh@princeton.edu}

\maketitle

\begin{abstract}
We introduce the conormal fan of a matroid $\M$, which is a Lagrangian analog of the Bergman fan of $\M$.
We use the conormal fan to give a Lagrangian interpretation of the Chern--Schwartz--MacPherson cycle of $\M$. This allows us to express the $h$-vector of the broken circuit complex of $\M$ in terms of the intersection theory of the conormal fan of $\M$.
We also develop general tools for tropical Hodge theory to prove that the conormal fan  satisfies Poincar\'e duality, the hard Lefschetz theorem, and the Hodge--Riemann relations.  The Lagrangian interpretation of the Chern--Schwartz--MacPherson cycle of $\M$,
when combined with the Hodge--Riemann relations for the conormal fan of $\M$, implies Brylawski's and Dawson's conjectures that the $h$-vectors of the broken circuit complex and the independence complex of $\M$ are log-concave sequences.
\end{abstract}

\tableofcontents

\section{Introduction} \label{sec:intro}

\subsection{Geometry of matroids}\label{sec:introgeom}

A \emph{matroid} $\M$ on a finite set $E$ is a nonempty collection of subsets of $E$, called \emph{flats} of $\M$, that satisfies the following properties:
\begin{enumerate}[(1)]\itemsep 5pt
\item The intersection of any two flats is a flat.
\item For any flat $F$, any element in $E- F$ is contained in exactly one flat that is minimal among the
flats strictly containing $F$.
\end{enumerate}
The set $\mathscr{L}(\M)$ of all flats of $\M$ is a geometric lattice, and all geometric lattices arise in this way from a matroid \cite[Chapter 3]{Welsh}.
 The theory of  matroids captures the combinatorial essence shared by  natural notions of independence in linear algebra, graph theory, matching theory, the theory
of field extensions, and the theory of routings, among others.

Gian-Carlo Rota, who helped lay down the foundations of the field,  was one of its most energetic ambassadors.
He rejected the ``ineffably cacophonous" name of matroids, preferring to  call them combinatorial geometries instead \cite{Crapo-Rota}. This alternative name never really caught on, but the geometric roots of the field have since grown much deeper, bearing many new fruits. The geometric approach to matroid theory has recently led to solutions of long-standing conjectures, and to the development of fascinating mathematics at the intersection of combinatorics, algebra, and geometry.

There are at least three useful polyhedral models of a matroid $\M$. For a short survey, see \cite{Ardilasurvey}.
The first one is the \emph{basis polytope} of $\M$ introduced by Edmonds in optimization and Gelfand--Goresky-MacPherson-Serganova in algebraic geometry. 
It reveals an intricate relationship of matroids with the Grassmannian variety and the special linear group.
The second model is the \emph{Bergman fan} of $\M$, introduced by Sturmfels and Ardila--Klivans in tropical geometry.
It was used by Adiprasito--Huh--Katz to prove the log-concavity of the $f$-vectors of the independence complex and the broken circuit complex of $\M$.
The third model, which we call the \emph{conormal fan} of $\M$, is the main character of this paper. We use its intersection-theoretic and Hodge-theoretic properties to prove conjectures of Brylawski \cite{Brylawski82}, Dawson \cite{Dawson83}, 
and Swartz \cite{Swartz03} on the $h$-vectors of  the independence complex and the broken circuit complex of $\M$.

\subsection{Conormal fans and their geometry}

Throughout the paper, we write $r+1$ for the rank of $\M$, write $n+1$ for the cardinality of $E$, and  suppose that $n$ is positive.\footnote{There are exactly two matroids on a single element ground set, the \emph{loop} and the \emph{coloop}, which are dual to each other.
These matroids will play exceptional roles in our inductive arguments.}
Following \cite{MaclaganSturmfels}, we define the \emph{tropical projective torus} of $E$ to be the $n$-dimensional vector space
\[
\N_E=\mathbb{R}^E/ \mathbb{R}\mathbf{e}_E, \qquad  \mathbf{e}_E=\sum_{i \in E} \mathbf{e}_i.
\]
The tropical projective torus is equipped with the  functions
\[
\alpha_j(z)=\max_{i \in E} ( z_j-z_i), \ \ \text{one for each element $j$ of $E$.}
\]
These functions are equal to each other modulo global linear functions on $\N_E$,
and we write $\alpha$ for the common equivalence class of $\alpha_j$.
The \emph{Bergman fan} of $\M$, denoted $\Sigma_{\M}$, is an $r$-dimensional fan  in the $n$-dimensional vector space $\N_E$ whose underlying set  is the \emph{tropical linear space}
\[
\trop(\M)=\Big\{z   \, | \, \text{$\min_{i \in C}(z_i)$ is achieved at least twice for every circuit $C$ of $\M$}\Big\} \subseteq \N_E.
\]
It is a subfan of the \emph{permutohedral fan} $\Sigma_E$ cut out by the hyperplanes $x_i = x_j$ for each pair of distinct elements $i$ and $j$ in $E$.
 This is the normal fan of the \emph{permutohedron} $\Pi_E$.
The functions $\alpha_j$ are piecewise linear on the permutohedral fan, and hence piecewise linear on the Bergman fan of $\M$.\footnote{A continuous function $f$ is said to be \emph{piecewise linear} on a fan $\Sigma$ if the restriction of $f$ to any cone in $\Sigma$ is linear.
In this case, we  say that the fan $\Sigma$ \emph{supports} the piecewise linear function $f$.}

Tropical linear spaces are central objects in tropical geometry:
For any linear subspace $V$ of  $\mathbb{C}^E$,
the tropicalization of the intersection of $\mathbb{P}(V)$ with the torus of $\mathbb{P}(\mathbb{C}^E)$ is the tropical linear space of the linear matroid on $E$ represented by $V$ \cite{Sturmfels}.
Furthermore, tropical linear spaces are precisely the tropical fans of degree one with respect to $\alpha$, that is, the tropical analogs of linear spaces \cite{Fink}.
Tropical manifolds are thus defined to be spaces that locally look like Bergman fans of matroids \cite{IKMZ}.

Adiprasito, Huh, and Katz showed that the Chow ring of the Bergman fan of ${\M}$
satisfies Poincar\'e duality, the hard Lefschetz theorem, and the Hodge--Riemann relations \cite{AHK15}. Furthermore, they interpreted the entries of the $f$-vector of the reduced broken circuit complex of $\M$  -- an invariant of the matroid generalizing the chromatic polynomial for graphs -- as intersection numbers in the Chow ring of $\Sigma_{\M}$.
The geometric interpretation then implied the log-concavity of the coefficients of the \emph{characteristic polynomial} and the \emph{reduced characteristic polynomial}
\[
\chi_{\M}(q)\coloneq \sum_{F \in \mathscr{L}(\M)} \mu(\varnothing,F) q^{\text{corank}(F)}, \qquad \overline{\chi}_{\M}(q) \coloneq  \chi_{\M}(q)/(q-1),
\]
where $\mu$ is the \emph{M\"obius function} on the geometric lattice $\mathscr{L}(\M)$ for a loopless matroid $\M$.\footnote{If $\M$ has a loop, by definition, the characteristic polynomial and the reduced characteristic polynomial of $\M$ are zero.}

The conormal fan  $\Sigma_{\M,\M^\perp}$ is an alternative polyhedral model for $\M$.
Its construction uses the dual matroid $\M^\perp$, the matroid on $E$ whose bases are the complements of  bases of $\M$.
We refer to \cite{Oxley} for background on matroid duality and other general facts on matroids.
A central role is played by  the \emph{addition map}
\[
\N_{E,E}\coloneq  \N_E \oplus \N_{E} \longrightarrow \N_E, \qquad (z,w) \longmapsto z+w.
\]
The  function $\alpha_j$ on $\N_E$ pulls back to a  function $\delta_j$ on $\N_{E,E}$ under the addition map.
Explicitly,
\[
\delta_j(z,w)=\max_{i \in E} ( z_j+w_j-z_i-w_i).
\]
The function $\delta_j$ is piecewise linear on a fan that we construct, called the \emph{bipermutohedral fan} $\Sigma_{E,E}$ (Section \ref{sec:ourfan}). This is the normal fan of a convex polytope $\Pi_{E,E}$ that we call the \emph{bipermutohedron}.
The functions $\delta_j$ for $j$ in $E$ are equal to each other modulo global linear functions on $\N_{E,E}$,
and we write $\delta$ for their common equivalence class.

The \emph{cotangent fan} $\Omega_E$ is the  subfan of the bipermutohedral fan $\Sigma_{E,E}$ whose underlying set  is the tropical hypersurface
\[
\trop(\delta)=\Big\{ (z,w)  \ | \  \text{$\min_{i \in E}\{z_i+w_i\}$ is achieved at least twice}\Big\} \subseteq \N_{E,E}.
\]
We show in Section \ref{ss:conormal} that, for any matroid $\M$ on $E$, we have
\[
\trop(\M) \times \trop(\M^\perp) \subseteq \trop(\delta).
\]
The \emph{conormal fan}  $\Sigma_{\M,\M^\perp}$ is defined to be the subfan of the cotangent fan $\Omega_E$ that subdivides
the product $\trop(\M) \times \trop(\M^\perp)$. For our purposes, it is necessary to work with the conormal fan of $\M$
instead of the product of the Bergman fans of $\M$ and $\M^\perp$, because the function $\delta_j$ need not be piecewise linear on the product of the Bergman fans.

The projections to the summands of $\N_{E,E}$ define morphisms of fans\footnote{A \emph{morphism} from a fan $\Sigma_1$ in $\N_1=\mathbb{R} \otimes \N_{1,\mathbb{Z}}$ to a fan $\Sigma_2$ in $\N_2 =\mathbb{R} \otimes \N_{2,\mathbb{Z}}$ is an integral linear map from $\N_1$ to $\N_2$ such that the image of any cone in $\Sigma_1$ is a subset of a cone in $\Sigma_2$. In the context of toric geometry, a morphism from $\Sigma_1$ to $\Sigma_2$ can be identified with a toric morphism from the toric variety of $\Sigma_1$ to the toric variety of $\Sigma_2$  \cite[Chapter 3]{CLS}.}
\[
\pi\colon \Sigma_{\M,\M^\perp} \longrightarrow \Sigma_{\M} \quad \text{and} \quad  \overline\pi\colon \Sigma_{\M,\M^\perp} \longrightarrow \Sigma_{\M^\perp}.
\]
Thus, in addition to the functions $\delta_j$, the conormal fan of $\M$ supports the pullbacks of $\alpha_j$ on $\Sigma_{\M}$ and $\overline{\alpha}_j$ on $\Sigma_{\M^\perp}$, which are the piecewise linear functions
\[
\gamma_j(z,w)= \max_{i \in E} (z_j-z_i) \quad \text{and} \quad
\og_j(z,w)= \max_{i \in E} (w_j-w_i).
\]
These define the equivalence classes $\gamma$ and $\og$ of functions on $\N_{E,E}$.

The conormal fan is a tropical analog of the incidence variety appearing in the classical theory of projective duality.
For a subvariety $X$ of a projective space $\mathbb{P}(V)$,
the incidence variety $\mathscr{I}_X$ is a subvariety of the product of $\mathbb{P}(V)$ with the dual projective space $\mathbb{P}(V^\vee)$ that projects onto $X$ and its dual $X^\vee$.
Over the smooth locus of $X$, the incidence variety $\mathscr{I}_X$ is the total space of the projectivized conormal bundle of $X$ and,
over the smooth locus of $X^\vee$, it is the total space of  the projectivized conormal bundle of $X^\vee$.\footnote{Thus, to be precise, the conormal fan is a tropical analog of the projectivized conormal variety and the cotangent fan is a tropical analog of the projectivized cotangent space.
We trust that the omission of the term ``projectivized'' will cause no confusion.}
It is the projectivization of a conic Lagrangian subvariety of $V \times V^\vee$,
and any conic Lagrangian subvariety of $V \times V^\vee$ arise in this way.
We refer to \cite[Chapter 1]{GKZ} for a modern exposition of the theory of projective duality.

We use the conormal fan of $\M$ to give a geometric interpretation of the polynomial $\overline{\chi}_{\M}(q+1)$,
whose coefficients form the $h$-vector of the broken circuit complex of $\M$ with alternating signs.
In particular, we give a geometric formula for \emph{Crapo's beta invariant}
\[
\beta_{\M} \coloneq  (-1)^{r}\ \overline{\chi}_{\M}(1).
\]
This new tropical geometry is inspired by the Lagrangian geometry of conormal varieties in classical algebraic geometry, as we now explain.

Consider the category of complex algebraic varieties with proper morphisms.
According to a conjecture of Deligne and Grothendieck,
there is a unique natural transformation ``$\csm$'' from the functor of constructible functions on complex algebraic varieties
to the homology of complex algebraic varieties such that, for any smooth and 
complete variety $X$,
\[
\csm(1_X) = c(TX) \cap [X]=(\text{the total homology Chern class of the tangent bundle of $X$}).
\]
The conjecture was proved by MacPherson \cite{MacPherson74},
and it was recognized later in \cite{BS} that the class $\csm(1_X)$, for possibly singular $X$, coincides
with a class constructed earlier by Schwartz \cite{Schwartz1}.
For any constructible subset $X$ of $Y$,   the $k$-th Chern--Schwartz--MacPherson class of $X$ in $Y$ is the homology class
\[
\csm_k(1_X) \in H_{2k}(Y).
\]
Aiming to introduce a tropical analog of this theory, L\'opez de Medrano, Rinc\'on, and Shaw
introduced the Chern--Schwartz--MacPherson cycle of the Bergman fan of ${\M}$ in  \cite{LRS17}:
The \emph{$k$-th Chern--Schwartz--MacPherson cycle} of $\M$ is, by definition, the weighted fan $\csm_k(\M)$ supported on the $k$-dimensional skeleton of $\Sigma_{\M}$
with the weights
\[
w(\sigma_\mathscr{F})=(-1)^{r-k}  \prod_{i=0}^{k} \beta_{\M(i)},
\]
where $\sigma_\mathscr{F}$ is the $k$-dimensional cone corresponding to a flag of flats $\mathscr{F}$ of $\M$ 
and $\M(i)$ is the minor of $\M$ corresponding to  the $i$-th interval in $\mathscr{F}$. 
This weighted fan  behaves well combinatorially and geometrically.
First,
the weights satisfy the balancing condition in tropical geometry \cite[Theorem 1.1]{LRS17},
so that we may view the Chern--Schwartz--MacPherson cycle as a Minkowski weight
\[
\csm_k(\M) \in \MW_k(\Sigma_{\M}).
\]
Second, when $\trop(\M)$ is the tropicalization of  the intersection $\mathbb{P}(V) \cap (\mathbb{C}^*)^E/\mathbb{C}^*$, the Minkowski weight can be identified with the $k$-th Chern--Schwartz--MacPherson class of $\mathbb{P}(V) \cap (\mathbb{C}^*)^E/\mathbb{C}^*$ in the toric variety of the permutohedron $\Pi_E$ \cite[Theorem 1.2]{LRS17}.
Third, the Chern-Schwartz-MacPherson cycles of $\M$ satisfy a deletion-contraction formula, a matroid version of the inclusion-exclusion principle \cite[Proposition 5.2]{LRS17}.
It follows that the degrees of these Minkowski weights determine the reduced characteristic polynomial of $\M$ by the formula
\[
\overline{\chi}_{\M}(q+1)=\sum_{k=0}^r \deg(\csm_k(\M)) q^k,
\]
where the degrees are taken with respect to the class $\alpha$ \cite[Theorem 1.4]{LRS17}.
Fourth,  the Chern-Schwartz-MacPherson cycles of matroids can be used to define Chern classes of smooth tropical varieties.
In codimension $1$, the class  agrees with the anticanonical divisor of a tropical variety defined by Mikhalkin in \cite{Mikhalkin}.
For smooth tropical surfaces, these classes agree with the Chern classes of tropical surfaces introduced in \cite{Cartwright} and \cite{Shaw}  to formulate  Noether's formula for tropical surfaces.

Schwartz's and MacPherson's constructions of $\csm$ for complex algebraic varieties are rather subtle. Sabbah  later observed that the Chern-Schwartz-MacPherson classes can be interpreted more simply as ``shadows'' of the characteristic cycles in the cotangent bundle.
Sabbah summarizes the situation in the following quote from \cite{Sabbah}:
\begin{quote}
\emph{la th\'eorie des classes de Chern de [Mac74] se ram\`ene \`a une th\'eorie de Chow sur $T^*X$, qui ne fait intervenir que des classes fondamentales.}
\end{quote}
The functor of constructible functions is replaced with a functor of Lagrangian cycles of $T^*X$,
which are exactly the linear combinations of the \emph{conormal varieties} of the subvarieties of $X$.
In the Lagrangian framework,   key operations on constructible functions become more geometric.

Similarly, L\'opez de Medrano, Rinc\'on, and Shaw's original definition of the Chern--Schwartz--MacPherson cycles of a matroid $\M$ is combinatorially intricate. We prove that they are ``shadows'' of  much simpler cycles under the pushforward map
\[
\pi_{*}\colon \MW_k(\Sigma_{\M, \M^\perp}) \longrightarrow \MW_k(\Sigma_{\M}) .
\]
See Section \ref{sec:homology} for a review of basic tropical intersection theory.

\begin{theorem}\label{CSMTheorem}
When $\M$ has no loops and no coloops, we have 
 \[
  \csm_k(\M)=(-1)^{r-k} \pi_{*} (\delta^{n-k-1} \cap 1_{\M,\M^\perp} ) \ \ \text{for $0 \le k \le r$},
 \]
where  $1_{\M,\M^\perp}$ is the top-dimensional constant Minkowski weight $1$ on the conormal fan of $\M$.
 \end{theorem}

It follows from Theorem \ref{CSMTheorem} and  the projection formula that the reduced characteristic polynomial of $\M$ can be expressed in terms of the intersection theory of the conormal fan as follows:

\begin{theorem}\label{degtheorem}
When $\M$ has no loops and no coloops, we have
\[
\overline{\chi}_{\M}(q+1)=\sum_{k=0}^r (-1)^{r-k} \deg( \gamma^k\, \delta^{n-k-1}) \, q^k,
\]
where the degrees are taken with respect to the top-dimensional constant Minkowski weight $1_{\M,\M^\perp}$ on the conormal fan.
\end{theorem}

When ${\M}$ is representable over $\mathbb{C}$,\footnote{
We say that $\mathrm{M}$ is \emph{representable} over a field $\mathbb{F}$ if there
exists a linear subspace $V\subseteq \mathbb{F}^E$ such that $S\subseteq E$ is independent in $\M$
if and only if the projection from $V$ to $\mathbb{F}^S$ is surjective.
Almost all matroids are not representable over any field \cite{Nelson}.}
the third author gave an algebro-geometric version of  Theorem \ref{CSMTheorem} in \cite{HuhML}.
The complex geometric version of the identity boils down to  the general fact that the Chern--Schwartz--MacPherson class of a smooth variety $X$ in its normal crossings compactification $Y$
is the total Chern class of the logarithmic tangent bundle:
\[
\csm(1_X)=c(\Omega_Y^1(\log Y - X)^\vee) \cap [Y].
\]
In fact, the logarithmic formula can be used to construct the natural transformation $\csm$ \cite{AluffiConstruction}. For precursors of the logarithmic viewpoint, see  \cite{Aluffi} and \cite{GoreskyPardon}. The current paper demonstrates that a similar geometry exists for arbitrary tropical linear spaces.

\subsection{Inequalities for matroid invariants}\label{sec:invariants}

Let $a_0, a_1, \ldots, a_n$ be a sequence of  nonnegative integers,
and let $d$ be the largest index with nonzero  $a_d$.
\begin{enumerate}[$\bullet$]\itemsep 5pt
\item The sequence is said to be \emph{unimodal} if
\[
 a_0 \leq a_1 \leq \cdots \leq a_{k-1} \leq a_k \geq a_{k+1} \geq \cdots \geq a_{n} \ \ \text{for some $0 \leq k \leq n$.}
\]
\item The sequence is said to be \emph{log-concave} if
\[
a_{k-1}a_{k+1} \leq a_k^2 \ \ \text{for all $0<k<n$.}
\]
\item The sequence is said to be \emph{flawless} if
\[
a_k \le a_{d-k} \ \ \text{for all $0 \le k \le d/2$.}
\]
\end{enumerate}
Many enumerative sequences are conjectured to have these properties, but proving them often turns out to be difficult.
Combinatorialists have been interested in these conjectures because their solution typically requires a fundamentally new construction or connection with a distant field, thus revealing hidden structural information about the objects in question.
For surveys of known results and open problems, see \cite{Brenti} and  \cite{StanleyLU,StanleyP}.

A  \emph{simplicial complex} $\Delta$ is a collection of subsets of a finite set, called \emph{faces} of $\Delta$,
that is downward closed.
The \emph{face enumerator} of $\Delta$ and the \emph{shelling polynomial} of $\Delta$ are the polynomials
\[
f_\Delta(q)=\sum_{S\in \Delta} q^{d-\abs{S}+1}=
\sum_{k \ge 0} f_k(\Delta) q^{d-k+1} \quad \text{and} \quad
h_\Delta(q)=f_\Delta(q-1)=\sum_{k \ge 0} h_k(\Delta)q^{d-k+1},
\]
where $d$ is the dimension of $\Delta$.
The \emph{$f$-vector} of a simplicial complex is the sequence of coefficients of its face enumerator,
and the \emph{$h$-vector} of a simplicial complex is the sequence of coefficients of its shelling polynomial.
When $\Delta$ is \emph{shellable},\footnote{An $r$-dimensional pure simplicial complex is said to be shellable if there is an ordering of its facets such that each facet intersects the simplicial complex generated by its predecessors in a pure $(r-1)$-dimensional complex.}
 the shelling polynomial of $\Delta$ enumerates the facets used in shelling $\Delta$,
and hence the $h$-vector of $\Delta$ is nonnegative.

We study the $f$-vectors and $h$-vectors of the following shellable simplicial complexes associated to $\M$.
For a gentle introduction, and for the proof of their shellability, see \cite{Bjorner}.
\begin{enumerate}[$\bullet$]\itemsep 5pt
\item The \emph{independence complex} $\IN(\M)$, the collection of subsets of $E$ that are independent in $\M$.
\item The \emph{broken circuit complex} $\bc(\M)$, the collection of subsets of $E$ which do not contain any broken circuit of $\M$.
\end{enumerate}
Here a \emph{broken circuit} is a subset  obtained from a circuit of $\M$ by deleting the least element relative to a fixed ordering of $E$.
The notion was developed by Whitney \cite{Whitney}, Rota \cite{Rota}, Wilf \cite{Wilf}, and Brylawski \cite{BrylawskiBroken}, for the ``chromatic'' study of matroids.
The $f$-vector and the $h$-vector of the broken circuit complex of $\M$ are determined by the characteristic polynomial of $\M$,
and in particular they do not depend on the chosen ordering of $E$:
\[
\chi_{\M}(q)=\sum_{k=0}^{r+1} (-1)^k f_k(\bc(\M)) q^{r-k+1}, \qquad \chi_{\M}(q+1)=\sum_{k=0}^{r+1} (-1)^k h_k(\bc(\M)) q^{r-k+1}.
\]

\begin{conjecture}\label{NumericalConjectures}
The following holds for any matroid $\M$.
\begin{enumerate}[(1)]\itemsep 5pt
\item The $f$-vector of $\IN(\M)$ is  unimodal, log-concave, and flawless.
\item\label{1.3.2} The $h$-vector of $\IN(\M)$ is  unimodal, log-concave, and flawless.
\item The $f$-vector of $\bc(\M)$ is  unimodal, log-concave, and flawless.
\item\label{1.3.4} The $h$-vector of $\bc(\M)$ is  unimodal, log-concave, and flawless.
\end{enumerate}
\end{conjecture}

Welsh \cite{Welsh71} and Mason \cite{Mason72} conjectured the log-concavity of the $f$-vector of the independence complex.\footnote{In \cite{Mason72}, Mason proposed a stronger conjecture that the $f$-vector of the independence complex of $\M$ satisfies
\[
\frac{f_k^{2}}{{n+1\choose k}^2} \ge \frac{f_{k-1}}{{n+1 \choose k-1}} \frac{f_{k+1}}{{n+1 \choose k+1}} \ \ \text{for all $k$.}
\]
In \cite{Brylawski82}, Brylawski conjectures the same set of inequalities for the $f$-vector of the broken circuit complex of $\M$.
Mason's stronger conjecture was recently proved in  \cite{ALOGVb} and  \cite{BHa,BHb}.
An extension of the same result to matroid quotients was obtained in \cite{EH}.
}
Dawson  conjectured the log-concavity of the $h$-vector of the independence complex in \cite{Dawson83},
and independently, Colbourn conjectured  the same  in \cite{Colbourn} in the context of network reliability.
Hibi conjectured  that the $h$-vector of the independence complex must be flawless \cite{Hibi92}.
The unimodality and the log-concavity conjectures for the $f$-vector of the broken circuit complex are due to Heron \cite{Heron}, Rota \cite{Rota70}, and Welsh \cite{Welsh}.
The same conjectures for the chromatic polynomials of graphs were given earlier by Read \cite{Read} and Hoggar \cite{Hoggar74}.
We refer to \cite[Chapter 8]{White} and \cite[Chapter 15]{Oxley} for overviews and historical accounts.
 Brylawski \cite{Brylawski82} conjectured the log-concavity of the $h$-vector of the broken circuit complex.\footnote{In \cite{Brylawski82}, Brylawski proposed a stronger conjecture that the $h$-vector of the broken circuit complex of $\M$ satisfies
 \[
\frac{h_k^{2}}{{n-k \choose n-r-1}^2} \ge \frac{h_{k-1}}{{n-k+1 \choose n-r-1}} \frac{h_{k+1}}{{n-k-1 \choose n-r-1}} \ \ \text{for all $k$.}
\]}
 That the $h$-vector of the broken circuit complex is flawless was
 stated as an open problem in \cite{Swartz03}
and reproduced in \cite{JKLe16} as a conjecture.
We deduce all the above statements using the geometry of conormal fans.

\begin{theorem}\label{thm:conjs}
Conjecture \ref{NumericalConjectures} holds.
\end{theorem}

We prove  the log-concavity of the $h$-vector of the broken circuit complex using Theorem \ref{CSMTheorem}.
This log-concavity  implies all other statements in Conjecture \ref{NumericalConjectures},
thanks to the following known observations:
\begin{enumerate}[$\bullet$]\itemsep 5pt
\item For any simplicial complex $\Delta$, the log-concavity of the $h$-vector implies the log-concavity of the $f$-vector  \cite[Corollary 8.4]{Brenti}.
\item For any pure simplicial complex $\Delta$, the $f$-vector of $\Delta$ is flawless.
 More generally, any pure O-sequence\footnote{A sequence of nonnegative integers $h_0,h_1,\ldots$ is an \emph{O-sequence} if there is an order ideal of monomials $\mathscr{O}$ such that $h_k$ is the number of degree $k$ monomials in $\mathscr{O}$. The sequence is a \emph{pure O-sequence} if the order ideal $\mathscr{O}$ can be chosen so that all the maximal monomials in $\mathscr{O}$ have the same degree. See \cite{pureO} for a comprehensive survey of pure O-sequences.}
  is flawless \cite[Theorem 1.1]{Hibi89}.  
\item For any shellable simplicial complex $\Delta$, the $h$-vector of $\Delta$ has no internal zeros, being an O-sequence \cite[Theorem 6]{StanleyCM}.
Therefore, if the $h$-vector of $\Delta$ is log-concave, then it is unimodal.
\item The broken circuit complex of $\M$ is the cone over the \emph{reduced broken circuit complex} of $\M$, and the two simplicial complexes share the same $h$-vector.
The independence complex of $\M$ is the reduced broken circuit complex of another matroid,
the \emph{free dual extension} of $\M$  \cite[Theorem 4.2]{BrylawskiBroken}. \item   If the $h$-vector of the broken circuit complex of $\M$ is unimodal for all $\M$, then
the $h$-vector of the broken circuit complex of $\M$ is flawless for all $\M$ \cite[Theorem 1.2]{JKLe16}.
\end{enumerate}

\noindent {\bf Previous work.}
The log-concavity of the $f$-vector of the broken circuit complex  was proved in \cite{Huh12} for matroids representable over a field of characteristic $0$.
The result was extended  to matroids representable over some field in \cite{HK12}
and to all matroids in \cite{AHK15}.
An alternative proof of the same fact using the volume polynomial of a matroid was obtained in \cite{BES}.
It was observed in \cite{Lenz13} that the log-concavity of the $f$-vector of the broken circuit complex implies that of the independence complex.

For  matroids representable over a field of characteristic $0$, the log-concavity of the $h$-vector of the broken circuit complex was proved in \cite{Huh15}.
The algebraic geometry behind the log-concavity of the $h$-vector, which became a model for the Lagrangian geometry of conormal fans in the present paper,
was explored in \cite{DGS} and \cite{HuhML}.
In \cite{JKLe16},  Juhnke-Kubitzke and Le used the result of \cite{Huh15}  to deduce that the $h$-vector of the broken circuit complex  is flawless for matroids representable over a field of characteristic $0$.
The flawlessness of the $h$-vector of the independence complex was first proved by Chari using a combinatorial decomposition of the independence complex \cite{Chari97}.
The result was recovered by Swartz \cite{Swartz03} and Hausel \cite{Hausel}, who obtained stronger algebraic results.
The other cases of Conjecture \ref{NumericalConjectures} remained open.

Our solution of Conjecture \ref{NumericalConjectures} was announced in \cite{Ardilasurvey}. 
Shortly after this paper appeared on \url{https://arxiv.org/abs/2004.13116}, Berget, Spink, and Tseng \cite{BST20} have
announced an alternative proof of the log-concavity of the $h$-vector of the independence complex (Conjecture \ref{NumericalConjectures}.\ref{1.3.2}). The relationship between our approach and theirs is still to be understood.
The $h$-vector of the broken circuit complex (Conjecture \ref{NumericalConjectures}.\ref{1.3.4}) is not currently accessible through the alternative method.

\subsection{Tropical Hodge theory}\label{sec:tropicalHodge}

Let us discuss in more detail the strategy of \cite{AHK15}
that led to the log-concavity of the $f$-vector of the broken circuit complex of $\M$.
For the moment, suppose that
there is a linear subspace $V \subseteq \mathbb{C}^E$  representing $\M$ over $\mathbb{C}$,
and consider the variety
\[
Y_V=\text{the closure of  $\PP(V)\cap (\mathbb{C}^*)^E/\mathbb{C}^*$
 in the toric variety of the permutohedron $X(\Sigma_E)$.}\footnote{Throughout the paper, the toric variety of a fan in $\N_E$ refers to the one constructed with respect to the lattice $\mathbb{Z}^E/\mathbb{Z}$.
Similarly, the toric variety of a fan in $\N_{E,E}$ refers to the one constructed with respect to the lattice $\mathbb{Z}^E/\mathbb{Z} \oplus \mathbb{Z}^E/\mathbb{Z}$.}
\]
 If nonempty, $Y_V$ is an $r$-dimensional smooth projective complex variety  which is, in fact,
 contained in the torus invariant open subset of $X(\Sigma_E)$ corresponding to the Bergman fan of $\M$:
\[
Y_V \subseteq X(\Sigma_{\M}) \subseteq X(\Sigma_E).
\]
The work of Feichtner and Yuzvinsky \cite{FY04}, which builds upon the work of De Concini and Procesi
\cite{DP95}, reveals that the inclusion maps induce isomorphisms between integral cohomology and Chow
rings:
\[
H^{2\bullet}(Y_V,\Z)\simeq A^\bullet(Y_V,\Z)\simeq
A^\bullet(X(\Sigma_{\M}),\Z).
\]
As a result, the Chow ring of the $n$-dimensional variety $X(\Sigma_{\M})$ has the
 structure of the even part of the cohomology ring of an $r$-dimensional smooth projective variety.
Remarkably,  this structure on the Chow ring of $X(\Sigma_{\M})$  persists for  any matroid $\M$, even
if $\M$ does not admit any representation over any field.
In particular, the Chow ring of $X(\Sigma_{\M})$ satisfies the Poincar\'e duality, the hard Lefschetz theorem, and the Hodge--Riemann relations \cite{AHK15}.
For a simpler proof of the three properties of the Chow ring, based on its semismall decomposition, see \cite{BHMPW}.

For any simplicial fan $\Sigma$, 
 let $A(\Sigma)$ be the ring of real-valued piecewise polynomial functions on $\Sigma$ modulo the ideal of the linear functions on $\Sigma$,
and let $\Kone(\Sigma)$ be the cone of \emph{strictly convex} piecewise linear functions on $\Sigma$ (Definition \ref{def:KSigma}).

\begin{definition}\label{def:lefschetz}
A $d$-dimensional simplicial fan $\Sigma$ is \emph{Lefschetz} if it satisfies the following.
\begin{enumerate}[(1)]\itemsep 5pt
\item (Fundamental weight) The group of $d$-dimensional Minkowski weights on $\Sigma$ is generated by a positive Minkowski weight $w$.
We write $\deg$ for the corresponding linear isomorphism
\[
\deg\colon A^d(\Sigma) \longrightarrow \R, \qquad \eta  \longmapsto  \eta \cap w.
\]
\item (Poincar\'e duality) For any $0 \le k \le d$, the bilinear map of the multiplication
  \[
    \begin{tikzcd}
      A^k(\Sigma) \times A^{d-k}(\Sigma) \ar[r] &
      A^d(\Sigma) \ar[r,"\deg"] & \R
    \end{tikzcd}
    \]
    is nondegenerate.

\item (Hard Lefschetz property) For any $0\leq k\leq \frac{d}{2}$ and any $\ell \in \Kone(\Sigma)$,
the multiplication map
\[
A^k(\Sigma) \longrightarrow A^{d-k}(\Sigma), \qquad \eta \longmapsto \ell^{d -2k} \eta
\]
is a linear isomorphism.

\item (Hodge--Riemann relations)   For any $0\leq k\leq \frac{d}{2}$ and any $\ell \in \Kone(\Sigma)$,
 the bilinear form
 \[
 A^k(\Sigma) \times A^k(\Sigma) \longrightarrow \R, \qquad (\eta_1,\eta_2) \longmapsto (-1)^k\deg(\ell^{d-2k}  \eta_1 \eta_2)
 \]
is positive definite when restricted to the kernel of the multiplication map $\ell^{d-2k+1}$.
\item (Hereditary property) For any $0< k \le d$ and any $k$-dimensional cone $\sigma$ in $\Sigma$, the star of $\sigma$ in $\Sigma$ is  a
Lefschetz fan of dimension $d-k$.
\end{enumerate}
\end{definition}

The Hodge--Riemann relations give analogs of the Alexandrov--Fenchel  inequality amongst degrees of products of convex piecewise linear functions $\ell_1,\ell_2,\ldots,\ell_d$ on $\Sigma$:
\[
\deg(\ell_1 \ell_2 \ell_3 \cdots \ell_d)^2 \ge \deg(\ell_1 \ell_1 \ell_3 \cdots \ell_d) \deg(\ell_2 \ell_2 \ell_3 \cdots \ell_d).
\]
The Bergman fan  of a matroid ${\M}$ is Lefschetz, and
the log-concavity
of the $f$-vector of the broken circuit complex of $\M$ follows from the Hodge--Riemann relations for the Bergman fan of $\M$ \cite{AHK15}.

We establish the log-concavity of the $h$-vector of the broken circuit complex of $\M$
in the same way,  using the conormal fan of $\M$ in place of the
Bergman fan of $\M$.
Theorem \ref{degtheorem} relates the intersection theory of the conormal fan of $\M$
 to the $h$-vector of the broken circuit complex of $\M$ via the Chern-Schwartz-MacPherson cycles of $\M$.
In order to proceed, we need to show that the conormal fan of $\M$ is  Lefschetz.
We obtain this from the following general result.

\begin{theorem}\label{thm:lefschetz}
  Let $\Sigma_1$ and $\Sigma_2$ be simplicial 
  fans  that have the same support $|\Sigma_1|=|\Sigma_2|$.
If $\Kone(\Sigma_1)$ and $\Kone(\Sigma_2)$
  are nonempty, then $\Sigma_1$ is Lefschetz if and only if $\Sigma_2$ is
  Lefschetz.
\end{theorem}

Theorem \ref{thm:lefschetz} implies, for example, that the reduced normal fan of any  simple polytope is Lefschetz, because the reduced normal fan of the standard simplex is Lefschetz.\footnote{McMullen gave an elementary proof of this fact in \cite{McMullen}.
See \cite{Timorin} and \cite{FK} for alternative presentations.
Our proof of Theorem \ref{thm:lefschetz} is  modeled on these arguments.
Theorem~\ref{thm:lefschetz} gives another
proof of the necessity of McMullen's
bounds~\cite{McMullen} on the face numbers of simplicial polytopes.
In the context of matroid theory, the authors of \cite{GS} used a similar argument to show that any unimodular fan whose support is a tropical linear space satisfies Poincar\'e duality.}
In the context of matroid theory, Theorem \ref{thm:lefschetz} implies that the conormal fan of $\M$ is Lefschetz,
because the Bergman fans of $\M$ and $\M^\perp$ are Lefschetz and the product of Lefschetz fans is  Lefschetz.
When $\Kone(\Sigma)$ is empty, the hard Lefschetz property and the Hodge--Riemann relations for $\Sigma$ hold vacuously. 
The proof of Theorem \ref{thm:lefschetz} shows that, if two simplicial 
fans $\Sigma_1$ and $\Sigma_2$ that have the same support $|\Sigma_1|=|\Sigma_2|$,
then $\Sigma_1$ satisfies Poincar\'e duality if and only if $\Sigma_2$ satisfies Poincar\'e duality.

\noindent {\bf Acknowledgments.} We are grateful to the reviewers for their careful reading of this manuscript and valuable feedback.
The first author thanks the Mathematical Sciences Research Institute, the Simons Institute for the Theory of Computing, the Sorbonne Universit\'e, the Universit\`a di Bolo\-gna, and the Universidad de Los Andes for providing wonderful settings to work on this project, and Laura Escobar, Felipe Rinc\'on, and Kristin Shaw for valuable conversations; his research is supported by NSF grant DMS-1855610 and Simons Fellowship 613384.
The second author thanks the University of Sydney School of Mathematics and
Statistics for hospitality during an early part of this project; his research is supported by NSERC of Canada.
The third author thanks Karim Adiprasito for helpful conversations; his research is
supported by NSF Grant DMS-1638352 and the Ellentuck Fund.

\section{The bipermutohedral fan}\label{sec:fan}

Let $E$ be a finite set of cardinality $n+1$.
For notational convenience, we often identify $E$ with the set of nonnegative integers at most $n$.
As before, we let $\N_E$  be the $n$-dimensional  space
\[
\N_E= \mathbb{R}^E/ \mathbb{R}\mathbf{e}_E, \qquad \mathbf{e}_E =\sum_{i \in E} \mathbf{e}_i.
\]
We write $\N_{E,E}$ for the $2n$-dimensional space $\N_E \oplus \N_{E}$, and 
 $\mu$ for  the addition map
\[
\mu\colon \N_{E,E}  \longrightarrow \N_E, \qquad (z,w) \longmapsto z+w.
\]
Throughout the paper,
all fans in $\N_E$ will be rational with respect to the lattice $\mathbb{Z}^E/\mathbb{Z}  \mathbf{e}_E$,
and all fans in $\N_{E,E}$ will be rational with respect to the lattice $\mathbb{Z}^E/\mathbb{Z}  \mathbf{e}_E\oplus \mathbb{Z}^E/\mathbb{Z}  \mathbf{e}_E$.
We follow \cite{CLS} when using the terms \emph{fan} and \emph{generalized fan}: A generalized fan is a fan if and only if each of its cone is strongly convex.
The notion of morphism of fans is extended to morphism of generalized fans in the obvious way.
For any subset $S$ of $E$, we write $\mathbf{e}_S$ and $\mathbf{f}_S$ for the vectors
\[
\mathbf{e}_S=\sum_{i \in S} \mathbf{e}_i, \qquad \mathbf{f}_S=\sum_{i \in S} \mathbf{f}_i,
\]
where $\mathbf{e}_i$ are the standard basis vectors of $\mathbb{R}^E$ defining the first summand of $\N_{E,E}$
and $\mathbf{f}_i$ are the standard basis vectors of $\mathbb{R}^E$ defining the second summand of $\N_{E,E}$.

In this section, we construct a complete simplicial fan $\Sigma_{E,E}$ in $\N_{E,E}$.
We offer five equivalent descriptions; each one of them will play a role for us. We call it the \emph{bipermutohedral fan} because it is the normal fan of a polytope which we call the \emph{bipermutohedron}. Before we begin defining the bipermutohedral fan $\Sigma_{E,E}$ in $\N_{E,E}$,
we recall some basic facts on the permutohedral fan $\Sigma_E$ in $\N_E$.

\subsection{The normal fan of the simplex}\label{subsec:simplex}

Consider the standard $n$-dimensional simplex
\[
\conv\{\mathbf{e}_i\}_{i \in E} \subseteq \R^E.
\]
 Its normal fan in $\R^E$ has a lineality space spanned by  $\mathbf{e}_E$.
For any convex polytope, we call the quotient of the normal fan   by its lineality space
 the \emph{reduced normal fan} of the polytope.\footnote{
 The  normal fan of a convex polytope in a real vector space  is a generalized fan in the dual space whose face poset is anti-isomorphic to the face poset of the polytope.
 Unlike the reduced normal fan, the normal fan of a polytope need not be a fan. We trust that the use of the term ``normal fan'' will cause no confusion.}
 For example, the reduced normal fan of the standard simplex, denoted $\Gamma_E$, is the complete fan in $\N_E$ with the cones
\[
\sigma_S \coloneq   \cone\{\mathbf{e}_i\}_{i \in S} \subseteq \N_E, \ \  \text{for every proper subset $S$ of $E$.}
\]
The cone $\sigma_S$ consists of the points $z \in \N_E$ such that $\displaystyle \min_{i \in E} z_i = z_s$ for all $s$ not in $S$.
For each element $j$ of $E$,
 the function $\displaystyle \alpha_j= \max_{i \in E} \{z_j-z_i\}$ is piecewise linear on the fan $\Gamma_E$.
 These piecewise linear functions are equal to each other modulo global linear functions on $\N_E$,
and we write $\alpha$ for the common equivalence class of $\alpha_j$.

\subsection{The normal fan of the permutohedron}\label{subsec:permutohedron}

Let $\Pi_E$ be  the  $n$-dimensional permutohedron
\[
\conv\Big\{(x_0,x_1,\ldots,x_n)\ | \ \text{$x_0,x_1,\ldots,x_n$ is a permutation of $0,1,\ldots,n$}\Big\} \subseteq \mathbb{R}^E.
\]
The \emph{permutohedral fan} $\Sigma_E$,  also known as the \emph{braid fan} or the \emph{type $A$ Coxeter complex}, is the reduced normal fan of the permutohedron $\Pi_E$.
It is the complete simplicial fan in $\N_E$ whose chambers are separated by the $n$-dimensional \emph{braid arrangement}, the real hyperplane arrangement in $\N_E$ consisting of the ${n+1 \choose 2}$ hyperplanes
\[
z_i = z_j,  \ \ \text{for distinct elements $i$ and $j$ of $E$.}
\]
The face of the permutohedral fan containing a given point $z$ in its relative interior is determined by the relative order of its homogeneous coordinates $(z_0,\ldots,z_n)$.
Therefore, the faces of the permutohedral fan correspond to the ordered set partitions
\[
\mathcal{P}=(E=P_1\sqcup \cdots \sqcup P_{k+1}),
\]
 which are in bijection with the strictly increasing sequences of nonempty proper subsets
 \[
\mathscr{S}=(\varnothing \subsetneq S_1 \subsetneq \cdots \subsetneq S_k \subsetneq E), \qquad S_m=\bigcup^{m}_{\ell=1} P_\ell.
 \]
The collection of ordered set partitions of $E$ form a poset under \emph{adjacent refinement}, where $\mathcal{P} \leq \mathcal{P}'$ if $\mathcal{P}$ can be obtained from $\mathcal{P}'$ by merging adjacent parts.

\begin{proposition}
The face poset of the permutohedral fan $\Sigma_{E}$ is  isomorphic to the poset of ordered set partitions of $E$.
\end{proposition}

Thus the permutohedral fan has $2(2^n-1)$ rays corresponding to the nonempty proper subsets of $E$
and $(n+1)!$  chambers corresponding to the permutations of $E$.

We now describe the permutohedral fan in terms of its rays.
Two subsets $S$ and $S'$ of $E$ are said to be \emph{comparable} if
\[
S\subseteq S'  \ \ \text{or} \ \  S\supseteq S'.
\]
A \emph{flag} in $E$ is a set of pairwise comparable subsets of $E$.
For any flag $\SS$ of subsets of $E$, we define
\[
\sigma_\mathscr{S}=\cone\{\mathbf{e}_{S}\}_{S \in \mathscr{S}} \subseteq \N_E.
\]
We identify a flag in $E$ with the strictly increasing sequence obtained by ordering the subsets in the flag.

\begin{proposition}\label{prop:permutohedralrays}
The permutohedral fan $\Sigma_E$ is the complete fan in $\N_E$ with the cones
 \[
\sigma_\mathscr{S}=\cone\{\mathbf{e}_{S}\}_{S \in \mathscr{S}}, \ \  \text{where $\mathscr{S}$ is a flag of nonempty proper subsets of $E$.}
 \]
 \end{proposition}

For example, the cone corresponding to the ordered set partition $25|013|4$ is
\[
 \cone(\e_{25}, \e_{01235}) =  \{z \in \N_E \ |\  z_2 = z_5 \geq z_0 = z_1 = z_3 \geq z_4\}.
\]
Proposition \ref{prop:permutohedralrays} shows that the permutohedral fan is a \emph{unimodular fan}: The set of primitive ray generators in any cone in $\Sigma_{E}$ is a subset of a basis of the free abelian group $\mathbb{Z}^E/\mathbb{Z}$.
It also shows that the permutohedral fan is a refinement of the fan $\Gamma_E$ in Section \ref{subsec:simplex}.

\begin{figure}[h]
\[
\begin{tikzpicture}[scale=0.75,baseline=(current bounding box.center),
plain/.style={circle,draw,inner sep=1.2pt,fill=white},
root/.style={circle,draw,inner sep=1.2pt,fill=black}]
\begin{scope}
\node (0) at (-3,0) [root,label=below:$569$] {};
\node (1) at (-1,0) [root,label=below:$7$] {};
\node (2) at (0,0) [root,label=below:$1$] {};
\node (3) at (1,0) [root,label=below:$04$] {};
\node (4) at (2,0) [root,label=below:$28$] {};
\node (5) at (4,0) [root,label=below:$3$] {};
\draw (0) -- (2) -- (5);
\draw (-4,0) -- (-3,0);
\draw (4,0) --  (5,0);
\end{scope}
\end{tikzpicture}
\qquad
\longmapsto
\qquad 3|28|04|1|7|569
\]
\caption{An $E$-tuple of points $p$ 
and its ordered set partition $\pi(p) = 3|28|04|1|7|569$.}\label{fig:configbraid}
\end{figure}
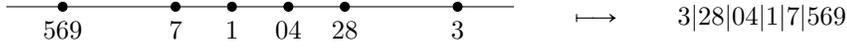

It will be useful to view the permutohedral fan as a configuration space as follows.
Regard $\N_E$ as the space of  $E$-tuples of points $(p_0,\ldots,p_n)$ moving in the real line, modulo simultaneous translation:
\[
p=(p_0, \ldots, p_n) = (p_0+\lambda, \ldots, p_n+\lambda) \ \  \text{for any $\lambda \in \R$.}
\]
The \emph{ordered set partition of $p$}, denoted $\pi(p)$, is obtained by reading the labels of the points in the real line from right to left, as shown in Figure \ref{fig:configbraid}. This model gives the permutohedral fan $\Sigma_E$ the following geometric interpretation.

\begin{proposition}\label{prop:braidorderedsetpartition}
The permutohedral fan $\Sigma_{E}$ is the configuration space of $E$-tuples of points in the real line modulo simultaneous translation, stratified according to their ordered set partition.
\end{proposition}

In Section \ref{ss:configplane}, we give an analogous description of the bipermutohedral fan $\Sigma_{E,E}$ as a configuration space of $E$-tuples of points in the real plane.

\subsection{The bipermutohedral fan as a subdivision} \label{sec:ourfan}

Denote a point in $\N_{E,E}$ by $(z,w)$. We construct the bipermutohedral fan $\Sigma_{E,E}$ in $\N_{E,E}$ as follows.

First, we subdivide $\N_{E,E}$ into the  \emph{charts}  $\mathscr{C}_0,\mathscr{C}_1, \ldots, \mathscr{C}_n$, where $\mathscr{C}_k$ is the cone
\[
\mathscr{C}_k = \Big\{(z, w)   \, | \, \min_{i \in E} (z_i+w_i) = z_k+w_k\Big\}.
\]
These  form the chambers of a complete generalized fan in $\N_{E,E}$, denoted $\Delta_E$.
The chamber $\mathscr{C}_k$ is the inverse image of the cone $\sigma_{E - k}$ under the addition map,
and hence
 $\Delta_E$ is the coarsest complete generalized fan in $\N_{E,E}$ for which the addition map is a morphism to the fan $\Gamma_E$ in Section \ref{subsec:simplex}.
To each chart $\mathscr{C}_k$ we associate  the linear functions
\[
Z_{i} =  z_i - z_k, \quad W_{i}=-w_i+w_k, \ \  \text{for every $i$ in $E$.}
\]
Omitting the zero function $Z_k=W_k$, we obtain a coordinate system $(Z,W)$ for $\N_{E,E}$ such that
\[
\mathscr{C}_k=\Big\{ (Z,W) \ | \ \text{$Z_i \ge W_i$ for every $i$ in $E$}\Big\}.
\]
This coordinate system depends on $k$, but we will drop $k$  from the notation for better readability.

Second, we consider the subdivision $\Sigma_k$ of the cone $\mathscr{C}_k$ obtained from the braid arrangement  of ${2n+1 \choose 2}$ hyperplanes
\[
Z_a = Z_b, \quad W_a = W_b, \quad Z_a=W_b, \quad  \text{for all $a$ and $b$ in $E$.}
\]
Note that the arrangement contains
 the $n$ hyperplanes that cut out $\mathscr{C}_k$ in $\N_{E,E}$.
One may view the subdivision $\Sigma_k$ of $\mathscr{C}_k$ as a copy of $1/2^n$-th of the $2n$-dimensional permutohedral fan, 
namely, the part of the permutohedral fan in $2n+1$ variables $Z_0, W_0, \ldots, Z_k=W_k, \ldots, Z_n, W_n$ where $Z_i \ge W_i$ for every $i \neq k$.
Figure \ref{FigureBipermutohedralFan} illustrates $\Sigma_0$ and $\Sigma_1$ when $n=1$.

\begin{proposition}\label{def:subdivision}
The union of the fans $\Sigma_i$ for $i \in E$ is a fan in $\N_{E,E}$. We call it the \emph{bipermutohedral fan} $\Sigma_{E,E}$.
\end{proposition}

\begin{proof}
To check that $\Sigma_{E,E}$ is indeed a fan,
we  need to check that the fans
$\Sigma_i$ glue compatibly along the boundaries of $\mathscr{C}_i$.
For this, we verify that $\Sigma_i$ and $\Sigma_j$ induce the same subdivision on $\mathscr{C}_i \cap \mathscr{C}_j$ for all $i \neq j$.

Consider the system of linear functions $(Z, W)$ for $\mathscr{C}_i$ and the system of linear functions $(Z', W')$ for $\mathscr{C}_j$.
It is straightforward to check that, for any point in $\N_{E,E}$,  we have
\[
Z_a-Z_b = Z'_a-Z'_b \ \ \text{and} \ \ W_a-W_b = W'_a-W'_b \ \ \text{for all $a$ and $b$ in $E$.}
\]
Furthermore, on the intersection of $\mathscr{C}_i$ and $\mathscr{C}_j$, where  $z_i+w_i=z_j+w_j$,  we have
\[
Z_a-W_b = (z_a-z_i) - (w_i - w_b) = (z_a-z_j) - (w_j - w_b) = Z'_a - W'_b.
\]
Thus the hyperplanes separating the chambers of  $\Sigma_i$ and $\Sigma_j$ have the same intersections with $\mathscr{C}_i \cap \mathscr{C}_j$.
\end{proof}

The following subfan of the bipermutohedral fan will serve as a guide toward Theorem \ref{CSMTheorem}.

\begin{definition}\label{def:cotangent}
The \emph{cotangent fan} $\Omega_{E}$ is the union of the fans $\Sigma_i \cap \Sigma_j$ for $i \neq j \in E$.
\end{definition}

In other words,  $\Omega_E$ is the subfan of $\Sigma_{E,E}$ whose  support  is the tropical hypersurface
\[
\trop(\delta)=\Big\{ (z,w)  \ | \  \text{$\min_{i \in E}(z_i+w_i)$ is achieved at least twice}\Big\} \subseteq \N_{E,E}.
\]
In Section \ref{ss:conormal}, we show that the cotangent fan contains the conormal fan of any matroid on $E$.

\subsection{The bipermutohedral fan as a configuration space}\label{ss:configplane}

It will  be useful to view the bipermutohedral fan $\Sigma_{E,E}$ as a configuration space as follows.
Regard $\N_{E,E}$ as the space of $E$-tuples of points $(p_0,\ldots,p_n)$ moving in the real plane, modulo simultaneous translation:
\[
(p_0, \ldots, p_n) = (p_0+\lambda, \ldots, p_n+\lambda) \ \  \text{for any $\lambda \in \R^2$.}
\]
The point $(z,w)$ in $\N_{E,E}$ corresponds to the points $p_i=(z_i,w_i)$ in $\R^2$ for $i$ in $E$.

\begin{definition}\label{def:Bseq}
A \emph{bisequence} on $E$ is a sequence $\mathscr{B}$ of nonempty subsets of $E$, called the \emph{parts} of $\mathscr{B}$, such that
\begin{enumerate}[(1)]\itemsep 5pt
\item  every element of $E$ appears in at least one part of  $\mathscr{B}$,
\item  every element of $E$ appears in at most two parts of  $\mathscr{B}$, and
\item some element of $E$ appears in exactly one part of $\mathscr{B}$.
\end{enumerate}
The \emph{trivial bisequence} on $E$ is the bisequence with exactly one part $E$.
A \emph{bisubset} of $E$ is a nontrivial bisequence on $E$ of minimal length $2$.
A \emph{bipermutation} of $E$ is a bisequence  on $E$ of maximal length $2n+1$.
\end{definition}

We will write bisequences by listing the elements of its parts, separated by vertical bars.
For example, the bisequence $\{2\},\{0,1\}, \{1\}, \{2\}$ on $\{0,1,2\}$ will be written $2|01|1|2$.

\begin{definition}\label{def:supportingline}
Let $p=(p_0,\ldots,p_n)$ be an $E$-tuple of points in $\mathbb{R}^2$.
\begin{enumerate}[(1)]\itemsep 5pt
\item  The \emph{supporting line} of $p$, denoted $\ell(p)$,  is the lowest  line of slope $-1$ containing a point in $p$.
\item For each point $p_i$, the vertical and horizontal projections of $p_i$ onto $\ell(p)$ will be labelled $i$.
\item The \emph{bisequence of $p$}, denoted $\mathscr{B}(p)$,  is obtained by reading the labels on $\ell(p)$ from right to left.
\end{enumerate}
\end{definition}

\begin{figure}[h]
\[
\begin{tikzpicture}[scale=0.75,baseline=(current bounding box.center),
plain/.style={circle,draw,inner sep=1.2pt,fill=white},
root/.style={circle,draw,inner sep=1.2pt,fill=black}]
\begin{scope}[rotate=-45]
\node (0) at (-1,0) [root,label=below left:$0$] {};
\node (1) at (1,0) [root,label=below left:$24$] {};
\node (2) at (2,0) [plain,label=below left:$1$,label=above right:$p_1$] {};
\node (3) at (3,0) [plain,label=below left:$035$,label=above right:$p_5$] {};
\node (4) at (4.5,0) [root,label=below left:$2$] {};
\node (5) at (6,0) [root,label=below left:$34$, label=below right:$\,\,\,\,\, \, \ell(p)$] {};
\draw (-1.5,0) -- (2) -- (3) -- (7,0);
\node (15) at (3.5,2.5) [plain, label=above right:$p_4$] {};
\node (03) at (1,2) [plain, label=above right:$p_0$] {};
\node (14) at (2.75,1.75) [plain, label=above right:$p_2$] {};
\node (35) at (4.5,1.5) [plain, label=above right:$p_3$] {};
\draw[style=very thin] (0) -- (03) -- (3);
\draw[style=very thin] (1) -- (14) -- (4);
\draw[style=very thin] (3) -- (35) -- (5);
\draw[style=very thin] (14) -- (15) -- (35);
\end{scope}
\end{tikzpicture}
\qquad
\longmapsto
\qquad 34|2|035|1|24|0
\]
\bigskip
\[
\begin{tikzpicture}[scale=0.75,baseline=(current bounding box.center),
plain/.style={circle,draw,inner sep=1.2pt,fill=white},
root/.style={circle,draw,inner sep=1.2pt,fill=black}]
\begin{scope}
\node (0) at (-1,0) [root,label=below:$0$] {};
\node (1) at (1,0) [root,label=below:$24$] {};
\node (2) at (2,0) [root,label=below:$1$] {};
\node (3) at (3,0) [root,label=below:$035$] {};
\node (4) at (4.5,0) [root,label=below:$2$] {};
\node (5) at (6,0) [root,label=below:$34$] {};
\draw (-1.5,0) -- (2) -- (3) -- (7,0);
\end{scope}
\end{tikzpicture} \hspace{5cm}
\]

\caption{An $E$-tuple of points $p=(p_0, \ldots, p_5)$ in the plane, their vertical and horizontal projections onto the supporting line $\ell(p)$, and the bisequence $\mathscr{B}(p)$.}\label{fig:config}
\end{figure}
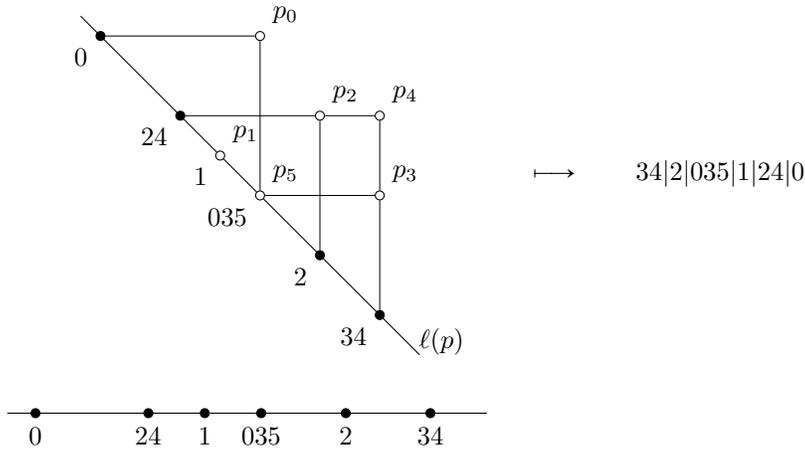

See Figure \ref{fig:config} for an illustration of Definition \ref{def:supportingline}.

\begin{remark}\label{rem:line}
One can recover any configuration $p$ from their projections onto the supporting line $\ell(p)$ and their labels. Therefore, modulo translations, we may also consider $p$ as a configuration of $2n+2$ points on the real line labeled $0,0,1,1,\ldots, n,n$ such that at least one pair of points with the same label coincide. This is illustrated at the bottom of Figure \ref{fig:config}.
\end{remark}

This model gives the bipermutohedral fan $\Sigma_{E,E}$ the following geometric interpretation.

\begin{proposition}\label{prop:configspace}
The bipermutohedral fan $\Sigma_{E,E}$ is the configuration space of $E$-tuples of points in the real plane modulo simultaneous translation, stratified according to their bisequence.
\end{proposition}

\begin{proof}
Consider a point $(z,w)$ in $\N_{E,E}$ and the associated configuration of points $p_i$ in the plane.
The chart $\mathscr{C}_k$ consists of configurations $p$ where $k$ appears exactly once in the bisequence $\mathscr{B}(p)$.
In other words, $p$ is in $\mathscr{C}_k$ if and only if $p_k$ is on the supporting line $\ell(p)$.
  We consider  the system of linear functions $(Z,W)$ for $\mathscr{C}_k$ discussed in Section \ref{sec:ourfan}.
The cones in the subdivision $\Sigma_k$ of $\mathscr{C}_k$ encode the relative order of $Z_0,\ldots,Z_n,W_0,\ldots,W_n$, where
\[
Z_k=W_k=0 \ \ \text{and} \ \ Z_i \ge W_i \ \ \text{for every $i$ in $E$.}
\]

On the other hand, the bisequence $\mathscr{B}(p)$ keeps track of the relative order of the vertical and horizontal projections of $p_i$ onto $\ell(p)$.
As shown in Figure \ref{fig:axisproj},
 after the translation by $(-z_k,-w_k)$, the vertical and horizontal projections of $p_i$ onto $\ell(p)$ are
\[
(z_i,z_k+w_k-z_i)-(z_k,w_k) = (Z_i, -Z_i) \ \ \text{and} \ \
(z_k+w_k-w_i,w_i) -(z_k,w_k) = (W_i, -W_i).
\]
Their relative order along $\ell(p)$ is given by the relative order of $Z_0, \ldots, Z_n, W_0, \ldots, W_n$.
\begin{figure}[h]
\[
\begin{tikzpicture}[scale=0.75,baseline=(current bounding box.center),
plain/.style={circle,draw,inner sep=1.2pt,fill=white},
root/.style={circle,draw,inner sep=1.2pt,fill=black}]
\begin{scope}[rotate=-45]
\node (0) at (-0.4,0) [root,label=below left:{$(W_i, -W_i) = (w_k-w_i,w_i-w_k)$}] {};
\node (2) at (0.8,0) [plain,label=above right:$p_k$] {};
\node (3) at (3,0) [root,label= below left:{$(Z_i, -Z_i) = (z_i-z_k, z_k-z_i)$},label= below right:{$\,\,\,\,\,\,\ell(p)$}] {};
\draw (-0.7, 0) -- (2) -- (4,0);
\node (03) at (1.3,1.7) [plain, label=above right:{$p_i$}] {};
\draw[style=very thin] (0) -- (03) -- (3);
\end{scope}
\end{tikzpicture}
\]
\caption{The vertical and horizontal projections of $p_i$ onto the supporting line $\ell(p)$, after the translation by $(-z_k,-w_k)$.}\label{fig:axisproj}
\end{figure}
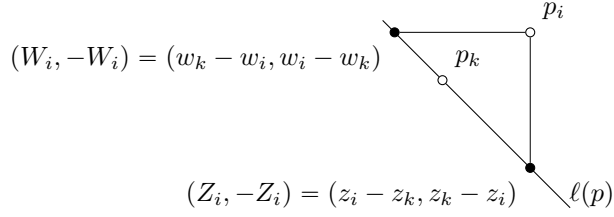
\end{proof}

The collection of bisequences on $E$ form a poset under \emph{adjacent refinement}, where $\mathscr{B} \leq \mathscr{B}'$ if $\mathscr{B}$ can be obtained from $\mathscr{B}'$ by merging adjacent parts.
 The poset of bisequences on $E$ is a graded poset.
 Its $k$-th level  consists of the bisequences of $k+1$ nonempty subsets of $E$,  and the top level consists of the bipermutations of $E$.

\begin{proposition}\label{p:fanBseqs}
The face poset of the  bipermutohedral fan $\Sigma_{E,E}$ is isomorphic to the poset of bisequences on $E$.
\end{proposition}

\begin{proof}
Remark \ref{rem:line} shows that, given any bisequence $\mathscr{B}$ on $E$, there is a configuration $p$ with $\mathscr{B}(p)=\mathscr{B}$.
Thus, by Proposition \ref{prop:configspace}, the cones in $\Sigma_{E,E}$ are in bijection with the bisequences on $E$.
If a configuration $p$ moves into more special position, then some adjacent parts of $\mathscr{B}(p)$ merge.
\end{proof}

\begin{figure}[h]
\centering
\begin{tikzpicture}[root/.style={circle,draw,inner sep=1.2pt,fill=black},every node/.style={scale=0.95}]
\draw[style=very thick,color=gray] (2,0) -- (-2,0);
\draw[style=thin,color=gray] (-1,-1.72) -- (1,1.72);
\draw[style=thin,color=gray] (1,-1.72) -- (-1,1.72);
\node at (0,-1.5) {$1|0|1$};
\node at (0,1.5) {$0|1|0$};
\node at (1.3,-0.7) {$1|1|0$};
\node at (-1.3,-0.7) {$0|1|1$};
\node at (1.3,0.7) {$1|0|0$};
\node at (-1.3,0.7) {$0|0|1$};
\node at (0,0) [root, label={[label distance=6pt]below:{$\,01$}}] {};
\node at (.8,-1.45) [label={[label distance=5pt]300:{$1|01$}}] {};
\node at (.8,1.45) [label={[label distance=5pt]60:{$01|0$}}] {};
\node at (-.8,-1.45) [label={[label distance=5pt]240:{$01|1$}}] {};
\node at (-.8,1.45) [label={[label distance=5pt]120:{$0|01$}}] {};
\node at (1.9,0) [label=right:{$1|0$}] {};
\node at (-1.9,0) [label=left:{$0|1$}] {};
\node (1) at (3,0) {};
\node (2) at (4.5,0) {};
\draw[->] (1) to node[above]{$\mu$} (2);
\node at (2.75,1) {$\Sigma_1$};
\node at (2.75,-1) {$\Sigma_0$};
\draw[color=gray] (5.5,-1.72) -- (5.5,1.72);
\node at (5.5,0) [root, label=right:{$01$}] {};
\node at (5.5,-0.7) [label=left:{$0|1$}] {};
\node at (5.5,0.7) [label=left:{$1|0$}] {};
\end{tikzpicture}
\caption{The map
$\mu\colon \Sigma_{\{0,1\},\{0,1\}}\to \Sigma_{\set{0,1}}$ from the bipermutohedral fan to the permutohedral fan, and the labelling of their cones with bisequences on $\{0,1\}$ and ordered set partitions on $\{0,1\}$.\label{fig:F2}}\label{FigureBipermutohedralFan}
\end{figure}
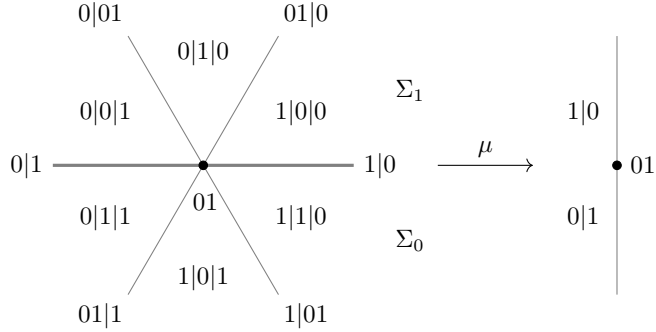

For a bisequence $\mathscr{B}$ on $E$, we write $\sigma_\mathscr{B}$ for the corresponding  cone defined by
\[
\sigma_\mathscr{B}=\closure\Big\{ \text{configurations $p$ satisfying $\mathscr{B}(p)=\mathscr{B}$}\Big\} \subseteq \N_{E,E}.
\]
In terms of the cones $\sigma_{\BB}$, the fan $\Sigma_i$ subdividing the chart $\C_i$  can be described as the subfan
\[
\Sigma_i=\{\sigma_{\BB} \, | \, \text{$i$ appears exactly once in the bisequence $\BB$}\} \subseteq \Sigma_{E,E}.
\]
See  Figure \ref{fig:F2} for an illustration of Proposition \ref{p:fanBseqs} when $n=1$.

\subsection{The bipermutohedral fan as a common refinement}\label{ss:resolves}

The importance of the bipermutohedral fan $\Sigma_{E,E}$ stems from its relationship with the normal fan $\Gamma_E$ of the standard simplex and the permutohedral fan $\Sigma_E$ described in Sections in  Section \ref{subsec:simplex}  and \ref{subsec:permutohedron}. Recall that a \emph{morphism} from a fan $\Sigma_1$ in $\N_1$ to a fan $\Sigma_2$ in $\N_2$ is an integral linear map from $\N_1$ to $\N_2$ that maps any cone in $\Sigma_1$ into a cone in $\Sigma_2$.

\begin{proposition}\label{prop:mu}
The bipermutohedral fan $\Sigma_{E, E}$
has the following properties.
\begin{enumerate}[(1)]\itemsep 5pt
\item
The projections $\pi(z,w) = z$ and $\overline \pi(z,w) = w$ are morphisms of fans from $\Sigma_{E, E}$ to  $\Sigma_E$.
\item The addition map $\mu(z,w) = z+w$ is a morphism of fans from $\Sigma_{E, E}$ to  $\Gamma_E$.
\end{enumerate}
\end{proposition}

\begin{proof}
That $\Sigma_{E, E}$ has the stated properties follows from the interpretation of $\Sigma_E$ and $\Sigma_{E, E}$ as  configuration spaces, as we now explain.
Suppose $(z,w)$ is a point in $\N_{E,E}$ and $p$ is the corresponding $E$-tuple of points in $\R^2$ modulo simultaneous translation, with corresponding bisequence $\mathscr{B}(p)$. Then
the smallest cone of $\Gamma_E$ containing $z+w$ is given by the entries that appear twice in $\mathscr{B}(p)$. The ordered set partition of $z$ in $\N_E$
is given by the first occurrence   of each $i$ in $\mathscr{B}(p)$.
Similarly, the ordered set partition of $w$ in $\N_E$ is given by the order of the last occurrence of each $i$ in $\mathscr{B}(p)$.
For example, if a point $(z,w)$ has the bisequence $34|2|035|1|24|0$, as in Figure \ref{fig:config},
then the sum $z+w$ is in the cone of ${0234}$ in $\Gamma_E$,
the first projection
$z$ is in the cone of ${34|2|05|1}$ in $\Sigma_E$, and the second projection
$w$ is in the cone of ${0|24|1|35}$ in $\Sigma_{E}$.
\end{proof}

We note, however, that the bipermutohedral fan is not the coarsest fan
structure for which projections and addition are morphisms: the reader
is referred to the discussion in Section~\ref{sec:whybiperm}.

\subsection{The bipermutohedral fan in terms of its rays and cones}\label{sec:combinfan}

The rays of the bipermutohedral fan $\Sigma_{E,E}$ correspond to the bisubsets of $E$.
In other words, the rays of $\Sigma_{E,E}$ correspond to the ordered pairs of nonempty subsets $S|T$ of $E$ such that
\[
S\cup T=E\ \ \text{and} \ \ S \cap T \neq E.
\]

\begin{proposition}\label{prop:rays}
The $3(3^n-1)$ rays of the bipermutohedral fan $\Sigma_{E,E}$ are generated by
\[
\mathbf{e}_{S|T}\coloneq  \e_S+\f_T, \ \ \text{where $S|T$ is a bisubset of $E$.}
\]
\end{proposition}

\begin{proof}
The configuration $p$ corresponding to $\e_{S|T}$ has  points with labels in $S \cap T$ located at $(1,1)$, the points with labels in $S - T$ located at $(1,0)$, and the points with labels in $T - S$ located at $(0,1)$. The  bisequence of $p$ is indeed $S|T$,
and hence the conclusion follows from Proposition \ref{prop:configspace}.
\end{proof}

\begin{proposition}\label{prop:chambers}
The bipermutohedral fan $\Sigma_{E,E}$ has $ (2n+2)!/2^{n+1}$  chambers.
\end{proposition}

\begin{proof}
By Proposition \ref{p:fanBseqs}, the chambers correspond to the bipermutations. These are obtained bijectively from the $(2n+2)!/2^{n+1}$ permutations of the multiset $\{0,0, \ldots, n,n\}$ by dropping the last letter in the one-line notation for permutations.
For example, the bipermutation $1|0|1|2|3|0|3$ correspond to the permutation $10123032$ of $\{0,0, 1,1,2,2, 3,3\}$.
\end{proof}

It is worth understanding Proposition \ref{prop:chambers} in a different way.
Recall that the bipermutohedral fan is obtained by gluing copies of $1/2^n$-th of the $2n$-dimensional permutohedral fan.
There are $(n+1)$ such copies, and each copy contains $(2n+1)!/2^n$ chambers,
producing the total of $ (2n+2)!/2^{n+1}$ chambers.
This viewpoint explains why Figure \ref{fig:F2} deceivingly looks like a permutohedral fan:
 For $n=1$, the bipermutohedral fan consists of two glued copies of half of the permutohedral fan.

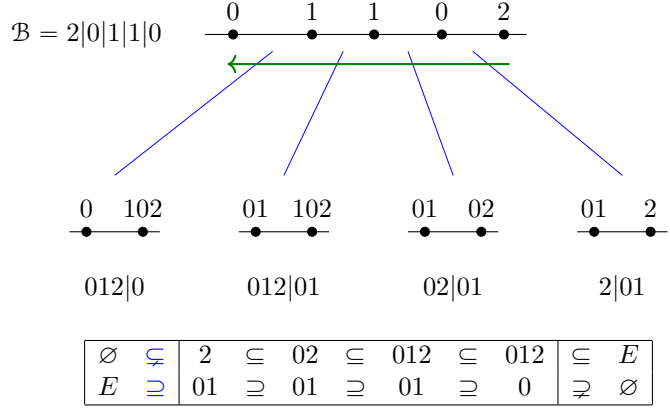
\begin{figure}[th]
\def\dy{3.5}
\begin{tikzpicture}[scale=0.75,baseline=(current bounding box.center),
root/.style={circle,draw,inner sep=1.2pt,fill=black}]
\node (0) at (-0.4,0) [root,label=above:$0$] {};
\node (1) at (1,0) [root,label=above:$1$] {};
\node (2) at (2.1,0) [root,label=above:$1$] {};
\node (3) at (3.3,0) [root,label=above:$0$] {};
\node (4) at (4.4,0) [root,label=above:$2$] {};
\draw (-.9,0) -- (4.8,0);
\node at (-3,0) {$\BB=2|0|1|1|0$};

\node (1T) at (-3,-\dy) [root,label=above:$0$] {};
\node (1S) at (-2,-\dy) [root,label=above:$102$] {};
\draw (-3.3,-\dy) -- (-1.7,-\dy);
\node at (-2.5,-\dy-1) {$012|0$};

\node (2T) at (0,-\dy) [root,label=above:$01$] {};
\node (2S) at (1,-\dy) [root,label=above:$102$] {};
\draw (-0.3,-\dy) -- (1.3,-\dy);
\node at (0.5,-\dy-1) {$012|01$};

\node (3T) at (3,-\dy) [root,label=above:$01$] {};
\node (3S) at (4,-\dy) [root,label=above:$02$] {};
\draw (2.7,-\dy) -- (4.3,-\dy);
\node at (3.5,-\dy-1) {$02|01$};

\node (4T) at (6,-\dy) [root,label=above:$01$] {};
\node (4S) at (7,-\dy) [root,label=above:$2$] {};
\draw (5.7,-\dy) -- (7.3,-\dy);
\node at (6.5,-\dy-1) {$2|01$};

\draw[blue,very thin] (0.3,-0.3) -- (-2.5,1-\dy);
\draw[blue,very thin] (1.55,-0.3) -- (0.5,1-\dy);
\draw[blue,very thin] (2.7,-0.3) -- (3.5,1-\dy);
\draw[blue,very thin] (3.85,-0.3) -- (6.5,1-\dy);
\draw[darkgreen,thick,->] (4.5,-0.15*\dy) -- (-0.5,-0.15*\dy) node[midway,above] {};

\node at (2,-6) {
$
\begin{array}{|cc|ccccccc|cc|}
\hline
\varnothing& \blue{\subsetneq} & 2 & \subseteq & 02 & \subseteq & 012& \subseteq & 012  &\subseteq&E\\
E & \blue{\supseteq} &01& \supseteq & 01 & \supseteq & 01 & \supseteq & 0  &\supsetneq&\varnothing\\
\hline
\end{array}
$
};
\end{tikzpicture}
\caption{The cone of $2|0|1|1|0$ has the rays generated by $\mathbf{e}_{2|01}, \mathbf{e}_{02|01}, \mathbf{e}_{012|01}, \mathbf{e}_{012|0}$.\label{fig:rays}
}
\end{figure}

We now describe the cones in the bipermutohedral fan in terms of their generating rays.
Let $\mathscr{B}=B_0| B_1|\cdots|B_k$ be a bisequence on $E$.
Propositions \ref{p:fanBseqs} and \ref{prop:rays} show that the rays of the $k$-dimensional cone $\sigma_\mathscr{B}$
are generated by the vectors
\[
\mathbf{e}_{S_1|T_1},\ldots,\mathbf{e}_{S_k|T_k}, \ \ \text{where} \ \  S_i=\bigcup_{j=0}^{i-1} B_j \ \ \text{and} \ \  T_i=\bigcup_{j=i}^k B_j.
\]
See Figure \ref{fig:rays} for an illustration.
We use the following table to record the rays of $\sigma_\mathscr{B}$:
\[
\begin{array}{|cc|ccccccc|cc|}
\hline
\varnothing & \subsetneq & S_1 & \subseteq & S_2 & \subseteq & \cdots & \subseteq & S_k  & \subseteq & E\\
E & \supseteq & T_1 & \supseteq & T_2 & \supseteq & \cdots & \supseteq & T_k & \supsetneq & \varnothing\\
\hline
\end{array}
\]
For each index $j$ such that
$S_j \subsetneq S_{j+1}$ and $T_j \supsetneq T_{j+1}$, we mark those two strict inclusions in blue.
We write $\SS(\mathscr{B})|\T(\mathscr{B})$  for the collection of bisubsets $S_i|T_i$
constructed from $\mathscr{B}$ as above by merging adjacent parts.
For convenience, we also refer to the pairs $S_0|T_0=\varnothing|E$ and  $S_{k+1}|T_{k+1}=E|\varnothing$.

Conversely, we may ask which subsets of $k$ rays in $\Sigma_{E,E}$ generate a $k$-dimensional cone in $\Sigma_{E,E}$.
To answer this question, we introduce the notion of a flag of bisubsets.

\begin{definition}\label{biflagProposition}
We say that two bisubsets $S|T$ and $S'|T'$ of $E$ are \emph{comparable} if
\[
\text{($S \subseteq S'$ and $T \supseteq T'$)} \ \  \text{or} \ \  \text{($S \supseteq S'$ and $T \subseteq T'$)}.
\]
A \emph{flag of bisubsets} in $E$, or a \emph{biflag} in $E$,  is a set $\mathscr{S}|\mathscr{T}$ of  pairwise comparable bisubsets of $E$ satisfying
\[
\bigcup_{S|T \in \mathscr{S}|\mathscr{T}} S \cap T \neq E.
\]
The \emph{length} of a biflag  is the number of bisubsets in it.
 \end{definition}

We have the following useful alternative characterization of biflags in $E$.

\begin{proposition}\label{prop:alt}
Let $\mathscr{S}$ be an increasing sequence of $k$ nonempty subsets of $E$, say
\[
\mathscr{S}=(\varnothing \subsetneq S_1 \subseteq \cdots \subseteq S_k \subseteq E),
\]
and let $\T$ be a decreasing sequence of $k$ nonempty subsets of $E$, say
\[
\mathscr{T} = (E \supseteq T_1 \supseteq \cdots \supseteq T_k \supsetneq \varnothing).
\]
Then the set $\mathscr{S}|\mathscr{T}$ consisting of the pairs $S_1|T_1,\ldots,S_k|T_k$ is a flag of bisubsets if and only if
 \[
 S_j \cup T_{j} = E  \ \ \text{for every  $1 \leq j \leq k$ and} \ \   S_j \cup T_{j+1} \neq E \ \ \text{for some  $0 \leq j \leq k$}.
 \]

 \end{proposition}

\begin{proof}
If $\mathscr{S}|\mathscr{T}$  is a biflag in $E$, then each $S_j|T_j$ is a bisubset of $E$, and hence $ S_j \cup T_{j} = E$ for all $j$.
Now let $e$ be an element not in the union of all $S_j \cap T_j$, and consider the largest index $i$ for which $e \notin S_i$.
Then $e \in S_{i+1}$, which implies $e \notin T_{i+1}$ by the definition of $e$.
Therefore,  $S_i \cup T_{i+1} \neq E$.

Conversely,  if $\SS$ and $\T$ satisfy the stated conditions, then the pairs  $S_j|T_j$ form a set of pairwise comparable bisubsets of $E$.
If $e$ is an element not in $S_j \cup T_{j+1}$ for some index $j$,
then $e$ is not in $S_k$ for all indices $k \le j$ and
$e$ is not in $T_k$ for all indices $k>j$.
Therefore, $e$ is not in the union of all $S_k \cap T_k$, as desired.
\end{proof}

Since $S_j \cup T_{j+1} \neq E$ implies $S_j \subsetneq S_{j+1}$ and $T_j \supsetneq T_{j+1}$, 
Proposition \ref{prop:alt} shows that the table of any biflag has at least one pair of strict inclusions marked in blue. 

For  a biflag $\SS|\T$ of length $k$, we write $\SS$ for the increasing sequence of $k$ nonempty subsets
\[
\mathscr{S}=(\varnothing \subsetneq S_1 \subseteq \cdots \subseteq S_k \subseteq E), \ \  \text{where $S_j$ are the first parts of the bisubsets in $\SS|\T$},
\]
and write $\T$ for the decreasing sequence of $k$ nonempty subsets
\[
\mathscr{T} = (E \supseteq T_1 \supseteq \cdots \supseteq T_k \supsetneq \varnothing), \ \  \text{where $T_j$ are the second parts of the bisubsets in $\SS|\T$}.
\]
We use $\SS$ and $\T$ to define $\BB({\SS|\T})$ as the sequence of $k+1$ nonempty sets
 \[
 B_0|B_1|\cdots|B_k, \ \ \text{where $B_j=(S_{j+1} - S_j) \cup (T_j - T_{j+1})$}.
 \]
The above construction gives an isomorphism between the poset of bisequences under adjacent refinement and the poset of biflags under inclusion.

\begin{proposition}\label{prop:cones}
The bisequences on  $E$ are in bijection with the biflags in $E$.
More precisely,
\begin{enumerate}[(1)]\itemsep 5pt
\item if $\BB$ is a bisequence on $E$, then $\SS(\mathscr{B})|\T({\BB})$ is a biflag in $E$,
\item if $\SS|\T$ is a biflag in $E$, then $\BB(\SS|\T)$ is a bisequence on $E$, and
\item the constructions $\SS(\mathscr{B})|\T({\BB})$ and $\BB(\SS|\T)$ are inverses to each other.
\end{enumerate}
\end{proposition}

Note that a bisubset $S|T$ corresponds to the biflag $\{S|T\}$ under the above bijection.
For  simplicity, we use the two symbols interchangeably.

\begin{proof}
Let $\mathscr{B}$ be a bisequence on $E$.
 Since every element of $E$ appears at least once in  $\mathscr{B}$, the increasing flag $\SS(\BB)$ and the decreasing flag $\T(\BB)$ satisfy  $S_j \cup T_j = E$ for all $j$.
In addition,  since some element of $E$ appears exactly once in $\BB$, say in $B_j$, we have $S_j \cup  T_{j+1} \neq E$ for some $j$.
 Therefore, by Proposition \ref{prop:alt},
 the pair $\SS(\mathscr{B})|\T(\BB)$ is a biflag in $E$.

Conversely, let $\SS|\T$ be a biflag in $E$.
Since $S_1|T_1,\ldots,S_k|T_k$ are pairwise distinct, $B_j$ must be nonempty for all $j$.
Clearly, every element in $E$ must appear in $B_j$ for some $j$.
In addition, each element $e$ in $E$ can occur at most twice in $\BB(\SS|\T)$, namely, in the parts $B_a$ and $B_{b}$ whose indices satisfy $e \in S_{a+1}- S_a$ and $e \in T_{b}- T_{b+1}$. Furthermore, by Proposition \ref{prop:alt},
there is an element $e$ not in $S_c \cup T_{c+1}$ for some index $c$,
 and in this case we must have $a = b = c$. That element $e$ can occur only in the part $B_a$ of $\BB(\SS|\T)$, and hence $\BB(\SS|\T)$ indeed is a bisequence.

It is straightforward to check that the constructions $\SS(\mathscr{B})|\T({\BB})$ and $\BB(\SS|\T)$ are inverses to each other.
\end{proof}

We identify a biflag $\SS|\T$ in $E$ with the sequence of bisubsets of $E$ obtained by ordering the bisubsets in  $\SS|\T$ as above.
For any sequence $\SS|\T$ of bisubsets of $E$, we define
\[
\sigma_{\mathscr{S}|\T}=\cone\{\mathbf{e}_{S|T}\}_{S|T \in \mathscr{S}|\T} \subseteq \N_{E,E}.
\]
Thus, for any bisequence $\BB$ on $E$, we have
$
\sigma_{\BB}=\sigma_{\SS(\BB)|\T(\BB)}$.

\begin{corollary}\label{cor:bipermutohedralrays}
The bipermutohedral fan $\Sigma_{E,E}$ is the complete fan in $\N_{E,E}$ with the cones
 \[
\sigma_{\mathscr{S}|\T}=\cone\{\mathbf{e}_{S|T}\}_{S|T \in \mathscr{S}|\T}, \ \  \text{for flags of bisubsets $\mathscr{S}|\T$ of $E$.}
 \]
 \end{corollary}

 \begin{proof}
 The statement is straightforward, given Propositions \ref{p:fanBseqs} and \ref{prop:cones}.
 \end{proof}

Corollary \ref{cor:bipermutohedralrays} can be used to show that the bipermutohedral fan is a unimodular fan.\footnote{Alternatively, one may appeal to the unimodularity of the $2n$-dimensional braid arrangement fan in $(Z,W)$-coordinates discussed in Section \ref{prop:braidorderedsetpartition}.}

\begin{proposition}\label{prop:unimodular}
The set of primitive ray generators of any chamber of $\Sigma_{E,E}$ is a basis of the free abelian group $\mathbb{Z}^E/\mathbb{Z}  \mathbf{e}_E \oplus \mathbb{Z}^E/\mathbb{Z}  \mathbf{f}_E$. In particular, $\Sigma_{E,E}$ is simplicial.
\end{proposition}

\begin{proof}
Let $\SS=\SS(\mathscr{B})$ and $\T=\T(\mathscr{B})$ for a bipermutation $\mathscr{B}$ of $E$.
If $0$ is the unique element of $E$ that appears exactly once in $\mathscr{B}$, then
\[
\Big\{\mathbf{e}_{S_{j+1}|T_{j+1}} -\mathbf{e}_{S_{j}|T_{j}} \ | \ \text{$0$ is contained in $S_j \cup T_{j+1}$}\Big\}=\Big\{\e_1,\ldots,\e_n,\f_1,\ldots,\f_n\Big\}.
\]
Therefore, the set of $2n$ primitive ray generators of $\sigma_{\mathscr{B}}$ generates $\mathbb{Z}^E/\mathbb{Z}  \mathbf{e}_E \oplus \mathbb{Z}^E/\mathbb{Z}  \mathbf{f}_E$.
\end{proof}

\subsection{The bipermutohedral fan as the normal fan of the bipermutohedron}\label{sec:polytopality}

We construct a polytope $\Pi_{E,E}$, called the \emph{bipermutohedron},
whose reduced normal fan is $\Sigma_{E,E}$.  The reader may skip
this subsection without interrupting the main logical flow of the paper (see
Remark~\ref{rem:shortcut}).  

For each bipermutation $\mathscr{B}$ of $E$, we construct a vertex $\mathbf{v}_\mathscr{B}$ in $\R^E \oplus \R^E$ as follows.
Let $k=k_\mathscr{B}$ be the element appearing exactly once in $\mathscr{B}$.
Consider the word $\pi_\BB$ obtained from $\mathscr{B}$ by replacing $k$ with $k \overline{k}$
and replacing the first and the second occurrences of each $j \neq k$ with $j$ and $\overline{j}$.
We identify this word with the bijection 
\[
\pi_\mathscr{B} \colon
E \cup \E
\longrightarrow
\{-2n-1,  \ldots, -3, -1, 1, 3, \ldots, 2n+1\}
\]
that sends the letters of the word to the odd integers in increasing order. 
For example, 
\[
\pi_{1|2|3|1|3|0|0} =
\left(
\begin{array}{rrrrrrrr}
1 & 2 & \2 & 3 & \1 & \3 & 0 & \0 \\
-7 & -5 & -3 & -1 & \,\,\,1 & \,\,\,3 & \,\,\,5 & \,\,\,7 \\
\end{array} \right).
\]
Let $\mathbf{u}_{\mathscr{B}} =(x,y)$ be the vector in  $\R^E \oplus \R^E$ with 
  $x_j = \pi_\BB(j)$ and $y_j = -\pi_\BB(\overline{j})$.
%
We define
\[
\mathbf{v}_\BB = \mathbf{u}_{\BB} - s_{\BB}(\mathbf{e}_k + \mathbf{f}_k), \ \ 
\text{where}
\ \ 
s_{\mathscr{B}} = \sum_{i \in E} x_i = \sum_{i \in E} y_i.
\]
For example,
writing $(x,y)$ as a matrix 
whose top and bottom rows are $x$ and $y$ respectively, 
\[
\mathbf{v}_{1|2|3|1|3|0|0} 
=
\left(
\begin{array}{rrrr}
5 & -7 & -5 & -1 \\
-7 & -1 & 3 & -3 \\
\end{array}
\right)
+ 8 \, \,
\left(
\begin{array}{rrrr}
0 & 0 & 1 & 0 \\
0 & 0 & 1 & 0 \\
\end{array}
\right).
\]
The row sums of $\mathbf{v}_\BB$ are both equal to $0$, so $\mathbf{v}_\BB$ is in $\M_E \oplus \M_E$, where $\M_E$ is the vector space dual to $\N_E$.

\begin{definition}
The \emph{bipermutohedron} of $E$ is
\[
\Pi_{E,E} \coloneq  \conv \{\mathbf{v}_\BB \, | \, \text{$\BB$ is a bipermutation of $E$} \} \subseteq  \M_E \oplus \M_E.
\]
\end{definition}

We refer to \cite{bipermutohedron} for a detailed study of this remarkable polytope.
In \cite[Proposition 2.8]{bipermutohedron}, it is shown that the bipermutohedron in $\M_E \oplus \M_E$ is defined by the inequalities
\[
\sum_{i \in S} x_i +\sum_{i \in T} y_i \ge - (|S|+|S-T|) \cdot  (|T|+|T-S|), \ \ \text{for each bisubset $S|T$ of $E$.}
\]
The inequality description is reminiscent of that of the permutohedron in $\M_E$, which reads 
\[
\sum_{i \in S} x_i \ge - |S| \cdot |E-S|, \ \ \text{for each nonempty proper subset $S$ of $E$.}
\]
The automorphism group of the permutohedron is the symmetric group $\mathfrak{S}_E$,
and the automorphism group of the bipermutohedron is the product $\mathfrak{S}_E \times \mathbb{Z}/2\mathbb{Z}$ \cite[Proposition 7.2]{bipermutohedron}.

\begin{proposition}\label{thm:ardilahedron}
The bipermutohedral fan $\Sigma_{E,E}$ is the normal fan of $\Pi_{E,E}$.
\end{proposition}

\begin{proof}
Let $\BB$ be a bipermutation on $E$. We claim that the cone of the normal fan of $\Pi_{E,E}$ corresponding to $\mathbf{v}_\BB$ is precisely the maximal cone $\sigma_\BB$ of the bipermutohedral fan $\Sigma_{E,E}$:
\[
\mathcal{N}_{\Pi_{E,E}}(\mathbf{v}_\BB) = \sigma_\BB.
\]
This will also show that each $\mathbf{v}_\BB$ is indeed a vertex of the bipermutohedron $\Pi_{E,E}$.

It is enough to show that the left-hand side is contained in the right-hand side, as the two fans we are comparing have the same support.
Let $\varphi=(z,w)$ be a linear functional such  that the $\varphi$-minimal face of $\Pi_{E,E}$ contains  $\mathbf{v}_\BB$, and let
 $k$ be the letter that is not repeated in $\BB$.
 We use the description of $\Sigma_{E,E}$ in Section \ref{sec:ourfan}.
We need to show that $(z,w)$ is in the chart $\C_k$, 
and that when we rewrite $(z,w)$ in the coordinate system 
\[
Z_i = z_i-z_k, \qquad W_i = -w_i+w_k, 
\]
the relative order of $Z_0, \ldots, Z_n, W_0, \ldots, W_n$ agrees with the order of $0, \ldots, n, \0, \ldots, \n$ in   $\pi_\BB$.

Let $i$ and $j$ be  any two adjacent letters in $\BB$ appearing in that order. 
When $i\neq j$, we write
$\BB'$ for the bipermutation obtained from $\BB$ by swapping $i$ and $j$:
\[
\BB = \cdots | i | j |  \cdots, \qquad \BB' = \cdots | j | i |  \cdots.
\]
When $i=j$, we write $\BB'$ for the bipermutation obtained from $\B$ by making $i$ occur only once and $k$ occur twice consecutively, as follows:
\[
\BB = \cdots | i | i |  \cdots | k | \cdots, \qquad \BB' = \cdots | i |  \cdots | k | k | \cdots. 
\]
We use  the inequality $\varphi(\mathbf{v}_{\BB}) \leq \varphi(\mathbf{v}_{\BB'})$ to 
determine  the relative order of $Z_0,\ldots,Z_n, W_0,\ldots,W_n$.
In what follows,  we consider $\mathbf{v}_\mathscr{B}$ and  $\mathbf{v}_{\mathscr{B}'}$ as  matrices with two rows whose columns are labeled by $E$.
Since  $\mathbf{v}_\BB$ and $\mathbf{v}_{\BB'}$ can only differ in columns labeled by $i$, $j$, or $k$,
we only display those columns in the computations below. 



First, we consider the case $i \neq j$ and $i , j \neq k$. 
There are four subcases.

\noindent
(1-1) Suppose both $i$ and $j$ are their first occurrences in $\mathscr{B}$. In this case, we have   $s_{\BB} = s_{\BB'}$,
and 
\[
\pi_\BB(i)=a-1,  \ \  \pi_\BB(j) = a+1, \ \   \pi_{\BB'}(i)=a+1, \ \   \pi_{\BB'}(j) = a-1 \ \ \text{for some  $a$.} 
\]
 Therefore, the condition that  the $\varphi$-minimal face of $\Pi_{E,E}$ contains  $\mathbf{v}_\BB$ implies
\[
\varphi
\left(
\begin{array}{ccc}
a-1&a+1&-s\\
b&c&-s
\end{array}
\right)
\le
\varphi
\left(
\begin{array}{ccc}
a+1&a-1&-s\\
b&c&-s
\end{array}
\right)
\ \ 
\text{for some $b$ and $c$.}
\]
We thus have
$(a-1)z_i + (a+1)z_j \leq (a+1)z_i + (a-1)z_j$,
and hence
 $Z_j \leq Z_i$.

\noindent
(1-2) Suppose both $i$ and  $j$ are their second occurrences in $\mathscr{B}$.
Similarly to the previous case, 
\[
\varphi
\left(
\begin{array}{ccc}
b &c &-s\\
-a+1&-a-1&-s
\end{array}
\right)
\le
\varphi
\left(
\begin{array}{ccc}
b&c&-s\\
-a-1&-a+1&-s
\end{array}
\right)
\ \ 
\text{for some $b$ and $c$.}
\]
We thus have $-(a-1)w_i - (a+1)w_j \leq -(a+1)w_i - (a-1)w_j$, and hence $W_j \le W_i$.


\noindent
(1-3) Suppose $i$ is its first occurrence in $\mathscr{B}$ and $j$ is its second occurrences in $\mathscr{B}$.
We have
\[
\pi_\BB(i)=a-1,  \ \  \pi_\BB(\overline{j}) = a+1, \ \   \pi_{\BB'}(i)=a+1, \ \   \pi_{\BB'}(\overline{j}) = a-1 \ \ \text{for some $a$,} 
\]
and hence   $s_{\BB'} = s_{\BB}+2$. 
The condition that  the $\varphi$-minimal face of $\Pi_{E,E}$ contains  $\mathbf{v}_\BB$ implies
\[
\varphi
\left(
\begin{array}{ccc}
a-1 &b &-s\\
c &-a-1&-s
\end{array}
\right)
\le
\varphi
\left(
\begin{array}{ccc}
a+1&b&-s-2\\
c&-a+1&-s-2
\end{array}
\right)
\ \ 
\text{for some $b$ and $c$.}
\]
We thus have
$(a-1)z_i - (a+1)w_j \leq (a+1)z_i - (a-1)w_j - 2z_k - 2w_k$,
and hence $W_j \le Z_i$.


\noindent
(1-4) Suppose $i$ is its second occurrence in $\mathscr{B}$ and $j$ is its first occurrences in $\mathscr{B}$.
Computing as in the previous case, we get  $Z_j \leq W_i$.

Second, we consider the case  $i \neq j$ and $i=k$.
There are two subcases.

\noindent
(2-1) Suppose $j$ is its first occurrence in $\mathscr{B}$.
In this case, for some  $a$, we have
\[
\pi_\BB(k)=a-2, \ \  \pi_{\BB}(\k)=a, \ \  \pi_\BB(j) = a+2, \ \  \pi_{\BB'}(k)=a, \ \  \pi_{\BB'}(\k)=a+2,\ \  \pi_{\BB'}(j) = a-2.
\]
and hence $s_{\BB'}=s_{\BB}-2$.
The condition that  the $\varphi$-minimal face of $\Pi_{E,E}$ contains  $\mathbf{v}_\BB$ implies
\[
\varphi
\left(
\begin{array}{cc}
a-s-2&a+2\\
-a-s&b
\end{array}
\right)
\le
\varphi
\left(
\begin{array}{cc}
a-s+2&a-2\\
-a-s&b
\end{array}
\right)
\ \ 
\text{for some $b$.}\]
We thus have $(a-s-2)z_i  +  (a+2)z_j \leq (a-s+2)z_i  +(a-2)z_j$,
and hence
$Z_j \le Z_i$.


\noindent
(2-2) Suppose $j$ is its second occurrence in $\mathscr{B}$.
Computing as above, we get  $W_j \leq W_i$.

Third, we consider the case  $i \neq j$ and $j =k$.
There are two subcases.

\noindent
(3-1) Suppose $i$ is its first occurrence in $\mathscr{B}$.
Computing as in (2-1), we get $Z_j \leq Z_i$.

\noindent
(3-2) Suppose $i$ is its second occurrence in $\mathscr{B}$.
Computing as in (2-1), we get $W_j \leq W_i$.

Last, we consider the case  $i=j$.
In this case, we have $\pi_\BB = \pi_{\BB'}$, and hence
\[
 \mathbf{u}_\BB = \mathbf{u}_{\BB'} \ \ \text{and}  \ \ s_\BB = s_{\BB'}.
 \]
Notice that, since $\ell$ precedes $\overline{\ell}$ in $\pi_\BB$ for all $\ell$, we have
\[
s_\BB \leq -(2n-1) - (2n-5) - \cdots + (2n-7) + (2n-3) < 0.
\]
Therefore, 
$
\varphi \big(\mathbf{u}_\BB - s_\BB(\e_k+\f_k) \big) \leq \varphi \big(\mathbf{u}_\BB - s_\BB(\e_i+\f_i)\big)$
implies $W_i \le Z_i$.

Applying the above analysis to each pair of adjacent letters of $\BB$, we conclude that, given $k$, the relative order of $Z_0, \ldots, Z_n, W_0, \ldots, W_n$ is determined by  $\pi_\BB$. In particular, since $i$ precedes $\i$ in $\pi_\BB$ for all $i$, we have that $Z_i \geq W_i$ for all $i$, that is, $\varphi$ belongs to the chart $\mathscr{C}_k$.
\end{proof}

\subsection{The bipermutohedral fan: an origin story}\label{sec:whybiperm}

The bipermutohedral fan can be approached from several different points of view, and it has many favorable properties, as the previous sections show. However, it may not yet be clear where this fan comes from, or why it is a good setting for the Lagrangian geometry of matroids. In this section we explain the geometric motivation for its construction, and the role its various properties play in the theory.

When $\M$ is the matroid of a subspace $V$ of $\mathbb{C}^E$, the conormal fan $\Sigma_{\M, \M^\perp}$ is a tropical model of the projectivized conormal bundle of $V$. Since $\M^\perp$ is the matroid of the orthogonal complement of $V$, we expect the conormal fan to be supported on $\trop(\M) \times \trop(\M^\perp)$. 
A desirable  fan structure $\Sigma$ on this support should have the following properties:
\begin{enumerate}[(1)]\itemsep 5pt
\item The classes $\gamma$ and $\delta$ can be defined in its Chow ring, so we can state Theorems \ref{CSMTheorem} and \ref{degtheorem}.
\item The Chow ring is tractable for computations, so we can prove Theorems \ref{CSMTheorem} and \ref{degtheorem}.
\item The fan is a subfan of the normal fan of a polytope, so its ample cone is nonempty.
\item The fan is Lefschetz, so we can derive Conjecture \ref{NumericalConjectures} in Theorem \ref{thm:conjs}.
\end{enumerate}

We resolve requirement (4) by showing in Theorem \ref{thm:lefschetz} that being Lefschetz only depends on the support $\trop(\M) \times \trop(\M^\perp)$ -- which is the support of a  Lefschetz fan by \cite{AHK15} -- and not on the fan structure that we choose. Thus we can focus on the first three.

Requirement (2) is stated imprecisely, but a very desirable initial property is that our fan $\Sigma$ is simplicial. When this is the case, the Chow ring $A(\Sigma)$ of the toric variety $X(\Sigma)$ has an algebraic combinatorial presentation due to Brion, and an interpretation in terms of piecewise polynomial functions due to Billera; see Section \ref{sec:homology}.
This will allow us to carry out intersection-theoretic computations in this Chow ring.
Thus the first fan structure on $\trop(\M) \times \trop(\M^\perp)$ that one might try is the product of Bergman fans $\Sigma_{\M} \times \Sigma_{\M^\perp}$, which is simplicial and does have a nice combinatorial structure. It is also a subfan of the normal fan of the product of permutohedra $\Pi_E \times \Pi_E$.

To address requirement (1), we rely on the geometry of the representable case, as studied in \cite{DGS, HuhML}, which tells us what the classes $\gamma$ and $\delta$ of Theorems \ref{CSMTheorem} and \ref{degtheorem} should be. If $\alpha$ is the piecewise linear function on the tropical projective torus defined in the introduction, 
then $\gamma$ and $\delta$ should be the pullbacks of $\alpha$ along  the following maps from $\N_E \times \N_E$ to $\N_E$:
\[
\pi\colon \Sigma \longrightarrow \Gamma_E, \ \ (z,w) \longmapsto z \ \ \text{and} \ \ 
\mu\colon \Sigma \longrightarrow \Gamma_E, \ \ (z,w) \longmapsto z+w,
\]
where $\Gamma_E$ is the reduced normal fan of the standard simplex. If $\Sigma = \Sigma_{\M} \times \Sigma_{\M^\perp}$ or any refinement of it, the first map is a map of fans, and $\gamma$ \emph{is} well defined. However, the second map is \emph{not} a map of fans for $\Sigma = \Sigma_{\M} \times \Sigma_{\M^\perp}$, as we will see in Example \ref{ex:criticalfan}. Thus the product fan structure will not serve our purposes; we need to subdivide it further. How might we do this simultaneously for all $\M$?

At this point, it is instructive to return to the case of tropical linear spaces, as used by Adiprasito, Huh, and Katz in \cite{AHK15}. In that case, one wants a similarly convenient fan structure for the tropical linear space $\trop(\M)$. Fortunately, one can do this for all matroids on $E$ at once, by intersecting $\trop(\M)$ with the permutohedral fan $\Sigma_E$. The result is the \emph{Bergman fan} $\Sigma_{\M}$ of $\M$ introduced by Ardila and Klivans in \cite{AK06}, where it is called the \emph{fine subdivision}.

Similarly, we might try to find a suitable fan structure of $\trop(\M) \times \trop(\M^\perp)$ for all matroids $\M$ on $E$ simultaneously, by intersecting them with an appropriate complete fan. There is a natural candidate: the coarsest common refinement of the product of permutohedral fans $\Sigma_E \times \Sigma_E$, which induces the fan structure $\Sigma_{\M} \times \Sigma_{\M^\perp}$, 
and $\mu^{-1}(\Gamma_E)$, the coarsest fan that guarantees that the class $\delta$ is well-defined. This is the reduced normal fan of a polytope
\[
H_{E,E} \coloneq (\Pi_E \times \Pi_E) + D_E,
\]
the Minkowski sum of the product of two permutohedra and the diagonal simplex $D_E = \conv\{\e_i + \f_i\}_{i \in E}$. The polytope $H_{E,E}$ does have an elegant combinatorial structure, as shown in \cite{ArdilaEscobar}. They call it the \emph{harmonic polytope} because its number of vertices is
\[
 |E|! \cdot |E|! \cdot  \Bigg(1+\frac{1}{2}+\cdots+\frac{1}{|E|}\Bigg).
\]
However, this polytope has a drawback for our purposes: it is not simple, so the resulting fan structure on $\trop(\M) \times \trop(\M^\perp)$ is not simplicial.
Thus we need to find a simplicial refinement of the corresponding \emph{harmonic fan}, with simple enough combinatorial structure that we can carry out computations.

The bipermutohedral fan $\Sigma_{E,E}$ is our answer to those requirements. It refines the harmonic fan by Proposition \ref{prop:mu}, so $\gamma$ and $\delta$ are well-defined. It is simplicial by Proposition \ref{prop:unimodular}, and it has an elegant combinatorial structure that makes explicit computations possible. It is the normal fan of the bipermutohedron, thanks to Proposition \ref{thm:ardilahedron}.

For the above reasons, we define the conormal fan $\Sigma_{\M, \M^\perp}$ to be the fan on $\trop(\M) \times \trop(\M)^\perp$ obtained by intersecting with the bipermutohedral fan $\Sigma_{E,E}$. 
The construction of  $\Sigma_{\M, \M^\perp}$ is a Lagrangian analogue of the construction of $\Sigma_{\M}$ in \cite{AK06}.
What remains is to understand the resulting combinatorial structure and carry out the necessary intersection-theoretic computations to prove Theorems \ref{CSMTheorem} and \ref{degtheorem} -- which is the goal of Sections \ref{sec:matroidfan} and \ref{sec:degreecomps}  --
and to prove Theorem \ref{thm:lefschetz} -- which we do in Section \ref{sec:hodge}.
We believe that 
the bipermutohedral fan will have applications beyond those presented in this paper.
For example, the bipermutohedral perspective could be a guide in finding useful tropical models for Lagrangian matroids studied in \cite{BGW},
of which the conormal fan of a matroid will be a particular case.


\section{The conormal intersection theory of a matroid}\label{sec:matroidfan}

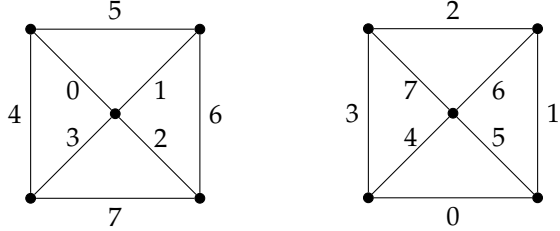
\begin{figure}[h]
\[
 \begin{tikzpicture}[baseline=(current bounding box.center),scale=0.75,
bend angle=25,vertex/.style={circle,draw,inner sep=1.3pt,fill=black}]

\node[vertex] (a) at (-0.5,-0.5) {};
\node[vertex] (b) at (2.5,-0.5) {};
\node[vertex] (c) at (2.5,2.5) {};
\node[vertex] (d) at (-0.5,2.5) {};
\node[vertex] (e) at (1,1) {};

\draw (a) to node[below] {7} (b);
\draw (b) to node[right] {6} (c);
\draw (c) to node[above] {5} (d);
\draw (d) to node[left] {4} (a);
\draw (a) to node[above] {3\ } (e);
\draw (b) to node[above] {\ 2} (e);
\draw (c) to node[below] {\ 1} (e);
\draw (d) to node[below] {0\ } (e);
\end{tikzpicture}
\qquad \qquad
\begin{tikzpicture}[baseline=(current bounding box.center),scale=0.75,
bend angle=25,vertex/.style={circle,draw,inner sep=1.3pt,fill=black}]

\node[vertex] (a) at (-0.5,-0.5) {};
\node[vertex] (b) at (2.5,-0.5) {};
\node[vertex] (c) at (2.5,2.5) {};
\node[vertex] (d) at (-0.5,2.5) {};
\node[vertex] (e) at (1,1) {};

\draw (a) to node[below] {0} (b);
\draw (b) to node[right] {1} (c);
\draw (c) to node[above] {2} (d);
\draw (d) to node[left] {3} (a);
\draw (a) to node[above] {4\ } (e);
\draw (b) to node[above] {\ 5} (e);
\draw (c) to node[below] {\ 6} (e);
\draw (d) to node[below] {7\ } (e);
\end{tikzpicture}
\]
\caption{The graph $G$ of the square pyramid and its dual graph $G^\perp$.}\label{fig:pyramid}
\end{figure}

In this section, we construct  the \emph{conormal fan}  of a matroid $\M$ on $E$,
and describe its Chow ring. Our running example will be the graphic matroid $\M(G)$ of the graph $G$ of the square pyramid, whose dual  is  the graphic matroid of the dual graph $G^\perp$ shown in Figure \ref{fig:pyramid}.

\subsection{Homology and cohomology}\label{sec:homology}

Throughout this section, we fix a simplicial rational fan $\Sigma$ in  $\N=\mathbb{R} \otimes \N_\mathbb{Z}$. For each ray $\rho$ in $\Sigma$, we write $\e_\rho$ for the primitive generator of $\rho$ in $\N_\mathbb{Z}$,
and introduce a variable $x_\rho$.
\begin{enumerate}[$\bullet$]\itemsep 5pt
\item Let $S(\Sigma)$ be the polynomial ring with real coefficients
that has $x_\rho$ as its variables, one for each ray $\rho$ of $\Sigma$.
\item Let $I(\Sigma)$ be the \emph{Stanley-Reisner} ideal of $S(\Sigma)$, generated by the square-free monomials indexing the subsets of rays of $\Sigma$ which do not generate a cone in $\Sigma$.
\item Let $J(\Sigma)$ be the ideal of $S(\Sigma)$ generated by the linear forms $\sum_\rho \ell(\e_\rho)x_\rho$, where
$\ell$ is any linear function on $\N$ and the sum is over all the rays in $\Sigma$.
\end{enumerate}

\begin{definition} \label{def:Chow}
The \emph{Chow ring} of $\Sigma$, denoted $A(\Sigma)$, is the  graded algebra $S(\Sigma)/(I(\Sigma)+J(\Sigma))$.
\end{definition}

Billera \cite{Bi89} constructed an isomorphism from the monomial quotient $S(\Sigma)/I(\Sigma)$  to the algebra of continuous piecewise polynomial functions on $\Sigma$ by identifying
  the variable $x_\rho$ with
  the piecewise linear  \emph{Courant function} on $\Sigma$ determined by the condition
\[
x_\rho(\e_{\rho'})=\begin{cases} 1, & \text{if $\rho$ is equal to $\rho'$}, \\ 0, &  \text{if $\rho$ is not equal to $\rho'$.} \end{cases}
\]
Thus, under this isomorphism, a piecewise linear function $\ell$ on $\Sigma$ is identified with the linear form
\[
\ell=\sum_\rho \ell(\e_\rho) x_\rho.
\]
We regard  the elements of $A(\Sigma)$ as equivalence classes of piecewise polynomial functions on $\Sigma$, modulo the restrictions of global linear functions to $\Sigma$.

Brion \cite[Theorem 6.7]{BrionEquivariant} showed that  the Chow ring
of the toric
variety $X(\Sigma)$ of $\Sigma$ with real coefficients
is isomorphic to $A(\Sigma)$.\footnote{In \cite{BrionEquivariant}, Brion identifies $A(\Sigma)$ with the Chow group of $X(\Sigma)$ with real coefficients. For the existence of the ring structure and the pullback, see \cite{Vistoli}.}
Under this isomorphism, the class of the torus orbit closure of a cone $\sigma$ in $\Sigma$
is identified with $\mult(\sigma) \, x_\sigma$,
where $x_\sigma$ is the monomial $\prod_{\rho \subseteq \sigma} x_\rho$ and $\mult(\sigma)$ is the index of the subgroup
\[
\Big( \sum_{\rho \subseteq \sigma} \mathbb{Z} \e_{\rho}\Big) \subseteq \N_\mathbb{Z} \cap \Big( \sum_{\rho \subseteq \sigma} \mathbb{R} \e_{\rho} \Big).
\]
All the fans appearing in this section will be unimodular, so $\mult(\sigma)=1$ for every  $\sigma$ in $\Sigma$.


We write $\Sigma(k)$ for the set of $k$-dimensional cones in $\Sigma$.
A \emph{$k$-dimensional Minkowski weight} on $\Sigma$ is a real-valued function
$\omega$ on $\Sigma(k)$ that satisfies the \emph{balancing condition}: For every $(k-1)$-dimensional cone $\tau$ in $\Sigma$,
\[
\text{$\sum_{\tau \subset \sigma} \omega(\sigma) \mathbf{e}_{\sigma/\tau}=0$ in the quotient space $\N/ \spn(\tau)$,}
\]
where $\mathbf{e}_{\sigma/\tau}$ is the primitive generator of the ray $( \sigma+ \spn(\tau) ) / \spn(\tau)$.
We say that $w$ is \emph{positive} if $w(\sigma)$ is positive for every $\sigma$ in $\Sigma(k)$.
We write $\mathrm{MW}_k(\Sigma)$ for the space of $k$-dimensional Minkowski weights on $\Sigma$, and set
\[
\MW(\Sigma)=\bigoplus_{k \ge 0} \MW_k(\Sigma).
\]
We will make use of the basic fact that the Chow group of  a toric variety is generated by the classes of torus orbit closures \cite[Lemma 12.5.1]{CLS}.
Thus, there is an injective linear map from  the dual of $A^k(\Sigma)$ to the space of $k$-dimensional weights on $\Sigma$,
whose image turns out to be  $\MW_k(\Sigma)$, as noted in \cite[Section 5]{AHK15}.\footnote{In \cite[Proposition 5.6]{AHK15}, all fans are unimodular and their Chow rings have integral coefficients. The same argument works for simplicial fans and their Chow rings with real or rational coefficients.}
Explicitly, the inverse isomorphism from the image is
\[
 \MW_k(\Sigma) \longrightarrow \Hom(A^k(\Sigma),\R), \qquad w \longmapsto (\mult(\sigma) x_\sigma \longmapsto w(\sigma)).
\]
Following \cite{AHK15}, 
we define
the \emph{cap product}, denoted $\eta \cap w$, using the composition
\[
A^\ell(\Sigma) \longrightarrow  \Hom(A^{k-\ell} (\Sigma),A^k(\Sigma)) \longrightarrow \Hom(\MW_k(\Sigma),\MW_{k-\ell}(\Sigma)),  \quad \eta \longmapsto (w \longmapsto \eta \cap w),
\]
where the first map is given by the multiplication in the Chow ring of $\Sigma$.
In short, $\MW(\Sigma)$ has the structure of a graded $A(\Sigma)$-module given by the isomorphism  $\MW(\Sigma) \simeq \Hom(A(\Sigma),\R)$.

Let  $f\colon\Sigma \to \Sigma'$ be a morphism of simplicial fans.
The pullback of functions define the \emph{pullback homomorphism} between the Chow rings
\[
f^* \colon A(\Sigma') \longrightarrow A(\Sigma),
\]
whose dual is the \emph{pushforward homomorphism} between the space of Minkowski weights
\[
f_* \colon \MW(\Sigma) \longrightarrow \MW(\Sigma').
\]
Since $f^*$ is a homomorphism of graded rings, $f_*$ is a homomorphism of graded modules.
\footnote{To see that the pullback between the Chow rings are determined by the pullback of piecewise linear functions, note that every divisor on simplicial toric variety is $\mathbb{Q}$-Cartier \cite[Proposition 4.2.7]{CLS} and that the pullback of Chern classes of line bundles corresponds to the pullback of piecewise linear functions \cite[Proposition 6.2.7]{CLS}.}
In other words, the pullback and the pushforward homomorphisms satisfy the \emph{projection formula}
\[
\eta \cap f_* w = f_*(f^* \eta \cap w).
\]

\subsection{The Bergman fan of a matroid}

The \emph{Bergman fan} of a matroid $\M$ on $E$, denoted $\Sigma_{\M}$, is the $r$-dimensional subfan of the $n$-dimensional permutohedral fan $\Sigma_E$
whose underlying set is the \emph{tropical linear space}
\[
\trop(\M)=\Big\{z  \, | \, \text{$\min_{i \in C}(z_i)$ is achieved at least twice for every circuit $C$ of $\M$}\Big\} \subseteq \N_E.
\]
The Bergman fan of $\M$ is equipped with the piecewise linear functions
\[
\alpha_j=\max_{i \in E} (z_j-z_i),
\]
and the space of linear functions on the Bergman fan is spanned by the differences
\[
\alpha_i-\alpha_j=z_i-z_j.
\]
Note that $\trop(\M)$ is nonempty if and only if $\M$ is  \emph{loopless}.
In the remainder of this section, we suppose that  $\M$ has no loops.
In this case, the Bergman fan of $\M$ is the induced subfan of $\Sigma_E$ generated by the rays corresponding to the nonempty proper flats of $\M$ \cite{AK06}.

\begin{proposition}\label{prop:BergmanFan}
The Bergman fan of ${\M}$ is the unimodular fan in $\N_E$ with the cones
\[
\sigma_{\F}=\cone\{\e_F\}_{ F \in \F}, \ \  \text{for flags of flats $\F$ of $\M$.}
\]
\end{proposition}

The most important geometric property of $\Sigma_{\M}$ is the following  description of its top-dimensional Minkowski weights.
For a proof, see, for example, \cite[Proposition 5.2]{AHK15}.

\begin{proposition}\label{PropositionUniqueBalancing}
An $r$-dimensional weight on $\Sigma_{\mathrm{M}}$ is balanced if and only if it is constant.
\end{proposition}

We write $1_{\M}$
for the \emph{fundamental weight} on $\Sigma_{\M}$,
the $r$-dimensional Minkowski weight on the Bergman fan that has the constant value $1$.

\subsection{The Chow ring of the Bergman fan}
In the context of matroids,
for simplicity, we set
\[
S_{\M}=S(\Sigma_{\M}), \quad I_{\M}=I(\Sigma_{\M}), \quad J_{\M}=J(\Sigma_{\M}), \quad A_{\M}=A(\Sigma_{\M}).
\]
We identify the elements of $S_{\M}/I_{\M}$ with the piecewise linear functions on $\Sigma_{\M}$ as before.

Let $x_F$ be the variable of the polynomial ring corresponding to  the ray generated by $\e_F$ in the Bergman fan.
For any set $\F$ of nonempty proper flats of $\M$, we write $x_{\F}$ for the monomial
\[
x_\mathscr{F}=
\prod_{F \in \mathscr{F}} x_F.
\]
The variable $x_F$, viewed as a piecewise linear function on the Bergman fan,  is given by
\[
x_F(\e_{F'})=\begin{cases} 1, & \text{if $F$ is equal to $F'$}, \\ 0, &  \text{if $F$ is not equal to $F'$,} \end{cases}
\]
and hence the piecewise linear function $\alpha_j$ on the Bergman fan satisfies  the identity
\[
\alpha_j=\sum_{F} \alpha_j(\e_F) x_F=\sum_{j \in F} x_F.
\]
Thus, in the above notation,
\begin{enumerate}[$\bullet$]\itemsep 5pt
\item $S_{\M}$ is the ring of polynomials in the variables $x_F$, where $F$ is a nonempty proper flat of $\M$,
\item $I_{\M}$ is the ideal  generated by the  monomials $x_{\F}$, where $\F$ is not a flag, and
\item $J_{\M}$ is the ideal  generated by the linear forms $\alpha_i-\alpha_j$, for any $i$ and $j$  in $E$.
\end{enumerate}
We write $\alpha$ for the common equivalence class of $\alpha_j$ in the Chow ring of the Bergman fan.

\begin{definition} \label{def:deg}
The fundamental weight $1_{\M}$  defines the \emph{degree map}
\[
\deg\colon A^r_{\M} \longrightarrow \R, \quad x_{\F} \longmapsto x_{\F} \cap 1_{\M}=\begin{cases} 1&\text{if $\F$ is a flag,} \\ 0& \text{if $\F$ is not a flag.}\end{cases}
\]
\end{definition}

By Proposition \ref{PropositionUniqueBalancing}, the degree map is an isomorphism.
In other words, for any maximal flag $\F$  of nonempty proper flats of $\M$, the class of the monomial
$x_{\F}$  in the Chow ring of the Bergman fan of $\M$ is nonzero and does not depend on $\F$.

\subsection{The conormal fan of a matroid}\label{ss:conormal}

The \emph{conormal fan} of a matroid $\M$ on $E$, denoted $\Sigma_{\M,\M^\perp}$, is the $(n-1)$-dimensional subfan of the $2n$-dimensional bipermutohedral fan $\Sigma_{E,E}$
whose support is the product of tropical linear spaces
\[
|\Sigma_{\M,\M^\perp}|=\trop(\M) \times \trop(\M^\perp).
\]
Equivalently, the conormal fan is the largest subfan of the bipermutohedral fan for which the projections to the factors are  morphisms of fans
\[
\pi\colon \Sigma_{\M,\M^\perp} \longrightarrow \Sigma_{\M} \quad \text{and} \quad
\overline \pi\colon \Sigma_{\M,\M^\perp} \longrightarrow \Sigma_{\M^\perp}.
\]
The addition map $(z,w) \mapsto z+w$ is also a morphism of fans $\Sigma_{\M,\M^\perp} \rightarrow \Gamma_E$.

The conormal fan of $\M$ is equipped with the piecewise linear functions
\[
\gamma_j=\max_{i \in E} (z_j-z_i),\qquad
\og_j=\max_{i \in E} (w_j-w_i), \qquad
\delta_j=\max_{i \in E} (z_j+w_j-z_i-w_i),
\]
which are the pullbacks of $\alpha_j$ under the projections $\pi$ and $\pi'$ and the addition map, respectively.
The space of linear functions on the conormal fan is spanned by the differences
\[
 \gamma_i-\gamma_j=z_i-z_j \quad \text{and}  \quad \og_i-\og_j=w_i-w_j.
  \]
Note that the support of the conormal fan of $\M$ is nonempty if and only if $\M$ is \emph{loopless} and \emph{coloopless}.
In the remainder of this section, we suppose that  $\M$ has no loops and no coloops.

\begin{definition}
A \emph{biflat} $F|G$ of $\M$ consists of a flat $F$ of $M$ and a flat $G$ of $\M^\perp$ that form a bisubset; that is, they are nonempty, they are not both equal to $E$, and their union is $E$. A \emph{biflag} of $\M$ is a flag of biflats.
\end{definition}

We give an analog of Proposition \ref{prop:BergmanFan} for conormal fans in terms of biflats.

\begin{proposition}
The conormal fan of  $\M$ is the unimodular fan in $\N_{E,E}$ with the cones
\[
\sigma_{\F|\G}=\cone\{\e_{F|G}\}_{F|G \in \F|\G}, \ \  \text{for flags of biflats $\F|\G$ of $\M$.}
\]
\end{proposition}

\begin{proof}
The proof is straightforward, given  Corollary \ref{cor:bipermutohedralrays} and Proposition \ref{prop:BergmanFan}:
If  $\F|\G$ is a flag of biflats of $\M$, then $\F$ is an increasing sequence of flats of $\M$ and $\G$ is a decreasing sequence of flats of $\M^\perp$,
and hence
\[
\sigma_{\F|\G}\subseteq \sigma_{\F} \times \sigma_{\G} \in \Sigma_{\M} \times \Sigma_{\M^\perp}.
\]
Therefore, the conormal fan of $\M$ contains the induced subfan of $\Sigma_{E,E}$ generated by the rays corresponding to the biflats of $\M$.
The other inclusion follows from the easy implication
\[
\text{$\e_{F|G}$ is in the support of the conormal fan of $\M$} \Longrightarrow \text{$F|G$ is a biflat of $\M$}. \qedhere
\]
\end{proof}

We also have the following analog of Proposition \ref{PropositionUniqueBalancing} for conormal fans.

\begin{proposition}\label{PropositionUniqueBalancing2}
An $(n-1)$-dimensional weight on $\Sigma_{\M,\M^\perp}$ is balanced if and only if it is constant.
\end{proposition}

We write $1_{\M,\M^\perp}$ for the \emph{fundamental weight} on $\Sigma_{\M,\M^\perp}$,
the top-dimensional Minkowski weight on the conormal fan that has the constant value $1$.

\begin{proof}
Proposition \ref{PropositionUniqueBalancing} applied to $\M$ and $\M^\perp$ shows that a top-dimensional weight on $\Sigma_{\M} \times \Sigma_{\M^\perp}$ satisfies the balancing condition if and only if it is constant.
This property of the fan
remains invariant under any subdivision of its support, as shown in \cite[Section 2]{GKM}.
\end{proof}

For our purposes, the product of the Bergman fans of $\M$ and $\M^\perp$ has a shortcoming:
The addition map need not be a morphism from the product to  the fan $\Gamma_E$.
Thus, in general, we cannot define the class of $\delta_j$ in the Chow ring of the product. This is our motivation for subdividing it further, to obtain the conormal fan $\Sigma_{\M,\M^\perp}$.

\begin{example}\label{ex:criticalfan}
Let $\M$ and $\M^\perp$ be the graphic matroids of the graphs in Figure \ref{fig:pyramid}.
Consider the cone $\sigma_{\F} \times \sigma_{\G}$ in the product  of Bergman fans of $\M$ and $\M^\perp$, where
\[
\F = (\varnothing \subsetneq 1 \subsetneq 015 \subsetneq 01345\subsetneq E) \quad \text{and} \quad \G = (\varnothing \subsetneq 2 \subsetneq 267 \subsetneq 12567 \subsetneq E).
\]
This cone is subdivided into the  chambers of $\Sigma_{\M,\M^\perp}$ corresponding to the  biflags
\begin{align*}
& \begin{array}{|cc|ccccccccccc|cc|}
\hline
\varnothing &\subsetneq&1 & \subseteq & 015 & \subseteq & 01345 & \subseteq &  01345 & \subseteq & 01345 & \blue{\subseteq} & E &\subseteq&E \\
E& \supseteq& E &  \supseteq  & E & \supseteq & E & \supseteq & 12567 & \supseteq & 267 & \blue{\supseteq} & 2 &\supsetneq& \varnothing\\
\hline
\end{array} \ , \\[1mm]
 &  \begin{array}{|cc|ccccccccccc|cc|}
\hline
\varnothing & \subsetneq & 1 & \subseteq & 015 & \blue{\subseteq} & 01345 & \subseteq &  01345 & \subseteq & \ \ E\  \ & \subseteq & E & \subseteq & E\\
E & \supseteq & E & \supseteq  & E & \blue{\supseteq} & 12567 & \supseteq & 267 & \supseteq & \ \ 267 \ \ & \supseteq & 2 & \supsetneq & \varnothing \\
\hline
\end{array} \ , \\[1mm]
 & \begin{array}{|cc|ccccccccccc|cc|}
\hline
\varnothing & \subsetneq & 1 & \subseteq & 015 & \blue{\subseteq} & 01345 & \subseteq &  E & \subseteq & E & \subseteq & E & \subseteq & E\\
E & \supseteq & E & \supseteq  & E & \blue{\supseteq} & 12567 & \supseteq & 12567  & \supseteq &  \ \ 267 \ \  & \supseteq & 2 & \supsetneq & \varnothing  \\
\hline
\end{array} \ .
\end{align*}
If $(z,w)$ is inside the first chamber, then the minimum of $z_i+w_i$ is attained by $z_6+w_6=z_7+w_7$, and hence $z+w$ is in the cone  $\sigma_{012345}$. If $(z,w)$ is inside the second or the third chamber, then the minimum of $z_i+w_i$ is attained by $z_3+w_3=z_4+w_4$, and hence $z+w$ is in the cone $\sigma_{012567}$.
Thus, the product cone  does not map into a cone in $\Gamma_E$ under the addition map. \end{example}

Recall from Definition \ref{def:cotangent} that the cotangent fan $\Omega_E$ is the subfan of $\Sigma_{E,E}$ with support
\[
\trop(\delta)=\Big\{ (z,w)  \ | \  \text{$\min_{i \in E}(z_i+w_i)$ is achieved at least twice}\Big\} \subseteq \N_{E,E}.
\]
In other words, the cotangent fan is the collection of cones $\sigma_{\BB}$ for bisequences $\BB$ on $E$, where at least two elements of $E$ appear exactly once in $\BB$.
We show that the cotangent fan contains all the conormal fans of matroids on $E$.

\begin{proposition}\label{prop:cotangent}
For any matroid $\M$ on $E$, we have
$
\trop(\M) \times \trop(\M^\perp) \subseteq \trop(\delta)$.
\end{proposition}

In other words, if the minimum of $(z_i)_{i \in C}$ is achieved at least twice for every circuit $C$ of $\M$ and
the minimum of $(w_i)_{i \in C^\perp}$ is achieved at least twice for every circuit $C^\perp$ of $\M^\perp$, then
the minimum of $(z_i+w_i)_{i \in E}$  is achieved at least twice.
We deduce Proposition \ref{prop:cotangent}  from Proposition \ref{prop:gap} below, a stronger statement on the flags of biflats of $\M$.
The notion of gaps introduced here for Proposition \ref{prop:gap} will be useful in Section \ref{sec:degreecomps}.

Let $\F|\G$ be a flag of biflats of $\M$.
As before, we write $\F$ and $\G$ for the sequences
\begin{align*}
\mathscr{F} &= (\varnothing \subsetneq F_1 \subseteq \cdots \subseteq F_k \subseteq E), \ \  \text{where $F_j$ are the first parts of the biflats in $\F|\G$}, \\
\mathscr{G} &= (E \supseteq G_1 \supseteq \cdots \supseteq G_k \supsetneq \varnothing), \ \  \text{where $G_j$ are the second parts of the biflats in $\F|\G$},
\end{align*}
where $k$ is the length of $\mathscr{F}|\mathscr{G}$.
Thus, the bisequence $\BB(\F|\G)$ from Proposition \ref{prop:cones} can be written
 \[
 B_0|B_1|\cdots|B_k, \ \ \text{where $B_j=(F_{j+1} - F_j) \cup (G_j - G_{j+1})$}.
 \]

 \begin{definition} \label{def:gaps}
The \emph{gap sequence} of $\F|\G$, denoted $\D(\F|\G)$,
is the sequence of \emph{gaps}
 \[
 D_0|D_1|\cdots|D_k, \ \ \text{where $D_j=(F_{j+1} - F_j) \cap (G_j - G_{j+1})$}.
 \]
 \end{definition}

Note that $D_j$ consists of the elements of $B_j$ that appear exactly once in the bisequence $\BB(\F|\G)$.

\begin{example}
The three maximal flags of biflats shown in Example \ref{ex:criticalfan} have the gap sequences
\[
\varnothing|\varnothing|\varnothing|\varnothing|\varnothing|67|\varnothing, \quad
\varnothing|\varnothing|34|\varnothing|\varnothing|\varnothing|\varnothing, \quad
\varnothing|\varnothing|34|\varnothing|\varnothing|\varnothing|\varnothing.
\]
We show in Proposition \ref{prop:doublejump}  that any maximal flag of biflats  has a unique nonempty gap.
\end{example}

\begin{lemma}\label{lem:gap}
The complement of the gap $D_j$ in $E$ is the union of $F_j$ and $G_{j+1}$.
\end{lemma}

Therefore, by Proposition \ref{prop:alt}, at least one of the gaps of $\F|\G$ must be nonempty.

\begin{proof}
Since $F_j|G_j$ and $F_{j+1}|G_{j+1}$ are bisubsets, we have $G_j^{c} \subseteq F_j$ and $F_{j+1}^c \subseteq G_{j+1}$. Thus,
\[
D_j^c = (F_{j+1} \cap F_j^c \cap G_j \cap G_{j+1}^c) ^c= F_{j+1}^c \cup F_j \cup G_j^c \cup G_{j+1} =F_{j} \cup G_{j+1}. \qedhere
\]
\end{proof}

\begin{lemma}\label{lem:nongaps}
Let $e\in E$. There exists an index $i$ for which $e\in F_i\cap G_i$ if and only if $e$ is not in any gap. In symbols, the union of the gaps of $\F|\G$ is
\[
\bigsqcup_{j=0}^k D_j = E - \bigcup_{i=1}^k \left(F_i \cap G_i\right).
\]
\end{lemma}
\begin{proof}
First suppose $e\in F_i\cap G_i$.  Then $e\in F_j$ for all $j\geq i$,
which means $e\not \in D_j$ for $i\leq j\leq k$.
Dually, $e\in G_j$ for all $j\leq i$, so $e\not\in D_j$ for all
$0\leq i\leq j-1$.

Now suppose $e$ is not in any gap, and consider the index $1\leq i\leq k+1$ for which $e\in F_i-F_{i-1}$.  Since $e\in F_{i-1}\cup G_i$, we must have $e\in G_i$ and hence $e \in F_i\cap G_i$.
\end{proof}

We will often use the following basic result.  Recall that $\abs{E}=n+1$.

\begin{lemma}\label{lem:big unions}
The union of a flat and a coflat cannot have exactly $n$ elements.
\end{lemma}
\begin{proof}
Let $F$ be a flat and $G$ be a coflat.  Recall that, for any matroid,
the complement of any hyperplane is a cocircuit
\cite[Proposition~2.1.6]{Oxley} and that any flat is an intersection of
hyperplanes \cite[Proposition 1.7.8]{Oxley}.  So we may write the complement of
$F\cup G$ as
\[
\Big( \bigcup_{C \in \C} C \Big)\cap\Big( \bigcup_{C^\perp \in
  \C^\perp} C^\perp \Big),
\]
where $\C$ is a collection of circuits and $\C^\perp$ is a collection
of cocircuits.  Thus, if the complement is nonempty,
there are $C \in \C$ and $C^\perp \in \C^\perp$ that intersect
nontrivially.  Now the conclusion follows from the classical fact that
the intersection of a circuit and a cocircuit is either empty or
contains at least two elements \cite[Proposition 2.11]{Oxley}.
 \end{proof}

 \begin{proposition}\label{prop:gap}
Every nonempty gap of a biflag $\F|\G$  of $\M$ has at least two distinct elements.
\end{proposition}

 \begin{proof}
   Since the complement of a gap of $\F|\G$ is the union of a flat and a coflat
   by Lemma~\ref{lem:gap}, the claim follows from Lemma~\ref{lem:big unions}.
 \end{proof}

For any biflag $\F|\G$, there are at least two elements of $E$ that appear exactly once in the bisequence $\BB(\F|\G)$; therefore
\[
\trop(\M) \times \trop(\M^\perp) \subseteq \trop(\delta),
\]
proving Proposition \ref{prop:cotangent}.

We will often use the following restatement of Proposition \ref{prop:gap}.  Recall that $|E|=n+1$.

For later use, we record here another elementary property of the flags of biflats of a matroid.

\begin{definition}\label{jumpsets}
The \emph{jump sets} of $\F$ and $\G$ are the sets of indices
\[
\mathsf{J}(\F)=\{j \, | \, \text{$0 \le j \le k$ and $F_j \neq F_{j+1}$} \} \ \ \text{and} \ \
\mathsf{J}(\G)=\{j \, | \, \text{$0 \le j \le k$ and $G_j \neq G_{j+1}$} \}.
\]
The elements of  $\mathsf{J}(\F) \cap \mathsf{J}(\G)$ are called the \emph{double jumps} of $\F|\G$.
\end{definition}

The double jumps are colored blue in the table of $\F|\G$, as shown in Example \ref{ex:criticalfan}. Clearly,  $j$ is a double jump whenever the corresponding gap $D_j$ is nonempty. We show that the converse holds when  $\F|\G$ is maximal.

\begin{proposition}\label{prop:doublejump}
Every maximal flag of biflats $\F|\G$ of $\M$ has a unique double jump.
Ignoring repetitions, $\F$ and $\G$ are complete flags of nonempty flats in $\M$ and $\M^\perp$, respectively.
\end{proposition}

In particular,
every maximal flag of biflats $\F|\G$ of $\M$ has a unique nonempty gap.

\begin{proof}
Recall that at least one of the gaps of $\F|\G$ is nonempty.
In addition, since tropical linear spaces are pure-dimensional,
the length of any maximal flag of biflats must be $n-1$.
Thus,
\[
|\mathsf{J}(\F)\cap \mathsf{J}(\G)| \ge 1 \ \ \text{and} \ \
|\mathsf{J}(\F)\cup \mathsf{J}(\G)|=n.
\]
On the other hand, writing $r+1$ for the rank of $\M$ as before, we have
\[
|\mathsf{J}(\F)| \le r+1 \ \ \text{and} \ \  |\mathsf{J}(\G)| \le n-r.
\]
Therefore,
$n+1 \le |\mathsf{J}(\F)\cup \mathsf{J}(\G)|+|\mathsf{J}(\F)\cap \mathsf{J}(\G)|=|\mathsf{J}(\F)| +|\mathsf{J}(\G)| \le n+1$, and hence
\[
|\mathsf{J}(\F)| = r+1, \ \ \ \  |\mathsf{J}(\G)| = n-r \ \ \text{and} \ \
|\mathsf{J}(\F)\cap \mathsf{J}(\G)| = 1
\]
which imply the desired results.
\end{proof}

By way of contrast, nonmaximal biflags have several double jumps, and they can have a double jump whose corresponding gap is empty:
\begin{example}
 For the graphic matroids of Figure~\ref{fig:pyramid} again, consider the biflag
  \[
  \F|\G \coloneq 
\begin{array}{|cc|ccccc|cc|}
\hline
\varnothing &\subsetneq& 1 & \blue{\subseteq} & 01345 & \blue{\subseteq} & E &\subseteq&E \\
E& \supseteq& E & \blue{\supseteq} & 12567 &  \blue{\supseteq} & 267 & \supsetneq & \varnothing\\
\hline
\end{array}\,.
\]
We see that $\F|\G$ has two double jumps, with gaps $034$ and $\varnothing$,
respectively.  In view of Proposition~\ref{prop:doublejump}, any maximal
biflag containing $\F|\G$ would necessarily contain another biflat $F|G$
satisfying
$01345\subseteq F\subseteq E$ and $12567\supseteq G\supseteq 267$.
\end{example}

\subsection{The Chow ring of the conormal fan}
For notational simplicity, we set
\[
\SSS=S(\Sigma_{\M,\M^\perp}), \quad \III=I(\Sigma_{\M,\M^\perp}), \quad \JJJ=J(\Sigma_{\M,\M^\perp}), \quad \PAA=A(\Sigma_{\M,\M^\perp}).
\]
We identify the elements of $\SSS/ \III$ with the piecewise linear functions on the conormal fan.

Let $x_{F|G}$ be the variable of the polynomial ring
 corresponding to the ray generated by $\e_{F|G}$ in the conormal fan.
 For any set  $\F|\G$ of biflats  of $\M$, we write $x_{\F|\G}$ for the monomial
\[
x_{\F|\G}=\prod_{F|G \in \F|\G} x_{F|G}.
\]
We note that the piecewise linear function $\delta_j$ on the conormal fan satisfies the identity
\[
\delta_j=\sum_{F | G} \delta_j(\e_{F|G})x_{F|G}
=\sum_{j \in F \cap G} x_{F|G}.
\]
Similarly, the piecewise linear functions $\gamma_j$ and $\og_j$ satisfy the identities
\[
\gamma_j=\sum_{j \in F \neq E} x_{F|G} \ \  \text{and} \ \
\overline \gamma_j=\sum_{j \in G \neq E} x_{F|G}.
\]
Thus, in the above notation,
\begin{enumerate}[$\bullet$]\itemsep 5pt
\item $\SSS$ is the ring of polynomials in the variables $x_{F|G}$, where $F|G$ is a biflat of $\M$,
\item $\III$ is the ideal  generated by the  monomials $x_{\F|\G}$, where $\F|\G$ is not a biflag, and
\item $\JJJ$ is the ideal  generated by the linear forms $\gamma_i-\gamma_j$ and $\og_i-\og_j$, for any $i$ and $j$  in $E$.
\end{enumerate}
We write $\gamma$, $\og$, and $\delta$, respectively, for the  equivalence classes of $\gamma_j$, $\og_j$, and  $\delta_j$
in the Chow ring of the conormal fan.

\begin{definition}\label{def:conormaldegree}
The fundamental weight $1_{\M,\M^\perp}$ of the conormal fan defines the \emph{degree map}
\[
\deg\colon A^{n-1}_{\M, \M^\perp} \longrightarrow \R, \quad x_{\F|\G} \longmapsto x_{\F|\G} \cap 1_{\M,\M^\perp}=\begin{cases} 1&\text{if $\F|\G$ is a biflag,} \\ 0& \text{if $\F|\G$ is not a biflag.}\end{cases}
\]
\end{definition}

By Proposition \ref{PropositionUniqueBalancing2}, the degree map is a linear isomorphism.
In other words,   for  maximal flag of biflats $\F|\G$ of $\M$,
the class of the monomial
$x_{\F|\G}$  in the Chow ring of the conormal fan of $\M$ is nonzero and does not depend on $\F|\G$.

Recall that the projection $\pi$ is a morphism from the conormal fan of $\M$ to the Bergman fan of $\M$.
The projection has the special property that the image of a cone in the conormal fan is a cone in the Bergman fan (and not just contained in one).
This property leads to the following simple description of the pullback  $\pi^*\colon A_{\M} \rightarrow {\PAA}$.

\begin{proposition}\label{prop:pullbacksum}
For any flag of nonempty proper flats $\F$ of $\M$,
\[
\pi^*(x_\F)=\sum_{\G} x_{\F|\G},
\]
where the sum is over all decreasing sequences $\G$ such that $\F|\G$ is a flag of biflats of $\M$.
\end{proposition}

Dually, the pushforward of any Minkowski weight $w$ on the conormal fan is given by
\[
\pi_*(w)(\sigma_{\F})=\sum_{\G} w(\sigma_{\F|\G}),
\]
where the sum is over all decreasing sequences $\G$ such that $\F|\G$ is a flag of biflats of $\M$.

\begin{proof}
Since $\pi(\e_{F|G})=\e_F$, the pullback of the piecewise linear function $x_F$ satisfies
\[
\pi^*(x_F)=\sum_{G} x_{F|G},
\]
where the sum is over all $G$ such that $F|G$ is a biflat of $\M$. Thus, for any given $\F$,
\[
\pi^*(x_{\F})=\prod_{F \in \F} \pi^*(x_F)=\sum_{\G} x_{\F|\G},
\]
where the sum is over all decreasing sequences $\G$ such that $\F|\G$ is a flag of biflats of $\M$.
\end{proof}

\section{Degree computations in  the conormal fan} \label{sec:degreecomps}

Throughout this section, we fix a  matroid $\M$ of rank $r+1$ on the ground set $E=\{0,1,\ldots,n\}$. 
We fix the usual ordering on the ground set 
\[
0<1<\cdots<n.
\]
We aim to evaluate various elements of  $A(\Sigma_{\M,\M^\perp})$ under the degree map (Definition \ref{def:conormaldegree}).

Let  $\F=(F_1 \subsetneq \cdots \subsetneq F_k)$ be a flag of nonempty proper flats of $\M$.
The \emph{beta invariant of $\F$} is
\[
\beta_{\M[\F]} =  \prod_{i=0}^{k} \beta_{\M[F_{i},F_{i+1}]},
\]
where $\beta_{\M[F_{i-1},F_i]}$ is the beta invariant of the minor $\M[F_{i-1},F_i] = M|F_i/F_{i-1}$.\footnote{We continue to use the convention that $F_0=\varnothing$ and $F_{k+1}=E$.}
The main goal of this section is to prove Proposition \ref{prop:betasF} in Section \ref{flagbeta}, which states
the identity 
\[
\deg(\pi^*(x_\F)\delta^{n-k-1})
= \beta_{\M[\F]}.
\]
In particular,  the degree of the conormal fan with respect to $\delta$ is the beta invariant:
\[
\deg(\delta^{n-1}) = \beta_{\M}.
\]
The result will be used  in Section \ref{sec:CSM} to prove Theorems \ref{CSMTheorem} and \ref{degtheorem} stated in the introduction.

Since  $\pi^*(x_\F)=\sum x_{\F|\G}$ by Proposition \ref{prop:pullbacksum},
it is enough to compute the degree of   $x_{\F|\G} \, \delta^{n-k-1}$ for all possible $\G$.
We will show in Lemma \ref{lem:vanishing} that, in fact, 
the degree is nonzero for at most one $\G$.
The degree computation will require us to study more closely the combinatorial structure of conormal fans, and develop algebraic combinatorial techniques for computing in their Chow rings. 

\subsection{Canonical expansions in  the conormal fan}

In order to compute 
 the degree of $x_{\F|\G} \, \delta^{n-k-1}$, 
 we seek to express it as a sum of square-free monomials, each of which have degree one.
One fundamental feature of this computation, which is simultaneously an advantage and a difficulty, is that there are many  ways to carry it out. 
We may choose any one of the 
different expressions for $\delta$ to compute, namely $\delta=\delta_i$ for each $i$ in $E$. To have control over the computation, we require some structure amidst that freedom. 

For every nonnegative integer $m$, we prescribe a canonical way of expressing 
$x_{\F|\G} \, \delta^m$ as a sum of square-free monomials.
Let  $e = e(\F|\G)$ be the largest gap element of $\F|\G$, which exists by Lemma \ref{lem:nongaps}.
In the notation of that lemma,  we have
\[
e = \max \Big(\bigsqcup_{j=0}^k D_j \Big) = \max \Big(E - \bigcup_{i=1}^k (F_i \cap G_i) \Big),
\]
where $D_0,\ldots,D_k$ is the gap sequence of $\F|\G$ defined in Definition \ref{def:gaps}.

\begin{definition}  \label{def:canonical} 
The \emph{canonical expansion of $x_{\F|\G}\,  \delta$} is the expression
\[
x_{\F|\G} \,\delta
=
x_{\F|\G} \,\delta_{e}
=
 \sum_{
e\in F \cap G
} x_{\F|\G} x_{F|G},
\]
where the sum is over all biflats $F|G$ such that $e \in F \cap G$.
We recursively obtain the \emph{canonical expansion of $x_{\F|\G}\,  \delta^m$}  by multiplying each monomial in the canonical expansion of $x_{\F|\G}\,  \delta^{m-1}$ by $\delta$, again using the canonical expansion.
\end{definition}

The canonical expansions are sums of square-free monomials in $A_{\M,\M^\perp}$.
Note that some or all of its summands may be zero in the Chow ring. The following lemma describes the nonzero terms.

\begin{lemma}\label{lem:canonexpansion}
If a summand $x_{\F|\G} x_{F|G}$ of the canonical expansion of $x_{\F|\G} \,\delta$ is nonzero 
and  $e$ is in  $D_j$, then  $F_j \subseteq F \subseteq F_{j+1}$ and $G_j \supseteq G \supseteq G_{j+1}$.
 \end{lemma}

\begin{proof}
If $\sigma_{\F \cup F | \G \cup G}$ is a cone in the conormal fan with $e \in F \cap G$, then $e \notin F_j$ and $e \notin G_{j+1}$.
Thus, the biflat  $F|G$ must be added to $\F |\G$ in between the indices $j$ and $j+1$. 
\end{proof}

We may think of the canonical expansion of $\delta^m$ as a recursive procedure to produce a list of $m$-dimensional cones in the conormal fan $\Sigma_{\M,\M^\perp}$, where each cone is built up one ray at a time according to the rules prescribed in Lemma \ref{lem:canonexpansion}.

\begin{example}\label{ex:pyramid2}
For the graph $G$ of the square pyramid in Figure \ref{fig:pyramid}, the canonical expansion of the highest nonzero power of $\delta$ in $A_{\M, \M^\perp}$ is 
\begin{align*}
\delta^6 =&  x_{6|E}\, x_{56|E}\, x_{4567|E}\, x_{E|23467}\, x_{E|347}\, x_{E|7} \\
& + x_{7|E}\, x_{67|E}\, x_{4567|E}\, x_{E|235}\, x_{E|35}\, x_{E|5}\\
& + x_{7|E}\, x_{57|E}\, x_{4567|E}\, x_{E|23467}\, x_{E|36}\, x_{E|6}. 
\end{align*}
This expression is deceivingly short. Carrying out this seemingly simple computation by hand is very tedious; if one were to do it by brute force, one would find that
the number of terms of the canonical expansions of $\delta^0, \ldots, \delta^6$ are the following:
\[
\begin{array}{|l||l|l|l|l|l|l|l|} \hline
& \delta^0 & \delta^1 & \delta^2 & \delta^3 & \delta^4 & \delta^5 & \delta^6\\
\hline \hline
\text{number of monomials counted with multiplicities} & 1 & 29 & 352 & 658 & 383 & 69 & 3\\ \hline
\text{number of distinct monomials} & 1 & 29 & 333 & 621 & 370 & 68 & 3\\ \hline
\end{array} \ .
\]
This example shows typical behavior: for small $k$ the number of cones in the
expansion of $\delta^k$ increases with $k$, but as $k$ approaches
$n-1$, increasingly many products $x_{\F|\G} \, \delta$ are zero, and the canonical expansions become shorter.
\end{example}

We record an explicit description of the canonical expansion of powers of $\delta$ in the following proposition. 

\begin{proposition}\label{prop:expansion}
For each $m$, the canonical expansion of $\delta^m$  is given by
\[
\delta^{m}=\sum_{(\F|\G, \e)
} x_{\F|\G},
\]
where the sum is over 
 all pairs 
$\F|\G=(F_1|G_1,\ldots,F_m|G_m)$ 
and
  $\e=(e_1, \ldots, e_m)$  satisfying
\[
e_i\in F_i\cap G_i \ \ \text{and} \ \ e_i=\max\Big(E-\bigcup_{ e_j > e_i} (F_j\cap G_j)\Big) \ \  \text{for all $1\leq i\leq m$}.
\]
\end{proposition}

\begin{proof}
The displayed formula for $\delta^m$ is an immediate consequence of the  construction of the canonical decomposition.
\end{proof}

It will be convenient to encode each summand of the canonical expansion 
 in a table:
\[
\begin{array}{|cc|ccccccccccc|cc|}
\hline
\varnothing & \subsetneq & F_1 & \subseteq & \cdots & \subseteq & F_i &
 \subseteq & \cdots & \subseteq  & F_m & \subseteq & E\\
E & \supseteq & G_1 & \supseteq & \cdots & \supseteq & G_i &
 \supseteq & \cdots & \supseteq  & G_m & \supsetneq & \varnothing\\
\hline
&&e_1 && \cdots && e_i && 
\cdots && e_m&&\\
\hline
\end{array} \,
\]
The canonical expansion of $\delta^m$ may contain repeated terms $x_{\F|\G}$ coming from tables that have the same biflag $\F|\G$ but different sequences $\e$,
as the numerics for the canonical expansions of $\delta^2,\delta^3,\delta^4,\delta^5$ in Example \ref{ex:pyramid2} show.
On the other hand, we will see in Proposition  \ref{prop:degreebetaM} that the canonical expansion of $\delta^{n-1}$ does not contain repeated terms.

\begin{example}
We revisit the canonical expansion of $\delta^6$ 
in Example~\ref{ex:pyramid2}. The first monomial arises from the following table: 
\[
 \begin{array}{|cc|ccccccccccc|cc|}
\hline
\varnothing & \subset&   6 & \subsetneq & 56 & \subsetneq &  4567 & \blue{\subsetneq} & E & = & E & = & E & = &  E \\
{}E & = & E & = & E & = & E & \blue{\supsetneq} & 23467 &
\supsetneq & 347 & \supsetneq &  7 & \supset&  \varnothing \\ \hline
&&6 & & 5 &&  4 &&  2 &&  3 && 7 &&\\
\hline
\end{array}
\]
The variables $x_{F_i|G_i}$ arrive to the monomial in the order $x_{E|{ 7}}x_{{ 6}|E}x_{{ 5}6|E}x_{{4}567|E}x_{E|{3}47}x_{E|{ 2}3467}$, in decreasing order of the $e_i$s.
The two other monomials are  $x_{{ 7}|E}\, x_{{6}7|E}\, x_{E|{5}}\, x_{{4}567|E}\, x_{E|{ 3}5}\, x_{E|{ 2}35}$ and $x_{ 7|E}\, x_{E|{6}}\, x_{{5}7|E}\, x_{{ 4}567|E}\, x_{E|{ 3}6} \, x_{E|{2}3467}$, where the terms are again listed in their order of arrival.
\end{example}


\subsection{The degree of the conormal  fan}
\label{sec:beta(F)}

We now  prove Proposition \ref{prop:degreebetaM}, which 
shows that the degree of $\delta^{n-1}$ in the Chow ring of the conormal fan of $\M$ is the beta invariant of $\M$.
Proposition \ref{prop:degreebetaM} will be used later to obtain the more general Proposition  \ref{prop:betasF}.

We write $\cl$ and $\cl^\perp$ for the closure operators  of $\M$ and $\M^\perp$.
For each basis $B$ of $\M$, denote the corresponding  basis of $\M^\perp$ by $B^\perp \coloneq  E-B$.

\begin{definition}
A \emph{broken circuit} of $\M$ is a set of the form $C - \min C$ where $C$ is a circuit of $\M$. 
\begin{enumerate}[(1)]\itemsep 5pt
\item An \emph{$\nbc$-basis} of $\M$ is a basis of $\M$ that contains no broken circuits of $\M$. 
\item A \emph{$\bnbc$-basis} of $\M$ is an $\nbc$-basis $B$ of $\M$ such that $B^\perp \cup \{0\} - \{1\}$ is an $\nbc$-basis of $\M^\perp$. 
\end{enumerate}
\end{definition}
The number of $\nbc$-bases of $\M$ is the M\"obius number $|\mu_{\M}|$, whereas the number of $\bnbc$-bases of $\M$ is the beta invariant $\beta_{\M}$.
The independence complex $\IN(\M)$ and the reduced broken circuit complex $\BC(\M)$ of  $\M$ are shellable, and hence homotopy equivalent to wedges of spheres. The $\nbc$-bases and $\bnbc$-bases of $\M$ naturally index the spheres in the lexicographic shellings of $\IN(\M)$ and $\BC(\M)$, respectively.
For the $\nbc$ and $\bnbc$ facts stated in this paragraph, we refer to  \cite{Bjorner} and \cite{Zie92}.

\begin{definition}
\label{def:betacones}
For  a $\bnbc$-basis $B$ of $\M$, we define a sequence $e_1,\ldots,e_{n-1}$ by setting
\[
B - 0 = \{e_1 > \cdots > e_{r}\}
 \ \ \text{and} \ \ 
 B^\perp - 1 = \{e_{r+1} < \cdots < e_{n-1}\}.
 \]
We write $\bcone(B)$ for the maximal cone in  $\Sigma_{\M,\M^\perp}$ corresponding to the table 
 \begin{equation*}
\hspace{-.5cm}
\begin{array}{|ccccccccccc|}
\hline
  \cl(e_1) & \subsetneq & \cdots & \subsetneq & \cl(e_1, \ldots, e_{r}) & \blue{\subsetneq} & E & = & \cdots & =  & E  \\
 E & = & \cdots & = & E &
\blue{\supsetneq} & \cl^\perp(e_{r+1}, \ldots, e_{n-1}) & \supsetneq & \cdots & \supsetneq  & \cl^\perp(e_{n-1}) \\
\hline
\end{array} \ .
\end{equation*}
The unique double jump of the displayed biflag is $r$, one less than the rank of $\M$.
\end{definition}

To see that the displayed table  indeed defines a biflag, we verify 
\[
\cl(B-0) \cup \cl^\perp(B^\perp - 1) \neq E.
\]
 Since $B^\perp$ is a basis of $\M^\perp$, we have $1 \notin \cl^\perp(B^\perp - 1)$; and if we had $1 \in \cl(B-0)$, then $B-0 \cup 1$ would contain a circuit $C$ whose minimum element is $1$, and hence $B$ would contain the broken circuit $C-1$, contradicting that $B$ is $\nbc$.

\begin{example}\label{ex:nbc cone}
The matroid of Figure \ref{fig:pyramid} has
three $\bnbc$-bases, namely
\[
B_1 = 0456, \ \ 
B_2 = 0457, \ \ 
B_3 = 0467.
\]
The corresponding maximal biflags are precisely the ones in the expansion of Example \ref{ex:pyramid2}.
\end{example}

We show that this is a general phenomenon.

\begin{proposition}\label{prop:degreebetaM}
Let $\M$ be a loopless and coloopless matroid on $E$. 
Then, in the Chow ring of the conormal fan of $\M$, we have the canonical expansion
\[
\delta^{n-1}=\sum_{B\in \bnbc(\M)} x_{\bcone(B)}
\]
where the sum is over the $\bnbc$-bases of $\M$.
Thus,  the degree of $\delta^{n-1}$ in the Chow ring of the conormal fan of $\M$ is the $\beta$-invariant of $\M$.
\end{proposition}

We proceed with a series of elementary lemmas. Proposition \ref{prop:expansion} describes the canonical expansion of $\delta^{n-1}$ in terms of pairs $(\F|\G, \e)$ of the form
\[
\begin{array}{|cc|ccccccccccc|cc|}
\hline
\varnothing & \subsetneq & F_1 & \subseteq & \cdots & \subseteq & F_d &
\blue{\subsetneq} & F_{d+1} & \subseteq & \cdots & \subseteq  & F_{n-1} & \subseteq & E\\
E & \supseteq & G_1 & \supseteq & \cdots & \supseteq & G_d &
\blue{\supsetneq} & G_{d+1} & \supseteq & \cdots & \supseteq  & G_{n-1} & \supsetneq & \varnothing\\
\hline
&&e_1 && \cdots && e_d && e_{d+1}&& \cdots && e_{n-1}&&\\
\hline
\end{array} \, ,
\]
which have a unique double jump $d$ by Proposition \ref{prop:doublejump}. 
A priori, the double jump could occur at any $d$. 
We will show that, in fact,  the double jump must occur at $d=r$.
In the remainder of this section, we fix a pair $(\F|\G,\e)$ that gives a nonzero summand of the canonical expansion of $\delta^{n-1}$, 
and  write $e_n$ and $e_{n+1}$
for the two elements  of $E$ missing from the sequence $\e$. 

\begin{lemma}\label{lem:ordering}
We consider the jump sets of $\F$ and $\G$ in Definition \ref{jumpsets}.
\begin{enumerate}[(1)]\itemsep 5pt
\item\label{4.9.1}  If $i$ is in the jump set of $\F$ but not in the jump set of $\G$, then $e_i>e_{i+1}$.  
\item\label{4.9.2}  If $i$ is in the jump set of $\G$ but not in the jump set of $\F$,  then $e_i<e_{i+1}$.
\item\label{4.9.3} If $i<j$ and $e_i<e_j$, then $e_i\not\in G_j$. 
\item\label{4.9.4} If $i<j$ and $e_i>e_j$, then $e_j\not\in F_i$.
\end{enumerate}
\end{lemma}

\begin{proof}
We prove the first statement. By way of contradiction, suppose  $i$ is in the jump set of $\F$, not in the jump set of $\G$, and $e_i<e_{i+1}$.
Then, in the canonical expansion of $\delta^{n-1}$, the variable $x_{F_i|G_i}$ arrives after the variable $x_{F_{i+1}|G_{i+1}}$ to the monomial of  $(\F|\G,\e)$. 
It follows that  $e_i$ is not in the intersection $F_{i+1} \cap G_{i+1}$, and this contradicts  
$e_i \in F_i \cap G_i \subseteq F_{i+1} \cap G_i = F_{i+1} \cap G_{i+1}$. 

We prove the third statement.
Suppose that $i<j$ and $e_i<e_j$. 
Then, in the canonical expansion of $\delta^{n-1}$,  
the variable $x_{F_i|G_i}$ arrives  after the variable $x_{F_j|G_j}$ to the monomial of $(\F|\G,\e)$. 
It follows that  $e_i \notin F_j \cap G_j$. Since $e_i \in F_i \subseteq F_j$, we have $e_i \notin G_j$.
\end{proof}

\begin{lemma}\label{lem:d=r}
The unique double jump is at $d=r$, and the table of $(\F|\G,\e)$ is of the form
\[
\begin{array}{|ccccccccccccccc|}
\hline
\varnothing & \subsetneq &
F_1  & \subsetneq & \cdots & \subsetneq & F_r & \blue{\subsetneq} & E &  = & \cdots & =  & E
& = & E \\
E & = &
E & = & \cdots & = & E & \blue{\supsetneq} & G_{r+1} & \supsetneq & \cdots & \supsetneq  & G_{n-1}
& \supsetneq & \varnothing \\
\hline
&&
e_1  & > & \cdots &> & e_r & & e_{r+1}&<& \cdots &<& e_{n-1}
&&
\\
\hline
\end{array}\, .
\]
\end{lemma}

\begin{proof} 
Since $\set{e_1,\ldots,e_d}\subseteq F_d$ and $\set{e_{d+1},\ldots,e_{n-1}}\subseteq G_{d+1}$,  Lemma~\ref{lem:big unions} gives
$F_d\cup G_{d+1}=\set{e_1,\ldots,e_{n-1}}$.
Therefore, the unique nonempty gap $D_d = E-(F_d \cup G_{d+1})$ is equal to $\{e_n, e_{n+1}\}$.

We  next show
$e_1 > e_2 > \cdots > e_d $  and $e_{d+1} < \cdots < e_{n-2} < e_{n-1}$.
By symmetry, it suffices to show the first set of inequalities.
For contradiction, suppose $e_j < e_{j+1}$ for a minimal choice of $j < d$.
If $j > 1$, then $e_{j-1} > e_j$ implies $e_j \notin F_{j-1}$ by Lemma \ref{lem:ordering}.\ref{4.9.4}; if $j=1$ this holds trivially.
On the other hand, $e_j < e_{j+1}$ implies $e_j \notin G_{j+1}$ by Lemma \ref{lem:ordering}.\ref{4.9.3}. However,  we have
\[
 \{e_1, \ldots, e_{j-1}\} \subseteq F_{j-1}, \ \  \{e_{j+1}, \ldots, e_{n-1}\} \subseteq G_{j+1},  \ \ \text{and} \ \  \{e_n, e_{n+1}\} \subseteq G_d \subseteq G_{j+1}. 
 \]
It follows that $F_{j-1} \cup G_{j+1} = E - e_j$, contradicting Lemma \ref{lem:big unions}. 

For $1 \le j <d$, the inequality $e_j > e_{j+1}$ implies  $e_{j+1} \in F_{j+1} - F_j$, and hence $j$ is in the jump set of $\F$. 
 It follows that the jump set of $\F$ is $\{0,  \ldots, d\}$, and similarly, the jump set of $\G$ is $\{d, \ldots, n-1\} $. 
\end{proof}

\begin{lemma}\label{lem:bis}
We have $\set{e_1,\ldots,e_{n-1}}=\set{2,3,\ldots,n}$ and $\set{e_n, e_{n+1}} = \{0,1\}$. 
Moreover,
\begin{align*}
&e_i=\min F_i \ \ \text{and} \ \   F_i = \cl(e_1, \ldots, e_i), \ \ \text{for $1 \le i \le r$},  \\
&e_i=\min G_i  \ \  \text{and} \ \ G_i= \cl^\perp(e_i, \ldots, e_{n-1}) \ \ \text{for $r<i<n$.}
\end{align*}
In particular, the biflag $\F|\G$ and the sequence $\e$  determine each other.
\end{lemma}

\begin{proof}
We may assume without loss of generality that $e_r<e_{r+1}$, so  that the variable $x_{F_r,G_r}$ is the last term to arrive in the monomial corresponding to $(\F|\G,\e)$.
By definition, we have
\[
e_r =\max\Big(E-\bigcup_{\substack{1 \leq j \leq n-1 \\ j \neq r}}
(F_j\cap G_j)\Big) = \max\Big(E-(F_{r-1}\cup G_{r+1})\Big).
\]
If we had $e_r\leq 1$, then $\abs{F_{r-1}\cup G_{r+1}}\geq n-1$, which would imply $\abs{F_r\cup G_{r+1}}=n$, contradicting Lemma~\ref{lem:big unions}. Thus $e_r=2$ and $F_r \cup G_{r+1} = E - \{0,1\}$.

We  now show that $e_i=\min F_i$ for $1 \leq i \leq r$. If this was not the case, then, we would have $\min F_i = e_j < e_i$ for some $j \neq i$, because $0$ and $1$ are not in  $F_i$.
Since $e_1 > \cdots > e_i$, this would imply $i<j$, and Lemma~\ref{lem:ordering}.\ref{4.9.4} would  tell us that $e_j \notin F_i$, contradicting $e_j =\min F_i$.
Similarly, we have $e_i = \min G_i$ for $r < i < n$.

Finally, since $F_i$ has rank $i$ and $ e_i \in F_i - F_{i-1}$ for $1 \le i \le r$ by Lemma \ref{lem:d=r}, the list $e_1, \ldots, e_i$ must be a basis of $F_i$ for $1 \le i \le r$. 
The analogous statement holds for $G_i$ as well.
\end{proof}

\begin{lemma}\label{lem:expansionbetanbc}
The set $\set{0, e_1,\ldots,e_r}$ is a $\bnbc$-basis of $\M$.
\end{lemma}

\begin{proof}
Since $e_r = \min F_r$ by the previous lemma, we have $0 \notin F_r$, and hence $B=\set{0, e_1,\ldots,e_r}$  indeed is a basis of $\M$.
We prove by contradiction that $B$ is $\nbc$. 
Assume that $B$ contains a broken circuit $C- \min C$. Since $\min C$ is not in $B$, one of the following must hold:

\noindent (i) $\min C = 1$. Let $C=\{1, e_{a(1)}, \ldots, e_{a(k)}\}$ where $1 \leq a(1) < \ldots < a(k) \leq r$. Then 
\[
1 \in \cl(e_{a(1)}, \ldots, e_{a(k)}) \subseteq F_{a(k)}\subseteq{F_r}.
\]
 This contradicts Lemma  \ref{lem:bis}, which shows that  the unique nonempty gap is $D_r=\{0,1\}$. 

\noindent (ii) $\min C = e_s$ for  $s \geq r+1$. Let $C=\{e_s, e_{a(1)}, \ldots, e_{a(k)}\}$ where $1 \leq a(1) < \ldots < a(k) \leq r$. Then 
\[
e_s \in \cl(e_{a(1)}, \ldots, e_{a(k)}) \subseteq F_{a(k)}.
\]
 This contradicts Lemma \ref{lem:ordering}.\ref{4.9.4}, since $a(k) \leq r < s$ and $e_{a(k)} > e_s$.

The same argument for $\M^\perp$ shows that $B^\perp - \{0\} \cup \{1\}$ is an $\nbc$-basis of $\M^\perp$. We conclude that $B$ is a $\bnbc$ basis of $\M$, as desired.
\end{proof}

We now have all the ingredients to complete the proof of Proposition \ref{prop:degreebetaM}.

\begin{proof}[Proof of Proposition \ref{prop:degreebetaM}]
Lemma \ref{lem:bis} tells us that each monomial $x_{\F|\G}$ that appears in the canonical expansion of $\delta^{n-1}$ has coefficient $1$. Combined with Lemma \ref{lem:expansionbetanbc}, it also tells us that every term that appears is of the form $x_{\bcone(B)}$ for a $\bnbc$-basis $B$.

Conversely, if $\F|\G$ is the biflag corresponding to the $\beta$-cone of a $\bnbc$-basis $B$, and if we define $\e$ by setting $B = \set{e_1 > \cdots >  e_r > 0}$ and $E-B = \set{e_{n-1} > \cdots > e_{r+1} > 1}$, then $(\F|\G, \e)$ satisfies the conditions of Proposition \ref{prop:expansion}, so it does  arise in the canonical expansion of $\delta^{n-1}$.
\end{proof}

\subsection{A vanishing lemma for the conormal fan}

Throughout the remainder of this section, we fix a  flag of $k$ nonempty proper flats
\[
\F=\{F_1 \subsetneq \cdots \subsetneq F_k\}, \ \ \text{keeping the convention that $F_0=\varnothing$ and $F_{k+1}=E$.}
\]
The interval $\M(i)=[F_{i}, F_{i+1}]$ is said to be \emph{short} if $|F_{i+1} - F_{i}| = 1$ and \emph{long} if $|F_{i+1} - F_{i}| > 1$. 

\begin{definition}
We define the \emph{orthogonal flag} $\F^\perp$  of $\F$ to be the flag of coflats
\[
\F^\perp = \{F_1^\perp \supseteq \cdots \supseteq F_k^\perp\}, \ \ \text{where}  \ \ F_i^\perp = \cl^\perp(E-F_i)  \text{ for } 1 \leq i \leq k.
\]
The orthogonal flag may contain repeated coflats, and it may contain the trivial coflat $E$.\footnote{Strictly speaking, the notation $F^\perp$ conflicts with the notation $B^\perp$ used in Section \ref{sec:beta(F)} for the dual basis of a basis $B$. We trust that no confusion will arise within a given context.}
\end{definition}

Our goal this section is to prove the following lemma, which shows that many monomials in the Chow ring of the conormal fan  vanish when multiplied by the highest possible power of $\delta$.

\begin{lemma}[Vanishing Lemma]\label{lem:vanishing}
Suppose $\F|\G$ is a biflag of length $k$ satisfying the condition
\[
x_{\F|\G} \, \delta^{n-k-1}  \ \ \text{is nonzero in the Chow ring of the conormal fan of $\M$.}
\]
Then $\G$ must be the orthogonal flag $\F^\perp$.
Furthermore,   the interval $\M(i)$ is either short or loopless and coloopless  for all $0 \le i \le k$.
\end{lemma}

Let  $x_{\F^+| \G^+}$ be a nonzero summand of the canonical expansion of $x_{\F|\G} \, \delta^{n-k-1}$,  and let
\[
\F|\G = \F_k|\G_k, \ \, \F_{k+1}|\G_{k+1},\, \ldots, \,  \F_{n-1}|\G_{n-1} = \F^+ | \G^+
\]
be some sequence of biflags obtained by recursively applying Lemma~\ref{lem:canonexpansion} in the  expansion. 
We write
 $D_{i,0} | \cdots | D_{i,i}$ for the gap sequence of $\F_i|\G_i$  in Definition~\ref{def:gaps}.
With Lemma \ref{lem:nongaps} in mind, we set
\begin{equation}\label{eq:Y}
Y_i =\bigsqcup_{j=0}^{i}D_{i,j} = E - \bigcup_{F|G \in \F_i|\G_i} (F \cap G).
\end{equation}
We write
 $D_0 |  \cdots |  D_k$ and $Y$ for the gap sequence and the union of the gaps of the initial flag $\F|\G$.
To prove the Vanishing Lemma \ref{lem:vanishing}, we need a preliminary result.

\begin{lemma} \label{lem:Y}
Suppose that the assumption of Lemma \ref{lem:vanishing} holds for $\F|\G$. 
\begin{enumerate}[(1)]\itemsep 5pt
\item If $\F|\G$ has $m$ empty gaps, then the union of its gaps has size
$|Y| = n+1-m$.
\item
For each empty gap $D_j$, we have $F_{j+1}-F_j = \{e_j\}$ for some $e_j \in E$. Furthermore,  
\[
Y = E - \{e_j \mid D_j = \varnothing\}.
\]
\item
For all $0 \leq i \leq k$, setting  $r_i = \text{rank}_{\M}(F_i)$ and $r^\perp_i = \text{rank}_{\M^\perp}(G_i)$, we have
\[
|F_{i+1}-F_i|  =
(r_{i+1}-r_i) + (r_i^\perp - r_{i+1}^\perp).
\]
\end{enumerate}
\end{lemma}

\begin{proof}[Proof of Lemma \ref{lem:Y}]
Let $m$ be the number of empty gaps of $\F|\G$.
We first prove the inequality
\begin{equation}
\label{eq:Y<}
|Y| \leq n+1-m.
\end{equation}
For each empty gap $D_{j}$, choose an element $e_j\in F_{j+1}-F_j$. Since $e_j\not\in D_{j} = E - (F_j \cup G_{j+1})$, we must have $e_j\in G_{j+1}$. This implies that $e_j \in F_{j+1}\cap G_{j+1}$, so the second equality in \eqref{eq:Y} gives $e_j \notin Y$.  There are $m$ such elements $e_j$, which are all distinct by construction; this implies \eqref{eq:Y<}.

To prove the first statement, it remains to show the opposite inequality
\begin{equation}
\label{eq:Y>}
|Y| \geq n+1-m.
\end{equation}
We obtained $\F_{i+1}|\G_{i+1}$ from $\F_i|\G_i$ by choosing the largest gap element $e = \max Y_i$, finding the unique gap $D_{i,j}$ of $\F_i|\G_i$ containing $e$, and inserting a new pair $F|G$ with $e \in F \cap G$ between the $j$-th and $(j+1)$-th biflats of $\F_i|\G_i$:
\[
\F_{i+1}|\G_{i+1}=\ 
\begin{array}{|ccccccccc|}
\hline
\cdots & \subseteq & F_{i,j} & \subseteq & F & \subseteq & F_{i,j+1} & \subseteq & \cdots \\
\cdots & \supseteq & G_{i,j} & \supseteq & G & \supseteq & G_{i,j+1} & \supseteq & \cdots \\
\hline
\end{array} \, 
\]
Thus the only difference between the gaps of $\F_i|\G_i$ and the gaps of $\F_{i+1}|\G_{i+1}$ is that we are replacing the gap $D_{i,j}$ with two smaller disjoint gaps $D_{i+1,j}$ and $D_{i+1,j+1}$ that do not contain $e$:
\begin{equation} \label{eq:e}
D_{i,j} \supseteq D_{i+1,j} \sqcup D_{i+1,j+1} \sqcup e.
\end{equation}
It is helpful to visualize this data as a graded forest of levels from $k$ to $n-1$. 
The vertices of the bottom level $k$ are the gaps $D_0, \ldots, D_k$ of the original biflag $\F|\G$; they are the roots of the  trees in the forest.
The vertices of the $i$-th level are the gaps 
of $\F_i|\G_i$. To go from level $i$ to level $i+1$, we connect the split gap $D_{i,j}$ with the gaps $D_{i+1, j}$ and $D_{i+1,j+1}$ that replace it.
Every other gap $D_{i,k}$ is connected to the gap in the next level that is equal to it; this is $D_{i+1,k}$ if $k<j$ and $D_{i+1,k+1}$ if $k>j$.
The final biflag $\F^+|\G^+$ has $n$ gaps, of which $n-1$ are empty and one of them, say $D$, has size at least $2$.
Each gap of $\F^+|\G^+$ 
 originates from one of the original gaps of $\F|\G$ through successive gap replacements. 
For each $0 \leq i \leq k$, we set
\[
d_i = \text{number of gaps of } \F^+|\G^+ \text{ that descend from the initial gap } D_{i} \text{ of } \F|\G.
\]
We give an upper bound of $d_i$ in terms of $|D_i|$, in each of the following three cases:

\noindent \emph{Case 1.} $D_{i} = \varnothing$:

\noindent In this case, the gap $D_{i}$ eventually becomes a single empty gap in $\F^+|\G^+$, so $d_i=1$.

\noindent \emph{Case 2.} $D_{i} \neq \varnothing$ is the progenitor of the unique nonempty gap $D$ of $\F^+|\G^+$:

\noindent Consider the gaps that descend from $D_{i}$ throughout the process.
By \eqref{eq:e}, every time one such gap gets replaced by two smaller ones, the size of the union of the gaps strictly decreases. In the end, this union has size $|D| \geq 2$. Therefore these gaps were split at most $|D_{i}|-2$ times, so $d_i \leq |D_{i}|-1$.

\noindent  \emph{Case 3.} $D_{i} \neq \varnothing$ is not the progenitor of the unique nonempty gap $D$ of $\F^+|\G^+$:

\noindent Again, every time a descendant of $D_{i}$ gets replaced by two smaller ones, the size of their union decreases. 
Furthermore, their union can never have size $1$ by Proposition~\ref{prop:gap}. Thus $d_i \leq |D_{i}|$.

\noindent Since the final number of gaps is $n$, we conclude that
\[
n  = \sum_{i=0}^k d_i \leq m + \left(\sum_{i : D_{i}\neq \varnothing} |D_{i}|\right) -1 = m + |Y| -1.
\]
This proves the opposite inequality \eqref{eq:Y>},
and hence
the first statement of the lemma.
Furthermore, every inequality we applied along the way must in fact have been an equality. 
We record these facts:

\noindent
(a)
For \eqref{eq:Y<} to be an equality, we must have $F_{j+1}-F_j = \{e_j\}$ for each empty gap $D_j$, and 
\[
Y = E - \{e_j \mid D_{j} = \varnothing\}.
\]
This proves the second statement of the lemma.

\noindent
(b)
For \eqref{eq:Y>} to be an equality, we must have
\[
\text{$d_i=1$ in \emph{Case 1},  \ \ 
$d_i=|D_{i}|-1$ in \emph{Case 2}, \ \ \text{and} \ \ 
$d_i=|D_{i}|$ in \emph{Case 3}.} 
\]

We use (a) and (b) to prove the third statement of the lemma, in two steps. 
First, we show
\begin{equation}\label{eq:d_l1}
d_i =
\begin{cases}
|F_{i+1}-F_i| & \text{in \emph{Case 1} and \emph{Case 3},} \\
|F_{i+1}-F_i| - 1 & \text{in \emph{Case 2}.}
\end{cases}
\end{equation}
If $D_i$ is empty, then  $d_i=1$ and $|F_{i+1}-F_i|=1$ by (a).
If $D_i$ is nonempty, we claim that
\[
D_{i} = F_{i+1}-F_i.
\]
The forward inclusion holds by definition. For the backward inclusion, let $e$ be an element of $F_{i+1}-F_i$. 
By (a), we must have $e \in Y$, and since $D_{i}$ is the only gap intersecting $F_{i+1}-F_i$, we must have $e \in D_{i}$.
Thus (b) implies \eqref{eq:d_l1}.
Second, we show
\begin{equation}\label{eq:d_l2}
d_i =
\begin{cases}
(r_{i+1}-r_i) + (r_i^\perp - r_{i+1}^\perp)  &  \text{in \emph{Case 1} and \emph{Case 3},} \\
(r_{i+1}-r_i) + (r_i^\perp - r_{i+1}^\perp) -1 &  \text{in \emph{Case 2}.}
\end{cases}
\end{equation}
%
If $D_{i}$ is not the progenitor of $D$, then
 the part of $\F^+|\G^+$ between $F_i|G_i$ and $F_{i+1}|G_{i+1}$ contains no double jumps. In each of the $d_i$ single jumps, either the rank increases by 1 or the corank decreases by 1, but not both. Therefore $d_i$ must equal the sum of the rank increase $r_{i+1}-r_i$ and the corank decrease $r^\perp_i - r^\perp_{i+1}$.
If 
%
$D_{i}$ is the progenitor of $D$, then
the part of $\F^+|\G^+$ between $F_i|G_i$ and $F_{i+1}|G_{i+1}$ contains one double jump.
 In each of the $d_i-1$ single jumps, either the rank increases by 1 or the corank decreases by 1, but not both. In the double jump, both changes occur.
 Therefore $d_i+1$ must equal the sum of the rank increase $r_{i+1}-r_i$ and the corank decrease $r^\perp_i - r^\perp_{i+1}$.
 Combining  \eqref{eq:d_l1} and \eqref{eq:d_l2}, we get the third statement of the lemma.
\end{proof}


\begin{proof}[Proof of the Vanishing Lemma \ref{lem:vanishing}]
First, we prove that $\G$ must be the orthogonal flag $\F^\perp$. 
We write $\rank$ and $\rank^\perp$ for the rank functions of $\M$ and $\M^\perp$.
One readily verifies that 
\[
(\rank(F_{i+1}) - \rank (F_i)) +
(\rank^\perp(E-F_i) -  \rank^\perp(E-F_{i+1})) =
|F_{i+1} - F_i|
\qquad \text{for $0 \leq i \leq k$. }
\]
By the third statement of Lemma \ref{lem:Y}, the sequences $\rank^\perp(E-F_i)$ and $\rank^\perp(G_i)$ satisfy the same recurrence; they also have the same initial value, so 
\[
\rank^\perp(E-F_i) =
\rank^\perp(G_i)
\quad \text{for $0 \leq i \leq k$. }
\]
Now $F_i \cup G_i = E$ implies $G_i \supseteq E-F_i$. 
Since $G_i$ is a coflat, we have  $G_i \supseteq \cl^\perp(E-F_i) = F_i^\perp$. 
 It follows that $G_i \supseteq  F_i^\perp$ are flats of $\M^\perp$ of the same rank, and hence $G_i = F_i^\perp$ for all $i$.

Next, we prove that every long interval $\M(i)$ must be loopless and coloopless. We argue by contradiction.

First assume that $\M(i) = (\M/F_{i})|(F_{i+1} - F_{i})$ has a loop $l$. Since restriction cannot create new loops, the element $l$ must also be a loop of $\M/F_{i}$. This contradicts the fact that $F_{i}$ is a flat.

Now assume that $\M(i) = (\M|F_{i+1})/F_{i}$ has a coloop $c$. Since contraction cannot create new coloops, the element  $c$ must also be a coloop of $\M|F_{i+1}$. Thus
$\rank(F_{i+1}-c) = \rank(F_{i+1})-1$,
which implies that
$\rank^\perp((E-F_{i+1}) \cup c) = \rank^\perp(E-F_{i+1})$.
This means that $c \in \cl^\perp(E-F_{i+1}) = F_{i+1}^\perp$.
Now, since $\M(i)$ is long, the second statement of Lemma \ref{lem:Y}  implies that $D_{i+1}$ is nonempty and that $c \in Y$. But then we must have $c \in D_{i+1} = (F_{i+1} - F_{i}) \cap (F_{i}^\perp - F_{i+1}^\perp)$, contradicting that $c \in F_{i+1}^\perp$. 
\end{proof}

\subsection{The beta invariant of a flag and the conormal intersection theory}\label{flagbeta}

In this section, we complete the proof that the degree of $\pi^*(x_\F)\delta^{n-k-1}$ is equal to the $\beta$-invariant $\beta_{\M[\F]}$. 
To prove by induction, we need a lemma relating the conormal fan of $\M$ with that of the contraction $\M/i$,
where $i$ is any fixed element  of $E$ that has no parallel elements in $\M$.

We may assume  that $\M$ has no loops and no coloops. 
Thus, we have $i^\perp$ is the ground set $E$ and $i|E$ is a biflat of $\M$.
We consider the simplicial fan
\[
\st_{i|E}\Sigma_{\M, \M^\perp} \subseteq  (\N_E / \e_i) \oplus \N_E.
\]
 We write  $\ebar_S$ for the image of $\e_S$ in $\N_E/\e_i$, and $\xbar_{F|G}$ for the variable in the Chow ring of the star corresponding to a biflat $F|G$; we set it equal to $0$ if $F|G$ does not correspond to a ray in this star.

\begin{lemma} \label{lem:contracti}
Consider the natural projection $\psi\colon (\N_E / \e_i) \oplus \N_E \longrightarrow \N_{E-i} \oplus \N_{E-i}$.
\begin{enumerate}[(1)]\itemsep 5pt
\item
 The projection $\psi$ induces a morphism of fans
\[
\psi\colon \st_{i|E}\Sigma_{\M, \M^\perp} \longrightarrow \Sigma_{\M/i, (\M/i)^\perp}.
\]
\item\label{4.17.2}
 The pullback $\psi^*$   between the Chow rings 
  is given by
\[
\psi^*(x_{P|Q}) = \xbar_{(P\cup i)|Q} + \xbar_{(P\cup i)|(Q \cup i)}, 
\]
where at least one of the terms in the right-hand side is nonzero.
\item\label{4.17.3}
 For any element $j$ of $E$, the pullback $\psi^*$ maps the class $\delta$ to the class 
 \[
\overline{\delta}=\overline{\delta}_j\coloneq  
 \sum_{ i \in F, \,  j \in F \cap G} \xbar_{F|G}. 
\]
\item\label{4.17.4}
 The pullback $\psi^*$ commutes with the degree maps of  the star and the conormal fan $\Sigma_{\M/i, (\M/i)^\perp}$: 
 \[
\deg_{\M/i} x_{\P|\Q} = 
  \deg_{\M}  (x_{i|E} \,\psi^*(x_{\P|\Q})).
\]
 \end{enumerate}
\end{lemma}

\begin{proof}
Let $F|G$ be a biflat of $\M$ with $i \in F$. The image of a ray corresponding to $F|G$ in the star is
\[
\psi(\ebar_F + \f_G) = \e_{F-i} + \f_{G-i}, 
\]
which is a ray of the conormal fan $\Sigma_{\M/i, (\M/i)^\perp}$  because $(F-i)|(G-i)$ is a biflat of $\M/i$:
\begin{align*}
\cl_{\M/i}(F-i) &= \cl_{\M}(F) - i = F - i, \text{ and}  \\
\quad \cl_{\M^\perp - i}(G-i) &= \cl_{\M^\perp}(G-i) - i \subseteq  \cl_{\M^\perp}(G) - i = G-i.
\end{align*}
Furthermore, if $i|E \cup \F|\G$ is a biflag of $\M$, its gaps occur to the right of $i|E$, and there will also be gaps in the corresponding positions of the biflag
\[
(\F-i)|(\G-i) \coloneq  \Big\{(F-i)|(G-i) \Big\}_{F|G \in \F|\G} \ \ \text{of $\M/i$}.
\]
 Therefore, the projection $\psi$ maps cones to cones.
This proves the first statement.

The value of the piecewise linear function $\psi^*x_{P|Q}$ on a ray $\ebar_F + \f_G$ of the star is
\[
\psi^*x_{P|Q}(\ebar_F + \f_G)
= x_{P|Q}(\e_{F-i} + \f_{G-i}))
= \begin{cases}
1 & \text{if $F=P \cup i$  and $G \in \{Q, Q \cup i\}$},  \\
0 & \text{if otherwise}.
\end{cases}
\]
Since $Q$ is a flat of $\M/i$, at least one of $Q$ and $Q \cup i$ is a flat of $\M$,
and we have the second statement.
Given the second statement, the third statement is a straightforward computation
%
\[
\psi^*(\delta_j)
= \sum_{ j \in P \cap Q} \Big(\xbar_{(P \cup i)|Q} + \xbar_{(P \cup i)| (Q \cup i)}\Big)
= \sum_{i \in F, \,  j \in F \cap G} \xbar_{F|G} = \overline{\delta}_j,
\]
where the first sum is over the biflats $P|Q$ of $\M/i$ and the second sum over the biflats $F|G$ of $\M$.

For the last  statement, we need to verify that, for any maximal biflag $\P|\Q$ of $\M/i$,
\[
\deg_{\M/i} x_{\P|\Q} =   \deg_{\M}  (x_{i|E} \,\psi^*(x_{\P|\Q})).
\]
Applying the second statement to
each $P|Q$ in $\P|\Q$, we may express $x_{i|E} \,\psi^*(x_{\P|\Q})$ as a sum of square-free monomials. One of the terms in this expression is $x_{i|E} \, x_{(\P \cup i)|\cl^\perp(\Q)}$, where 
\[
(\P \cup i)|\cl^\perp (\Q) \coloneq \Big\{(P \cup i)|\cl^\perp(Q)\Big\}_{ P|Q \in \P|\Q}.
\]
 We need to prove  that this is the only nonzero term.

Consider any term $x_{i|E} x_{\F|\G}$ that arises in the expression for $x_{i|E} \,\psi^*(x_{\P|\Q})$. We automatically have $F_j = P_j \cup i$ for all $j$, so it remains to prove that $G_j$ is the closure of $Q_j$ in $\M^\perp$ for all $j$.
Let $k$ be the largest index satisfying $i \in \cl^\perp(Q_k)$, so that 
\[
 \cl^\perp(Q_j) = Q_j \cup i \ \  \text{for $j \leq k$} \ \ \text{and} \ \  \cl^\perp(Q_j)=Q_j \ \  \text{for $j > k$}.
 \]
For $ j \leq k$, the set $Q_j$ is not a flat in $\M^\perp$, and hence $G_j = Q_j \cup i = \cl^\perp(Q_j)$.
Since $Q_k$ and $Q_{k+1}$ are flats of consecutive ranks in $(\M/i)^\perp = \M^\perp - i$,  the flats $Q_k \cup i =\cl^\perp(Q_k)$ and $Q_{k+1} = \cl^\perp(Q_{k+1})$ of $\M^\perp$ also have consecutive ranks. Since $Q_{k+1} \cup i$ is strictly between the two, it cannot be a flat of $\M^\perp$.
Thus,  we must have $G_{k+1} = Q_{k+1}$, and hence $G_j = Q_j = \cl^\perp(Q_j)$ for $j > k$. We conclude that $\F|\G = (\P \cup i) | \cl^\perp(\Q)$ as desired.
\end{proof}

We can now give an intersection-theoretic interpretation of the beta invariant of a flag.
Together with the Vanishing Lemma \ref{lem:vanishing}, it gives the identity $\deg(\pi^* (x_{\F}) \, \delta^{n-k-1}) = \beta_{\M[\F]}$.

\begin{proposition}\label{prop:betas}
For any strictly increasing flag $\F$ of 
$k$ nonempty proper flats of $\M$, we have
\[
\deg(x_{\F|\F^\perp} \, \delta^{n-k-1}) = \beta_{\M[\F]}.
\]
\end{proposition}

\begin{proof}
We proceed by induction on $k$. The base case $k=0$ follows from Proposition \ref{prop:degreebetaM}:
\[
\deg(\delta^{n-1})=\deg\Big(\sum_{B\in \bnbc(\M)} x_{\bcone(B)} \Big)=\beta_{\M}.
\]
When $k$ is positive,
write $F|F^\perp$ for the first biflat in $\F|\F^\perp$, and write $\F|\F^\perp$ as the disjoint union $F|F^\perp \,  \sqcup \, \G|\G^\perp$. 
We consider the contraction $\M/F$, its flag  of flats $\G-F \coloneq  \{G-F  \}_{ G \in \G}$, and 
the corresponding biflag 
\[
(\G-F)|(\G-F)^\perp  \coloneq  \Big\{(G-F)|(G^\perp - F)  \Big\}_{G \in \G}. 
\]
The displayed description is justified by 
$G^\perp - F = \cl_{(\M/F)^\perp}((E-F) - (G-F))$. 
Clearly, 
\[
\beta_{\M[\F]} = \beta_{\M|F} \cdot \beta_{(\M/F)[\G-F]}.
\]
We separately consider two cases, depending on the shortness of the first interval $[\varnothing,F]$:

\noindent \emph{Case 1.} The flat $F$ contains exactly one element $i \in E$.

Recall that the beta invariant of $[\varnothing,i]$ is equal to $1$. Therefore, $\beta_{\M[\F]}$ is equal to
\[
\setstretch{1.5}
\begin{array}{lcll}
 \beta_{(\M/i)[\G-i]} &=& \deg_{\M/i}\big(x_{\G-i |(\G-i)^\perp}\, \delta_{\M/i}^{(n-1)-(k-1)-1}\big)&\text{by the inductive hypothesis,}   \\
&=& \deg_{\M}\big(x_{i|E}\, \psi^*(x_{\G-i | (\G-i)^\perp}) \, \delta_{\M}^{n-k-1})\big) &\text{by Lemma \ref{lem:contracti}.\ref{4.17.3}  and \ref{lem:contracti}.\ref{4.17.4},}   \\
&=& \displaystyle \deg_{\M}\big(x_{i|E}\prod_{G \in \G}(x_{G | (G^\perp-i)} + x_{G | G^\perp}) \,  {\delta}_{\M}^{n-k-1}\big) \qquad  & \text{by Lemma \ref{lem:contracti}.\ref{4.17.2}.} 
\end{array}
\]
By  the Vanishing Lemma \ref{lem:vanishing}, the right-hand side simplifies to 
\[
 \deg_{\M}\big(x_{i|E} \, x_{\G | \G^\perp} \,  {\delta}_{\M}^{n-k-1}\big) =\deg_{\M}\big(x_{\F|\F^\perp} \, {\delta}_{\M}^{n-k-1}\big).
 \]

\noindent \emph{Case 2.} The flat $F$ contains more than one element.

%

We may assume the interval $[\varnothing, F]$ is coloopless by  Lemma \ref{lem:vanishing}. This means that the flat $F$ is \emph{cyclic}, that is,   $F^\perp = E-F$. We then have the natural bijections
\begin{align*}
\phi_1\colon
\Big\{\text{biflats of $\M|F$}\Big\}
&\longrightarrow 
\Big\{\text{biflats $F'|G'$ of $\M$ with $F' \subseteq F$ and $G' \supseteq E-F$}\Big\} \ \ \text{and} \\
\phi_2\colon
\Big\{\text{biflats of $\M/F$}\Big\}
&\longrightarrow 
\Big\{\text{biflats $F'|G'$ of $\M$ with $F' \supseteq F$ and $G' \subseteq E-F$}\Big\},
\end{align*}
where $\phi_1(A| B) = A|(B \cup (E-F))$ and $\phi_2(A|B) = (A \cup F)|B$.
The bijection $\phi_1$ extends to the bijection between the biflags of $\M|F$ and the biflags of $\M$ that are supported on
the corresponding set of biflats, and have a gap to the left of $F|(E-F)$.
Similarly, the bijection
 $\phi_2$ extends to the bijection between the biflags of $\M/F$ and the biflags of $\M$ that are supported on
the corresponding set of biflats, and have a gap to the right of $F|(E-F)$.

We now compute the degree of $x_{\F|\F^\perp} \, \delta^{n-k-1}$ using the following variant of the canonical expansion of Definition \ref{def:canonical}, which proceeds in two stages:

\noindent \emph{Stage 1}. At each step, choose $e$ to be the largest gap element that is in $F$, if there is one.

\noindent \emph{Stage 2}. At each step, choose $e$ to be the largest gap element in $E-F$.

The first $|F|-2$ steps of this computation will give $x_{\F|\F^\perp}$ times the image under $\phi_1$ of the canonical expansion of $\delta_{M|F}^{|F|-2}$. By Proposition \ref{prop:degreebetaM}, there will be $\beta_{\M|F}$ square-free monomials.

Each such monomial will have a unique nonempty gap before $F$; say it is $D_j$, between biflats $F_j|G_j$ and $F_{j+1}|G_{j+1}$ of $\M$.
The flats $F_j$ and $F_{j+1}$ have consecutive ranks,
and the coflats  $G_j$ and $G_{j+1}$ have consecutive  coranks. In step $|F|-1$ of the computation, this gap $D_j$ will be filled in a unique way by the biflat $F_{j+1}|G_j$. There will no longer be gap elements in $F$. 

In step $|F|$, the computation will enter \emph{Stage 2} for each of the resulting $\beta_{\M|F}$ monomials. The following $(|E-F|-1) - (k-1) - 1$ steps will compute the image under $\phi_2$ of the canonical expansion of $x_{(\G-F) | (\G-F)^\perp} \delta_{M/F}^{|F|-2}$. This expansion has $\beta_{(\M/F)[\G-F]}$ square-free monomials, by the inductive hypothesis.

This concludes the computation of $x_{\F|\F^\perp} \delta^{n-k-1}$. The result will be the sum of $ \beta_{\M[\F]}$ square-free monomials, as we wished to prove.
\end{proof}

\begin{proposition}\label{prop:betasF}
Let $\F=\{F_1 \subsetneq \cdots \subsetneq F_k\}$ be a strictly increasing flag of
flats of $\M$.  We have
\[
\deg( \pi^*(x_{\F}) \, \delta^{n-k-1}) = \beta_{\M[\F]}.
\]
\end{proposition}

\begin{proof}
Since $\pi^*(x_{\F}) = \sum_{\F | \G \text{ biflag}} x_{\F|\G}$, this follows from Lemma \ref{lem:vanishing} and Proposition \ref{prop:betas}.
\end{proof}

\subsection{A conormal interpretation of the Chern--Schwartz--MacPherson cycles}
\label{sec:CSM}

Recall that
the $k$-dimensional Chern--Schwartz--MacPherson cycle of $\M$ is the Minkowski
weight
$\csm_k(\M)$ on the Bergman fan of $\M$ defined by the
formula
\[
\csm_k(\M)(\sigma_{\F})= (-1)^{r-k}  \beta_{\M[\F]},
\]
where $\sigma_{\F}$ is the  $k$-dimensional cone corresponding to  a  flag  of flats $\F$ of $\M$.



\begin{reptheorem}{CSMTheorem}
When $\M$ has no loops and no coloops, we have 
\[
\csm_k(\M)=(-1)^{r-k}\pi_*(\delta^{n-k-1}\cap 1_{\M,\M^\perp}) \ \ \text{for $0 \le k \le r$.}
\]
\end{reptheorem}


\begin{proof}
 By  Proposition \ref{prop:betasF}, Definition \ref{def:deg}, and the projection formula, we have
\[
\beta_{\M[\F]} = \deg(\pi^*(x_{\F})\delta^{n-k-1}) 
=\pi^*(x_{\F})\delta^{n-k-1}\cap 1_{\M,\M^\perp}
=\pi_*(\delta^{n-k-1}\cap 1_{\M,\M^\perp})(\sigma_{\F}).
\]
The result then follows from the definition of the Chern--Schwartz-MacPherson cycle of $\M$.
\end{proof}

\begin{reptheorem}{degtheorem}
When $\M$ has no loops and no coloops, we have
\[
\overline{\chi}_{\M}(q+1)=\sum_{k=0}^r (-1)^{r-k} \deg( \gamma^k\, \delta^{n-k-1}) \, q^k.
\]
\end{reptheorem}

\begin{proof}
We use \cite[Theorem 1.4]{LRS17}, which states that
\[
\overline{\chi}_{\M}(q+1)=\sum_{k=0}^r \alpha^k \cap \csm_k(\M) \, q^k.
\]
The authors of \cite{AB20} give
a non-recursive proof of the identity using tropical intersection theory.
For representable matroids, the identity was given earlier in 
 \cite[Theorem 1.2]{Alu13}.
By Theorem \ref{CSMTheorem} and the projection formula, the $k$-th coefficient of the displayed polynomial is
\[
 \alpha^k \cap \pi_*(\delta^{n-k-1} \cap 1_{\M,\M^\perp})
 = \pi_* \big(\pi^*\alpha^k \cap (\delta^{n-k-1} \cap 1_{\M,\M^\perp}) \big) 
 = \pi_* \big(\gamma^k \delta^{n-k-1} \cap 1_{\M,\M^\perp} \big).
\]
This proves the desired formula for the reduced characteristic polynomial.
%
%
\end{proof}



\section{Tropical Hodge theory}\label{sec:hodge}



Throughout this section, we fix a rational simplicial fan $\Sigma$ in  $\N=\mathbb{R} \otimes \N_\mathbb{Z}$. 
Our goal is to prove Theorem \ref{thm:lefschetz}, which says that the property of $\Sigma$ being Lefschetz only depends on the support of $\Sigma$.
We deduce the Lefschetz property of $\Sigma_{\M,\M^\perp}$ from the Lefschetz property of $\Sigma_{\M} \times \Sigma_{\M^\perp}$,
and use it to prove Theorem \ref{thm:conjs}.

\subsection{Convexity of  piecewise linear functions}\label{ss:K}
A piecewise linear
function $\phi\colon\Sigma\to\R$ is said to be \emph{positive} on $\Sigma$ if
$\phi(x)$ is positive for all nonzero $x\in\abs{\Sigma}$. 
A class in $A^1(\Sigma)$ is said to be \emph{positive} if it has a positive representative.
We write $\Eff^\circ(\Sigma)\subseteq A^1(\Sigma)$ for the open cone of positive classes.

For each cone $\sigma$ of $\Sigma$, the projection
$ \N\to \N/\spn(\sigma)$ defines a morphism from the closed star 
\[
\pi_\sigma\colon \clst_{\Sigma}(\sigma) \longrightarrow \st_{\tallSigma}(\sigma).
\]
It is straightforward to check that  the associated pullback map $\pi_\sigma^*$ between their Chow rings is an isomorphism, and that for a ray $\nu$ of
$\st_{\Sigma}(\sigma)$, we have
\[
\pi_\sigma^*(x_\nu)=\frac{\mult(\sigma\cup\set{\nu})}{\mult(\sigma)} x_\nu.
\]
Thus, the inclusion of fans $i_\sigma\colon \clst_{\Sigma}(\sigma) \to \Sigma$ defines a ring homomorphism
\[
i_\sigma^*\colon A(\Sigma)\longrightarrow  A(\clst_{\Sigma}(\sigma)) \simeq A(\st_{\tallSigma}(\sigma)),
\]
where the first factor is given by the restriction of piecewise linear functions and the second factor is the inverse of $\pi^*_\sigma$.

We use the pullback homomorphism $i^*_\sigma$ to define \emph{strict convexity} of piecewise linear functions  on $\Sigma$.
The notion agrees with the one used in \cite[Section 4]{AHK15}.

\begin{definition}\label{def:KSigma}
 The cone 
$\Kone(\Sigma)\subseteq A^1(\Sigma)$ is defined by the following conditions:
\begin{enumerate}[(1)]\itemsep 5pt
\item If $\Sigma$ is at most $1$-dimensional,  $\Kone(\Sigma)=
\Eff^\circ(\Sigma)$.  
\item If otherwise, 
$
\Kone(\Sigma) =  \set{\ell\in A^1(\Sigma)\colon \ell\in \Eff^\circ(\Sigma)
\text{~and~} i_\sigma^*(\ell)\in \Kone(\st_{\tallSigma}(\sigma))
\text{~for all nonzero $\sigma\in \Sigma$}}$.
\end{enumerate}
The piecewise linear functions on $\Sigma$ whose classes are in $\Kone(\Sigma)$ are said to be \emph{strictly convex}.
\end{definition}

Clearly, $\ell$ belongs to $\Kone(\Sigma)$ if and only if $i_\sigma^*(\ell)$ belongs to
 $\Eff^\circ(\st_{\tallSigma}(\sigma))$ for all $\sigma\in \Sigma$.\footnote{When $\Kone(\Sigma)$ is nonempty, its closure in $A^1(\Sigma)$ is the cone $\mathscr{L}(\Sigma)$ introduced in \cite[Definition 2.5]{GM12}. 
The cone $\mathscr{L}(\Sigma)$ consists of divisor classes on the toric variety $X_\Sigma$ of $\Sigma$ whose pullback to any  torus orbit closure is effective.
When $\Sigma$ is complete, $\Kone(\Sigma)$ is the ample cone of $X_\Sigma$ \cite[Theorem 6.1.14]{CLS}.} 
 Geometrically,
 $\ell$ belongs to $\Kone(\Sigma)$ if and only if, for each
cone $\sigma$, the class $\ell$ has a piecewise linear representative  
which is zero
on $\sigma$ and positive on the cones containing $\sigma$.
When $\Sigma$ has convex support of full dimension, the notion coincides with the usual notion of strict convexity of piecewise linear functions \cite[Section 6.1]{CLS}. 
In general,  $\Kone(\Sigma)$ is an open polyhedral cone, and $i^*_\sigma
\Kone(\Sigma)\subseteq \Kone(\st_{\tallSigma}(\sigma))$ for all $\sigma\in \Sigma$.


\begin{remark}
A fan $\Sigma$ is \emph{quasiprojective} if it is a subfan of the normal fan
of a convex polytope. 
 When $\Sigma$ is quasiprojective, the cone   $\Kone(\Sigma)$ is nonempty.
More generally, for  simplicial fans  $\Sigma_1 \subseteq \Sigma_2$,
the restriction of piecewise linear functions maps $\Kone(\Sigma_2)$ to $\Kone(\Sigma_1)$.
\end{remark}


A map of fans $f\colon \Sigma\to\Sigma'$ is said to be \emph{projective} if
the induced map of toric varieties $X_\Sigma\to X_{\Sigma'}$ is projective
in the sense of Grothendieck \cite[D\'efinition 5.5.2]{EGAII}.  
According to  \cite[Theorem\ 7.2.12]{CLS}, the map $f$ is projective if and only if 
 $f$ is proper and there exists a piecewise linear function
$\eta$ on $\Sigma$ for which $\eta$ is strictly convex on 
$\abs{f^{-1}(\sigma')}$ for each cone $\sigma'$ of $\Sigma'$,
where $f^{-1}(\sigma')$ denotes the subfan of $\Sigma$ consisting
of cones mapping into $\sigma'$ under $f$.  If $f$ is
induced by the identity map $\N\to \N$, then $f$ is proper if and only if
$\Sigma$ subdivides $\Sigma'$.  In this case, we will call
$\Sigma$ a \emph{projective refinement} of $\Sigma'$ if $f$ is moreover
projective.

\begin{proposition}\label{prop:proj morphism}
  Let $f\colon \Sigma\to\Sigma'$ be a projective map of simplicial  fans.
  If $\Kone(\Sigma')$ is nonempty, then $\Kone(\Sigma)$ is nonempty.
\end{proposition}

\begin{proof}
  We proceed by induction on the dimension of $\Sigma$.  If $\dim\Sigma=1$,
  then $\Kone(\Sigma)\neq\emptyset$.  Otherwise, we choose any
  $\ell\in\Kone(\Sigma')$, and 
%
  let $\eta$ be a piecewise linear  function given by the projectivity of $f$.

  First, since $\ell$ is strictly convex, $\ell$ has
  a representative which is positive
  on $\abs{\Sigma'}-\set{0}$.  Thus $f^*\ell$ is nonnegative on $\abs{\Sigma}$,
  and positive outside of $\abs{f^{-1}(0)}$.  Modulo a global linear function,
  we may choose $\eta$ to be positive on $\abs{f^{-1}(0)}-\set{0}$.  It
  follows that
  $f^*\ell+\epsilon_0\cdot\eta$ is positive on $\abs{\Sigma}-\set{0}$ for
  sufficiently small $\epsilon_0>0$.

  For a nonzero cone $\sigma$, 
  let $\sigma'$ be the smallest cone of $\Sigma'$ containing $f(\sigma)$.
  A fortiori, the restriction
  $f\colon\st_{\Sigma}(\sigma)\to\st_{\Sigma'}(\sigma')$ is projective, 
  and $\dim\st_{\Sigma}(\sigma)<\dim\Sigma$.
  Since $\ell$ is strictly convex, we have 
  $i^*_{\sigma'}\ell\in \Kone(\st_{\Sigma'}(\sigma'))$.  Therefore
  the restriction of $f^*\ell+\epsilon_{\sigma}\cdot\eta$ is in
  $\Kone(\st_{\Sigma}(\sigma))$ for sufficiently small $\epsilon_\sigma$, by
  the induction hypothesis.  We conclude that $f^*\ell+\epsilon\cdot\eta\in
  \Kone(\Sigma)$ for all positive $\epsilon\leq\min\set{\epsilon_{\sigma}}$.
\end{proof}

We now focus on how subdividing a cone by adding a ray affects the
Chow ring.  For $\sigma \in \Sigma$, let $\rho$ denote a new ray
spanned by a primitive vector
\[
  \e_\rho \coloneq \sum_{\nu\in\sigma(1)}a_\nu\mathbf{e}_\nu,
  \]
for some positive rational coefficients $\set{a_\nu}$.
The \emph{stellar subdivision} of $\Sigma$ by $\rho$, denoted
$\stellar_{\rho}\Sigma$,  is obtained from $\Sigma$
by setting
\[
\stellar_{\rho}\Sigma \coloneq
\big(\Sigma - \set{\tau\in \Sigma\colon \tau\supseteq \sigma}\big)
\cup \big(\rho+\partial\sigma+\link_{\Sigma}(\sigma)\big),
\]
where the right-hand $+$ is an internal direct sum of fans.%
\footnote{If $\Sigma\subseteq \N$ and $\Sigma'\subseteq \N'$ are fans for
  which $\N\cap \N'=\set{0}$, the internal direct sum, by definition,
  consists of cones $\sigma+\sigma'$ for all $\sigma\in \Sigma$ and $\sigma'\in \Sigma'$,
  where $+$ denotes Minkowski sum.}
  In the special case when $X_\Sigma$ is smooth and each $a_\nu=1$, 
  the toric variety of the stellar subdivision is the blowup of
  $X_\Sigma$ along the torus orbit closure $V(\sigma)$ \cite[Section 3.3]{CLS}.
In the remainder of this section, we write $\TSigma$ for the stellar subdivision of $\Sigma$ by $\rho$, and write $\bl\colon\TSigma\to \Sigma$ for the map of fans given by the
identity map of $\N$.

Any stellar subdivision is projective. 
In fact, the function $\eta=-x_\rho$ is strictly convex on the
preimage of each cone of $\Sigma$ \cite[Proposition 11.1.6]{CLS}.
The proof of Proposition~\ref{prop:proj morphism} then shows the following.

\begin{proposition}\label{lem:blowup K}
If $\ell\in\bl^*(\Kone(\Sigma))$, then $\ell-\epsilon\cdot x_\rho
  \in \Kone(\TSigma)$ for sufficiently small $\epsilon>0$.
\end{proposition}


%
%

We will  distinguish two cases of stellar subdivisions. 
 The first is the case when 
every closed orbit in $X_{\Sigma}$ meets $V(\sigma)$.  In terms of
fans,
this means that $\Sigma$ is the closed star of $\sigma$
and $\TSigma$ is the closed star of the new ray $\rho$.
In this case, we have 
\[
A(\Sigma)\simeq A(\st_\tallSigma(\sigma)) \ \ \text{and} \ \ 
A(\TSigma)\simeq A(\st_{\TSigma}(\rho)),
\]
which are Chow rings of fans of dimensions $\dim(\Sigma)-\dim(\sigma)$ and
$\dim(\Sigma)-1$.  We will call this a \emph{star-shaped subdivision}.
If otherwise, we will call the stellar subdivision  \emph{ordinary}.

\begin{remark}
In general, the quotient map $\N/\spn(\rho)\to \N/\spn(\sigma)$
induces a map between the stars $\st_{\TSigma}(\rho)\to \st_{\tallSigma}(\sigma)$, and 
the corresponding map of toric varieties is a projective bundle. 
If the stellar subdivision is star-shaped, then $\Sigma$ and $\TSigma$ cannot be Lefschetz,
as their Chow rings 
 vanish in degree $\dim \Sigma$.
 However, the smaller-dimensional fans 
 $\st_{\tallSigma}(\sigma)$ and $\st_{\TSigma}(\rho)$, whose Chow rings are isomorphic to the Chow rings of $\Sigma$ and $\TSigma$ respectively, 
 can be Lefschetz.
\end{remark}

 In the star-shaped case, we will freely use the isomorphisms $A(\Sigma)\simeq A(\st_\tallSigma(\sigma))$ and
$A(\TSigma)\simeq A(\st_{\TSigma}(\rho))$
 in the arguments that follow. 
 This allows us to think of  the bundle map $\st_{\TSigma}(\rho)\to \st_{\tallSigma}(\sigma)$ as the stellar subdivision $\TSigma \to \Sigma$.

\subsection{Lefschetz fans}

Recall from Definition \ref{def:lefschetz} that a  $d$-dimensional Lefschetz fan $\Sigma$
 has 
a $d$-dimensional fundamental weight $w$ which induces \emph{Poincar\'e
duality}. We shall abbreviate this statement  by $\PD(\Sigma)$.  
A Lefschetz fan also satisfies the \emph{hard Lefschetz
property} and \emph{Hodge--Riemann relations} in Definition \ref{def:lefschetz}.  
We will call
these statements $\HL^k(\Sigma,\ell)$ and $\HR^k(\Sigma,\ell)$, respectively, for $0\leq k\leq \frac{d}{2}$ and  $\ell\in \Kone(\Sigma)$.
We say that  $\HL^k(\Sigma)$ holds if $\HL^k(\Sigma,\ell)$ holds for all
$\ell\in \Kone(\Sigma)$, and that $\HL(\Sigma)$ holds if  
$\HL^k(\Sigma)$ holds  for all $k$.  We will use the symbols $\HR^k(\Sigma)$ and $\HR^k(\Sigma,\ell)$ analogously.


\begin{definition}\label{def:mixed lefschetz}
Let $\Sigma$ be a rational simplicial fan. 
\begin{enumerate}[(1)]\itemsep 5pt
\item The fan $\Sigma$ satisfies the \emph{mixed hard Lefschetz
  property} if, for  $0\leq k\leq \frac{d}{2}$ and
  all $\ell_1,\ldots,\ell_{d-2k}\in \Kone(\Sigma)$, the multiplication map
\[
A^k(\Sigma) \longrightarrow A^{d-k}(\Sigma),
\qquad
\eta \longmapsto \Bigg(\prod_{i=1}^{d-2k}\ell_{i} \Bigg)\eta
\]
is a linear isomorphism.  
\item The fan $\Sigma$ satisfies the \emph{mixed Hodge--Riemann
  relations} if, for all $0\leq k\leq \frac{d}{2}$ and  $\ell_0,\ell_1,\ldots,\ell_{d-2k}\in \Kone(\Sigma)$,  
  the bilinear form 
 \[
A^k(\Sigma) \times A^k(\Sigma) \longrightarrow \R, \qquad   (\eta_1,\eta_2) \longmapsto (-1)^k\deg  \Bigg(\prod_{i=1}^{d-2k}\ell_{i} \Bigg)\eta_1 \eta_2
\]
    is positive definite when restricted to the kernel of the multiplication map $ \prod_{i=0}^{d-2k}\ell_{i}$.
\end{enumerate}
\end{definition}

Clearly, the mixed properties imply the ordinary ones.  
Using results from \cite{CKS87},
Cattani showed that the converse is true as well  \cite{Cat08}.
Since the mixed \HR\ property
is particularly
convenient for applications such as Theorem~\ref{thm:conjs},
we include a self-contained
proof that Lefschetz fans also possess the ``mixed'' properties; see Theorem~\ref{thm:pure->mixed}.

\begin{example}
We remark that any complete simplicial fan $\Sigma$ is Lefschetz.
    In  this case, $\Kone(\Sigma)$ is the cone of K\"ahler classes on the compact complex variety $X_\Sigma$, and 
  there are isomorphisms
  \[
  A^k(\Sigma)\simeq H^{2k}(X_{\Sigma},\R)\simeq \textit{IH}^{2k}(X_{\Sigma},\R).
  \]
The Lefschetz property of $\Sigma$ follows from Poincar\'e duality, the hard Lefschetz theorem, and the Hodge--Riemann relations for the intersection cohomology of $X_\Sigma$ \cite[Section 12.5]{CLS}.
 Alternatively,  one may use Theorem \ref{thm:lefschetz} to deduce the Lefschetz property of $\Sigma$.
\end{example}

\subsection{The weak factorization theorem}

Alexander proved that any 
subdivision of a simplicial complex can be expressed as a sequence of stellar subdivisions of edges and their
inverses \cite{Alexander30}. 
We will use a refined version of his result for simplicial fans.
We continue to write $\TSigma$ for the stellar subdivision of $\Sigma$ by
$\rho$.  For brevity,  we adopt the language of simplicial complexes and call
$\TSigma$ an \emph{edge subdivision} of $\Sigma$ if  the cone $\sigma$ containing $\rho$ in its relative interior is
two-dimensional.

\begin{lemma}\label{lem:codim 2 factors}
There exists a sequence of simplicial fans
  $(\Sigma_0,\Sigma_1,\ldots,\Sigma_n)$ such that
  \begin{enumerate}[(1)]\itemsep 5pt
  \item the initial entry is the fan  $\TSigma$, the final entry is the fan  $\Sigma$, and,
  \item  for each $i$, either $\Sigma_i$ is an edge subdivision of
    $\Sigma_{i+1}$ or $\Sigma_{i+1}$ is an edge subdivision of
    $\Sigma_{i}$.
  \end{enumerate}
  Moreover, if $\Sigma$ is a projective refinement of some fan $\Delta$,
  then so is $\Sigma_i$ for every $i$.
\end{lemma}

\begin{proof}
  We use induction on the dimension of $\sigma$.  The composition of projective maps
  $\TSigma\to\Sigma\to\Delta$ is projective, so the statement when $\dim \sigma=2$ is
  trivial.

  Suppose $\dim \sigma >2$ and $\rho$ is a ray in the relative interior of
  $\sigma$.
  Let $\sigma'$ be any maximal cone of the boundary $\partial\sigma$, and let
  $\mu$ be the unique ray in $\sigma(1)-\sigma'(1)$.  The span of $\set{\rho,
    \mu}$ intersects $\sigma'$ along a ray which we call $\rho'$, and we let
  $\Sigma'=\stellar_{\rho'}\Sigma$.
  Let $\sigma''$ be the $2$-dimensional cone spanned by $\mu$ and $\rho'$.
  The ray $\rho$ lies in $\sigma''$, and we let
  $\Sigma''=\stellar_{\rho}\Sigma'$.
  We claim that $\Sigma''=\stellar_{\rho'}\TSigma$.

  By construction, the fans have the same rays.  To compare the remaining
  cones, it is sufficient  to characterize those
  subsets of rays $\Sigma''(1)$ which fail to span a cone.  We recall
  notation from 
   Section \ref{sec:matroidfan}: If $\Delta$ is a simplicial fan, $I(\Delta)$
  is the Stanley-Reisner ideal of $\Delta$ in a polynomial ring $S(\Delta)$,
  generated by square-free monomials $x_A$, where $A\subseteq\Delta(1)$ runs
  over subsets of rays that do not span a cone.  For cones $\sigma\in \Delta$,
  we will continue to write $x_\sigma$ for the monomial indexed by the
  rays of $\sigma$.
  If we subdivide a cone $\tau$ of a fan $\Delta$ by a ray $\rho$,
  some cones of $\Delta$ are unchanged.  We denote that set by
  \[
  U(\Delta,\tau) \coloneq \Delta - \clst_{\Delta}(\tau).
  \]
By definition of the stellar subdivision, in the polynomial ring $S(\stellar_{\tau}(\Delta))$,  we have 
\[
    I(\stellar_{\rho}(\Delta))=I(\Delta)+(x_\tau)+(x_\rho x_\upsilon\colon
    \upsilon\in U(\Delta,\tau)).
    \]
    Thus, it follows that 
  \begin{align*}
    I(\Sigma'') & = I(\Sigma')+(x_{\sigma''})+
    (x_\rho x_\upsilon\colon \upsilon\in U(\Sigma',\sigma''))\nonumber \\
    &=I(\Sigma)+(x_{\sigma'},x_{\rho'}x_\mu)+
    (x_{\rho'} x_\upsilon\colon \upsilon\in U(\Sigma,\sigma'))+
    (x_\rho x_\upsilon\colon \upsilon \in U(\Sigma',\sigma'')).
  \end{align*}
  On the other hand, we have
  \begin{align*}
    I({\stellar_{\rho'}\TSigma}) &= I(\TSigma)+(x_{\sigma'})+
    (x_{\rho'} x_\upsilon\colon \upsilon\in U(\TSigma,\sigma'))\nonumber \\
    &= I(\Sigma)+(x_{\sigma'},x_{\sigma})+
    (x_{\rho}x_\upsilon\colon \upsilon\in U(\Sigma,\sigma))+
    (x_{\rho'} x_\upsilon\colon \upsilon\in U(\TSigma,\sigma'))\nonumber \\
    &=I(\Sigma)+(x_{\sigma'})+
    (x_{\rho'} x_\upsilon\colon \upsilon\in U(\TSigma,\sigma'))+
    (x_{\rho}x_\upsilon\colon \upsilon\in U(\Sigma,\sigma)),
  \end{align*}
  noting that $x_\sigma=x_\mu x_{\sigma'}$.
  To conclude that $\Sigma''=\stellar_{\rho'}(\TSigma)$, we check that their ideals
  have the same generators,
  using two observations:
  \begin{enumerate}[$\bullet$]\itemsep 5pt
  \item
    There is a
      bijection $\clst_{\Sigma}(\sigma)\simeq\clst_{\Sigma'}(\sigma'')$: if
      $\tau\in\clst_{\Sigma}(\sigma)$, we can write $\tau=\tau'+\sigma'+\mu$
        for some cone $\tau'$: then $\tau'+\rho'+\mu$ is a cone of
        $\clst_{\Sigma'}(\sigma'')$ (because $\sigma''=\rho'+\mu$).  This map
        is easily seen to be invertible.  It follows that
        $U(\Sigma',\sigma'')=U(\Sigma,\sigma)$.
      \item
        There is a bijection $\clst_{\tallSigma}(\sigma')
        \simeq\clst_{\TSigma}(\sigma')$: suppose a cone $\tau\in\Sigma$ contains
        $\sigma'$.  If $\tau\not\supseteq \mu$, then $\tau\in\TSigma$.
        Otherwise, $\tau=\tau'+\mu$ for some $\tau'$, and $\tau'+\rho$ is in
        $\TSigma$.  It follows that
        \[
        U(\TSigma,\sigma')=U(\Sigma,\sigma')\cup\set{\tau\in\TSigma\colon
          \tau\supseteq\mu}.
        \]
  \end{enumerate}
  This gives a commuting diagram of refinements of $\Delta$:
\[
    \begin{tikzpicture}
      \path (0,0) node(a) {$\TSigma$} -- (3,0) node(b) {$\Sigma$};
      \path (1,1.2) node(c) {$\Sigma''$} -- (2,0.6) node(d) {$\Sigma'$};
      \path (a) -- (1.5,-1) node(e) {$\Delta$};
      \draw[->] (a) -- node[above] {$\scriptstyle d$} (b);
      \draw[->] (c) -- node[above left] {$\scriptstyle d-1$} (a);
      \draw[->] (c) -- node[above right] {${\scriptstyle 2}$} (d);
      \draw[->] (d) -- node[above right] {$\scriptstyle d-1$} (b);
      \draw[dotted,->] (a) -- (e);
      \draw[dotted,->] (b) -- (e);
      \draw[dotted,->] (c) -- (e);
      \draw[dotted,->] (d) -- (e);
    \end{tikzpicture}
\]
Since stellar subdivisions are projective,
so are the refinements $\Sigma'\to\Delta$ and $\Sigma''\to\Delta$.
The stellar subdivisions $\TSigma\leftarrow\Sigma''\rightarrow\Sigma'\rightarrow
\Sigma$ take place over cones of dimension $<d$, so by induction there are
sequences of fans from $\TSigma$ to $\Sigma''$ and from $\Sigma'$ to
$\Sigma$ so that each step is an edge subdivision, and each fan is a
projective refinement of $\Delta$.
\end{proof}

\begin{figure}[h]
\centering
\includegraphics[height=6cm]{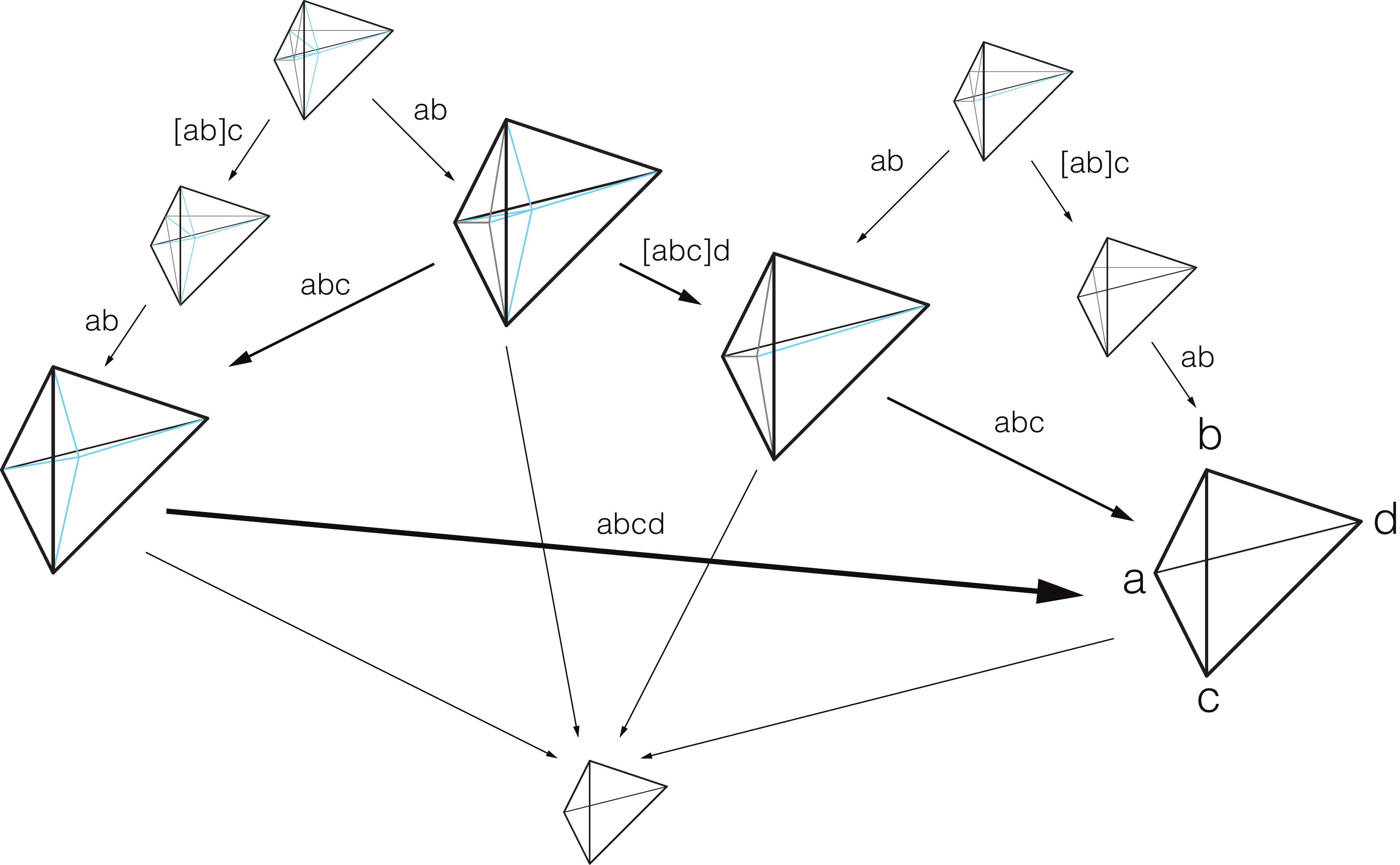}
\caption{Factoring a codimension-$4$ subdivision into edge subdivisions}
\end{figure}

\begin{theorem} \label{thm:fan MFT}
  If $\Sigma$ and $\Sigma'$ are simplicial fans with the same support,  
  there exists a sequence of simplicial fans $\Sigma_0,\Sigma_1,
  \ldots,\Sigma_n$  such that 
    \begin{enumerate}[(1)]\itemsep 5pt
  \item the initial entry is  the fan $\Sigma$, the final entry is the fan  $\Sigma'$, and,
  \item  for each $i$, either $\Sigma_i$ is an edge subdivision of
    $\Sigma_{i+1}$ or $\Sigma_{i+1}$ is an edge subdivision of
    $\Sigma_{i}$.
  \end{enumerate}
  Furthermore, the entries can be chosen in such a way that
  there is an index $i_0$ for which $\Sigma_i$ is a projective refinement
  of $\Sigma$ for all $i\leq i_0$, and $\Sigma_i$ is a projective refinement
  of $\Sigma'$ for all $i\geq i_0$.
\end{theorem}
\begin{proof}
  By \cite[Theorem\ 11.1.9]{CLS}, there exists a sequence of stellar
  subdivisions of both $\Sigma$ and $\Sigma'$ that refine those fans,
  respectively, to unimodular fans. Thus, by Proposition \ref{lem:blowup K} and Lemma \ref{lem:codim 2 factors}, 
  we can reduce to the case where $\Sigma$ and $\Sigma'$ are
  both unimodular fans.
  
  By \cite[Theorem A]{Wlo97}, there is a sequence of simplicial fans
  for which $\Sigma_i$ and $\Sigma_{i-1}$ differ by a stellar subdivision,
  for each $i$.  The assertion that these can be chosen to
  have the second property follows from \cite[Theorem\ 2.7.1]{AKMW02}.
  As stated,  \cite[Theorem\ 2.7.1]{AKMW02} applies to the more general setting of
  toroidal embeddings and polyhedral complexes, and the proof given in \cite{AKMW02}
   specialize to the case of toric varieties and fans.\footnote{The birational cobordism used in the proof can
  be chosen to be toric, using a toric resolution of singularities.}

  Suppose now that $i\leq i_0$ and $\Sigma_i=\stellar_\rho(\Sigma_{i-1})$.
  Since $\Sigma_{i-1}$ is a projective refinement of $\Sigma$, we use
  Lemma~\ref{lem:codim 2 factors} to obtain a sequence of fans between
  $\Sigma_i$ and $\Sigma_{i-1}$ in which consecutive fans differ
  by edge subdivisions,
  and which are also projective refinements of $\Sigma$.  The remaining
  cases are treated by exchanging $i$ with $i-1$, and $\Sigma$ with $\Sigma'$.
\end{proof}

In terms of toric varieties, edge subdivisions correspond to morphisms that are \emph{semismall} in the sense of Goresky--MacPherson.
In projective geometry, the semismallness is particularly convenient for transferring the Lefschetz property,
as  pullbacks of ample line bundles by semismall maps  satisfy the hard Lefschetz property and the Hodge--Riemann bilinear relations \cite{dCM}. 
We will use Theorem \ref{thm:fan MFT} to prove Theorem \ref{thm:lefschetz} by establishing  an analogous property in the context of Lefschetz fans.

\subsection{Chow rings of stellar subdivisions}
We continue to write $\TSigma$ for the stellar subdivision of $\Sigma$ by
$\rho$, a ray in the relative interior of $\sigma$.
The primitive ray generator of $\rho$ is given by
\[
  \e_\rho \coloneq \sum_{\nu\in\sigma(1)}a_\nu\mathbf{e}_\nu,
  \]
for some positive rational coefficients $\set{a_\nu}$.
We relate the Chow rings of $\Sigma$ and $\TSigma$.
For economy of notation, we
abbreviate $E\coloneq \clst_{\TSigma}(\rho)$ and $Z\coloneq \clst_{\Sigma}(\sigma)$.  We
let $j\colon E\to\TSigma$ denote the inclusion of fans, and
let $q\colon E\to Z$ denote the restriction
of $p\colon\TSigma\to\Sigma$.  We will write $i$ in place of $i_\sigma$
through the rest of this section:
\[
\begin{tikzcd}
  E\ar[d,"q"]\ar[hookrightarrow,r,"j"] &  
  \TSigma\ar[d,"\bl"] \\
  Z\ar[hookrightarrow,r,"i"] & \Sigma
\end{tikzcd}
\]

A straightforward calculation shows
that the pullback homomorphism $\bl^*\colon A(\Sigma)\to A(\TSigma)$ 
is determined by the formula 
\[
p^*(x_\nu) =  \begin{cases}
x_\nu & \text{if $\nu\not\in\sigma(1)$,}\\
x_\nu+a_\nu x_\rho & \text{if $\nu\in \sigma(1)$.}
\end{cases}
\]
Since $\bl$ is a proper map of fans \cite[Theorem 3.4.11]{CLS},
there is a pushforward 
$\bl_*\colon A(\TSigma)\to A(\Sigma)$, which is a homomorphism
of $A(\Sigma)$-modules.  
By definition \cite[Section 1.4]{Fulton}, for any $\widetilde\tau \in \TSigma$, we have 
\[
\bl_*(x_{\widetilde\tau})=\begin{cases} x_\tau & \text{if $\bl(\widetilde\tau)
\subseteq \tau$ and $\dim\widetilde\tau=\dim\tau$,}\\
0 & \text{if otherwise.}
\end{cases}
\]

\begin{proposition}\label{thm:blowup chow}
If $\sigma$ is two-dimensional, then 
the map
\[
\bl_*\oplus q_*j^*\colon
A(\TSigma) \rightarrow A(\Sigma)\oplus A(\st_{\tallSigma}(\sigma))[-1]
\]
is an isomorphism of graded $A(\Sigma)$-modules.
\end{proposition}

The proof, given below, uses some preliminary calculations.  Let $\nu_1$,
$\nu_2$ be the rays of $\sigma$ and $\e_1$, $\e_2$ the primitive ray
generators.  Then $\e_\rho=a_1\e_1+a_2\e_2$ for some positive rational coefficients
$a_1,a_2$.  For $i,j\in\set{1,2,\rho}$, we let $m_{i,j}$ denote the
index of $\Z\set{\e_i,\e_j}$ inside $\R\set{\e_i,\e_j}\cap \N_{\Z}$.
Computing determinants, we see that
\[
  a_1=m_{2,\rho}/m_{1,2} \ \ \text{and} \ \ 
    a_2=m_{1,\rho}/m_{1,2}.
  \]


\begin{lemma}\label{lem:qj_formula}
  We have
  \[
  q_*j^*j_*q^*=-\frac{m_{1,2}}{m_{1,\rho}m_{2,\rho}}.
  \]
\end{lemma}
\begin{proof}
Consider an element $v\in A(Z)$.  Since $j^*$ is surjective,
we have $q^*(v)=j^*(u)$ for some $u\in A(\TSigma)$.  Then, by the projection formula,
\begin{align*}
  q_*j^*j_*q^*(v) & = q_*j^*j_*q^*(1_Z\cdot v) \\
  &= q_*j^*j_*\big(q^*(1_Z)\cdot q^*(v)\big)\\ 
  &= q_*j^*j_*\big(q^*(1_Z) \cdot j^*(u)\big)\\
  &= q_*j^*\big(j_*q^*(1_Z)\cdot u\big)\\
  &= q_*\big(j^*j_*q^*(1_Z)\cdot j^*(u)\big)\\
  &= q_*\big(j^*j_*q^*(v)\big)\\
  &= q_*j^*j_*q^*(1_Z)\cdot v,
\end{align*}
so it is enough to verify the claim on the fundamental class $1_Z$.

Let $h=j^*(x_\rho)\in A^1(E)$.  By definition, we have
\[
  (q_*j^*j_*q^*)(1_Z) = q_*j^*j_*(1_E)
  = q_*j^*(x_\rho)
  = q_*(h).
  \]
We extend $\set{\e_1,\e_2}$
to a basis for $\N$, and write $\set{\e_1^*,\e_2^*,\ldots}$ for the
dual basis.  The piecewise-linear function $x_\rho-a_1^{-1}\e_1^*$ 
is equivalent to $x_\rho$,
and its values on
$\e_1$, $\e_2$, and $\e_\rho$ are $-a_1^{-1}$, $0$, and $0$, respectively.
This is to say that $h=j^*(-a_1^{-1}x_1+g)$ for some piecewise linear function $g$
on $\TSigma$ which is zero on $\e_1$, $\e_2$, and $\e_\rho$, and hence
\[
h=-a_1^{-1}m_{1,\rho}^{-1}x_1+j^*(g).
\]
As $j^*(g)$ is a linear combination of Courant functions $x_\nu$ for 
rays $\nu$ in $\st_{\TSigma}(\rho)$ not contained in the support of $\sigma$, we have $q_*j^*(g)=0$.
On the other hand, 
$q_*(x_1)=1_Z$, so 
\[
q_*(h)=-a_1^{-1}m_{1,\rho}^{-1} 1_Z=-m_{1,2}/(m_{1,\rho}m_{2,\rho})1_Z. \qedhere
\]
\end{proof}
\begin{proof}[Proof of Proposition~\ref{thm:blowup chow}]
  Let $\psi(u)=(p_*(u),q_*j^*(u))$ for $u\in A(\TSigma)$, and let
  $\phi(u,v)=\bl^*(u)+j_*q^*(v)$ for $(u,v)\in A(\Sigma)\oplus A(Z)[-1]$.  We
  first check that $\psi\circ \phi$ is an isomorphism.

 Observe that  $p_*p^*=1$, because $p$ is birational, and $q_*q^*=0$, because $q$ has positive relative dimension. Therefore, we have
  \begin{align*}
    \psi\circ\phi(u,v) &= \bl_*\big(\bl^*(u)+j_*q^*(v)\big) +
    q_*j^*\big(\bl^*(u)+j_*q^*(v)\big)\\
    &= \bl_*\bl^*(u)+\bl_*j_*q^*(v) + q_*j^*\bl^*(u)+q_*j^*j_*q^*(v)\\
    &= \bl_*\bl^*(u)+ i_*q_*q^*(v) + q_*q^*i^*(u) + q_*j^*j_*q^*(v)\\
    &= u + q_*j^*j_*q^*(v),
  \end{align*}
  which is invertible by Lemma~\ref{lem:qj_formula}.  It follows that
  $\phi$ is injective, and we now argue that it is also surjective.

  Since squarefree monomials span $A(\TSigma)$, it is enough to show that 
  each monomial $x_\tau$ is of the form $\bl^*(u)+j_*q^*(v)$ for suitable $u\in A(\Sigma)$
  and $v\in A(Z)$.  If none of $\nu_1$, $\nu_2$ or $\rho$ is contained in
  $\tau$, clearly $x_\tau=\bl^*(x_\tau)$.  Noting that no cone of
  $\TSigma$ contains both $\nu_1$ and $\nu_2$, it remains to consider
  the following three cases.
  \begin{itemize}
  \item[Case 1:] Suppose $\set{\nu_1,\nu_2,\rho}\cap\tau(1)=\set{\rho}$.
   If we set $\tau'=\tau-\set{\rho}$, then
    \[
    j_*q^*(x_{\tau'})=j_*(x_{\tau'})=x_\rho x_{\tau'},
    \]
    so we may take $u=0$ and $v=x_{\tau'}$.
  \item[Case 2:] Suppose $\set{\nu_1,\nu_2,\rho}\cap\tau(1)=\set{\nu_1}$.
    If we set $\tau'=\tau-\set{\nu_1}$, then
    \begin{align*}
      x_1x_{\tau'} &= (x_1+a_1x_\rho)x_{\tau'} - a_1x_\rho x_{\tau'}.
    \end{align*}
    The first summand equals $\bl^*(x_1x_{\tau'})$, and the second summand
    is in the image of $\phi$ in view of the previous case.
  \item[Case 3:]
    Suppose $\set{\nu_1,\nu_2,\rho}\cap\tau(1)=\set{\nu_1,\rho}$.
    We set $\tau'=\tau-\set{\nu_1,\rho}$, and
    extend the vectors $\e_1,\e_2$ to a basis for $N$.  Then
    the linear function $\e_2^*$ may be written
    \[
    \e_2^*=x_2+a_2x_\rho+g,
    \]
    where $g$ is a piecewise linear function vanishing on the rays $\set{\nu_1,\nu_2,\rho}$.
    Since the class of $\e_2^*$ is zero in the Chow ring of $\TSigma$, multiplying it by $x_1x_{\tau'}$ gives
    \[
    0 = 0\cdot x_1x_{\tau'}=0+a_2x_1x_\rho x_{\tau'}+gx_1x_{\tau'},
    \]
    because $x_1x_2=0$.  Thus $a_2x_\tau=-gx_1x_{\tau'}$, which
    is in the image of $\phi$ by the previous case.
    \end{itemize}
    This completes the proof of Proposition~\ref{thm:blowup chow}.
\end{proof}

\begin{corollary}\label{cor:pi* inj} 
  The pullback homomorphism  $\bl^*$ is injective, and it restricts to an
  isomorphism in degree $d=\dim \Sigma$.
\end{corollary}
\begin{proof}
  The isomorphism $\phi$ restricts to $\bl^*$ on $A(\Sigma)$, so $\bl^*$ is
  injective.  
  Since $\st_{\Sigma}(\sigma)$ is $(d-2)$-dimensional,   $A^{d-1}(\st_{\Sigma}(\sigma))=0$, and hence the isomorphism $\phi$ agrees with
  $\bl^*$ in degree $d$.
\end{proof}

\subsection{Hodge--Riemann forms and their signatures}\label{ss:signature}

Our goal in the next few pages is to understand how the Lefschetz property
behaves under edge subdivisions. 
In  this subsection, we fix a $d$-dimensional simplicial fan $\Sigma$ that satisfies Poincar\'e duality (Definition \ref{def:lefschetz})
and 
$k \le \frac{d}{2}$.
Suppose that the multiplication by  $L\in A^{d-2k}(\Sigma)$ is
an isomorphism in degree $k$.  
Using  Poincar\'e duality, one can check directly that  the multiplication by  $L\in A^{d-2k}(\Sigma)$ is
an isomorphism in degree $k$ if and only if the corresponding \emph{Hodge--Riemann 
form} 
\[
\hr^k(\Sigma,L)\colon A^k(\Sigma) \times A^k(\Sigma) \longrightarrow \mathbb{R}, \qquad (\eta_1,\eta_2) \longmapsto (-1)^k \deg(L \eta_1 \eta_2)
\]
is nondegenerate.
 Thus, in this case, $\hr^k(\Sigma,L)$ has
$b_k^+$ positive eigenvalues and $b_k^-$ negative eigenvalues, where
$b_k^++b_k^-$ is the dimension of $A^k(\Sigma)$.
We use its signature $b_k^+-b_k^-$ can be used to characterize the \HR\ property.
This characterization appears as \cite[Proposition  7.6]{AHK15} and 
\cite[Theorem  8.6]{McMullen} in the case when $L=\ell^{d-2k}$ for $\ell \in A^1(\Sigma)$.

In what follows, we write  $b_k(\Sigma)$ for the dimension of $A^k(\Sigma)$.
Given $L\in A^{d-2k}(\Sigma)$ and $\ell_0 \in A^1(\Sigma)$, we define the \emph{primitive part} of $A^k(\Sigma)$ to be the subspace
\[
\PA^k(\Sigma,\ell_0,L)\coloneq \Big\{ \eta \in A^k(\Sigma) \, \, | \, \,  \ell_0  \cdot L  \cdot \eta=0\Big\}.
\]
For simplicity, when $L=\ell_0^{d-2k}$ and there is no possibility of confusion, we write $\hr^k(\Sigma,\ell_0)$ for $\hr^k(\Sigma,L)$
and $\PA^k(\Sigma,\ell_0)$ for $\PA^k(\Sigma,\ell_0,L)$.

\begin{proposition}\label{prop:fixed sig}
  If $U\subseteq A^{d-2k}(\Sigma)$
  is a connected subset in the Euclidean topology and
  if the Hodge--Riemann form  $\hr^k(\Sigma,\lef)$ is nondegenerate 
  for all
  $\lef\in U$, then the signature of $\hr^k(\Sigma,\lef)$ is constant for
  all $\lef\in U$.
\end{proposition}
\begin{proof}
  The eigenvalues of $\hr^k(\Sigma,\lef)$ are real, and they vary continuously
  with $\lef$.  By hypothesis, they are all nonzero for $\lef\in U$, so their
  signs are constant on $U$, because $U$ is connected.
\end{proof}

We write $\Sym^{d-2k} \Kone(\Sigma) \subseteq A^{d-2k}(\Sigma)$ for the subset of products of elements of $\Kone(\Sigma) \subseteq A^1(\Sigma)$.

\begin{proposition}[$\HR$ signature test]\label{thm:HR as signature}
  Suppose that $\Sigma$ satisfies the conditions 
  \begin{enumerate}[(1)]\itemsep 5pt
  \item   $\hr^i(\Sigma,\lef)$ is nondegenerate for all $0\leq i\leq k$
    and all $\lef\in \Sym^{d-2i} \Kone(\Sigma)$, and
    \item $\hr^i(\Sigma,\lef)$ is
  positive definite on the kernel of the multiplication by $\ell_0 L$ 
  for all $\ell_0\in \Kone(\Sigma)$, all $\lef\in \Sym^{d-2i} \Kone(\Sigma)$, and all $i<k$.
  \end{enumerate}
  Then 
  $\hr^k(\Sigma,L)$ is positive definite on the kernel of the multiplication by $\ell_0 L$ 
  for all $\ell_0\in \Kone(\Sigma)$ and all  $L\in\Sym^{d-2k}\Kone(\Sigma)$ if and only if its signature 
  equals
  \[
  \sum_{i=0}^k (-1)^{k-i}\big(
  b_i(\Sigma)- b_{i-1}(\Sigma)\big).
  \]
\end{proposition}

\begin{proof}
The proof is the same with the one given in 
\cite[Proposition 7.6]{AHK15}   for the special case 
  $L=\ell_0^{d-2k}$.
The result follows from the induction on $k$ and the Lefschetz decomposition
\[
A^k(\Sigma)=\PA^k(\Sigma,\ell_0,L) \oplus \ell_0 A^{k-1}(\Sigma),
\]   
which is orthogonal for the Hodge--Riemann form  $\hr^k(\Sigma,L)$.
\end{proof}



\begin{corollary}\label{cor:HR+1}
  If $\Sigma$ satisfies mixed $\HR^i$ for all $i<k$, mixed $\HL^k$,
  as well as  $\HR^k(\lef')$ for some $\lef'\in  \Sym^{d-2k}\Kone(\Sigma)$,
  then $\Sigma$ satisfies $\HR^k$.
\end{corollary}
\begin{proof}
  Let $L\in \Sym^{d-2k}\Kone(\Sigma)$ be any element.  By the
  hypothesis mixed $\HL^k$, the Hodge--Riemann form $\hr^k(\Sigma,\lef)$ is nondegenerate.  By
  Proposition~\ref{prop:fixed sig}, it has the same signature as
  $\hr^k(\Sigma,\lef')$.  Since  $\Sigma$ satisfies mixed
  $\HR^i$ for $i<k$, Proposition \ref{thm:HR as signature} shows that 
  $\HR^k(\lef)$ and $\HR^k(\lef')$ are equivalent.
\end{proof}

 Let $\TDelta=\st_{\TSigma}(\rho)$ and $\Delta=\st_{\tallSigma}(\sigma)$.
 In the case of a star-shaped blowup, the signature test simplifies
slightly. 
In this case, by Propositions  \ref{thm:blowup chow} and \ref{thm:HR as signature},
 the signature of $\hr^k(\TDelta,L)$  satisfying the Hodge--Riemann relations is   $b_k(\Delta)-b_{k-1}(\Delta)$ for $k<\frac{d}{2}$.


\subsection{The Lefschetz property and  edge subdivisions}
\label{ss:bl lefschetz}

With these preparations, we now set out to show that
the Lefschetz property of a fan is
unaffected by edge subdivisions and their inverses.  The precise statements
and their proofs appear in Section \ref{ss:blowup theorems}
as Theorems \ref{thm:lefschetz I} and \ref{thm:lefschetz II}.
Here, we first consider Poincar\'e duality, and we first do so for
star-shaped subdivisions.

Let $\TSigma$ be the stellar subdivision of a $d$-dimensional simplicial fan $\Sigma$ by a ray $\rho$ in a two-dimensional cone $\sigma$.
As before, we set 
$\TDelta=\st_{\TSigma}(\rho)$ and $\Delta=\st_{\tallSigma}(\sigma)$.

\begin{proposition}\label{prop:star PD}
  Poincar\'e duality holds for $\TDelta$ if and only
  if it holds for $\Delta$.
\end{proposition}

\begin{proof}
  Assume that \PD\ holds for at least one of $\TDelta$ and $\Delta$.
  By Proposition \ref{thm:blowup chow}, for all positive $i$,
\[
A^i(\TDelta)\simeq A^i(\Delta)\oplus x_\rho A^{i-1}(\Delta).
\]
  We see that $A^{d-2}(\Delta)\simeq A^{d-1}(\TDelta)$, so if one of $\Delta$
  or $\TDelta$ has a fundamental weight, they both do.
  By inspection, $b_i(\Delta)=b_{d-2-i}(\Delta)$ for all $i$ if and only if
  $b_i(\TDelta)=b_{d-1-i}(\TDelta)$ for all $i$.  So we may assume both sets
  of equalities hold.

  For any $u\in A^i(\TDelta)$ and $v\in A^{d-1-i}(\TDelta)$, we write
  $u=u_0+u_1x_\rho$ and $v=v_0+v_1x_\rho$ where $u_0$, $u_1$, $v_0$, $v_1$
  are elements of $A(\Delta)$ of degrees $i$, $i-1$, $d-1-i$, and $d-2-i$,
  respectively.  Then $u_0v_0\in A^{d}(\Delta)=0$,
  and $x_\rho^2=c_1\cdot x_\rho+c_2$ for some $c_1,c_2\in A(\Delta)$.
  With respect to the decomposition above, the matrix of the
  multiplication pairing has the form
\[
M^i(\TDelta)=
\Biggl(\mkern-5mu
\begin{tikzpicture}[baseline=-0.65ex]
  \matrix[
    matrix of math nodes,
    column sep=1ex,
  ] (m)
         {
           0 & -M^{i-1}(\Delta) \\
           -M^i(\Delta) & *\\
         };
   \coordinate (top) at ([yshift=14pt,xshift=0.5ex]m-2-1.north east);
   \path (top) +(0,-1.1) coordinate (bottom);  
   \coordinate (left) at (m-1-1.south west);
   \path (left) +(3.3,0) coordinate (right);
   \draw[dotted] (top) -- (bottom);
   \draw[dotted] (left) -- (right);
\end{tikzpicture}
\mkern-5mu\Biggr),
\]
where $M^i(\Delta)$ denotes the matrix of the pairing $A^i(\Delta)\times
A^{d-2-i}(\Delta)\to\R$.
  Thus if each matrix $M^i(\TDelta)$ is invertible, so is each matrix
  $M^i(\Delta)$,
  and conversely. Therefore, if either $\Delta$ or $\TDelta$ has \PD, then they
  both do.
\end{proof}

\begin{proposition}\label{prop:PD}
  Suppose that
  Poincar\'e duality holds for $\Delta$. 
  Then Poincar\'e duality holds for $\TSigma$ if and only if it holds for $\Sigma$.
\end{proposition}
\begin{proof}
  Let us assume that at least one of $\Sigma$ and $\TSigma$ has Poincar\'e
  duality, then show that they both do.  For dimensional reasons, 
  $\TSigma$ must be an ordinary subdivision.  
  By Corollary~\ref{cor:pi* inj}, we have
  $A^d(\TSigma)\simeq A^d(\Sigma)$, and they have the common degree map.

  By Proposition \ref{thm:blowup chow} 
  and Poincar\'e duality for $\st_{\tallSigma}(\sigma)$,
  we have $b_i(\Sigma)=b_{d-i}(\Sigma)$
  and  $b_i(\TSigma)=b_{d-i}(\TSigma)$ for all $0\leq i\leq d$.
  Since $A^{s}(\Sigma)\times A^{t}(\st_{\tallSigma}(\sigma))\to
  A^{s+t}(\st_{\tallSigma}(\sigma))$ is the zero map when $s+t>d-2$, ordering
  bases compatibly with the decomposition in Proposition \ref{thm:blowup chow} gives a block-diagonal matrix:
\[
M^i(\TSigma)=
\Biggl(\mkern-5mu
\begin{tikzpicture}[baseline=-0.65ex]
  \matrix[
    matrix of math nodes,
    column sep=1ex,
  ] (m)
         {
           M^i(\Sigma) & 0\\
           0 & M^{i-1}(\Delta)\\
         };
   \coordinate (top) at ([xshift=0.5ex]m-1-1.north east);
   \path (top) +(0,-1.1) coordinate (bottom);  
   \coordinate (left) at (m-1-1.south west);
   \path (left) +(3.3,0) coordinate (right);
   \draw[dotted] (top) -- (bottom);
   \draw[dotted] (left) -- (right);
\end{tikzpicture}
\mkern-5mu\Biggr)
\]
Clearly $M^i(\TSigma)$ has full rank if and only if $M^i(\Sigma)$
and $M^{i-1}(\Delta)$ both do as well, which completes the proof.
\end{proof}

\begin{lemma}\label{lem:raysum}
  Suppose that Poincar\'e duality holds for $\Sigma$. If  $I\subseteq
  \Sigma(1)$ is a subset of rays for which
  $\{x_\nu\}_{\nu\in I}$ spans $A^1(\Sigma)$,  then 
  $
  \oplus \, i^*_\nu\colon A^i(\Sigma)\to\bigoplus_{\nu\in
    S}A^i(\st_{\tallSigma}(\nu))
  $
  is injective for all $0\leq i<d$.
\end{lemma}

\begin{proof}
Suppose $i^*_\nu(u)=0$ for each ray $\nu$.  Then $i_*^\nu i^*_\nu(u)=
x_\nu u=0$ for a set of generators $x_\nu$ of $A(\Sigma)$.  Since
$A(\Sigma)$ has no nonzero socle in degree $<d$ by Poincar\'e duality, the  element $u$ must be zero. 
\end{proof}

\begin{proposition}\label{prop:HR=>HL+1}
  Suppose that Poincar\'e duality holds for $\Sigma$. 
  If
  $\st_{\tallSigma}(\nu)$ satisfies mixed \HR\ for each ray $\nu\in\Sigma(1)$,
  then $\Sigma$ satisfies mixed \HL.
\end{proposition}

\begin{proof}
  Let $\lef\coloneq \ell_1\cdots\ell_{d-2k}$ be an element of $\Sym^{d-2k} \Kone(\Sigma)$, and
  consider the map $\lef \, \cdot\colon A^k(\Sigma)\to A^{d-k}(\Sigma)$.
  By Poincar\'e duality, we know that the domain and the target have the same dimension, 
  so it is enough to show that $\lef \, \cdot$ is injective.  Suppose, then, that
  $\lef \cdot u=0$ for some $u\in A^k(\Sigma)$.

  Let $\lef' \coloneq  \ell_2\cdots \ell_{d-2k}$.
 Note that,  for each index $i$ and each ray $\nu$ in $\Sigma$, the pullback $i^*_\nu(\ell_i)$  belongs to $ \Kone(\st_{\tallSigma}(\nu))$. 
  Furthermore, since $\lef \, \cdot u=0$, the pullback of
  $u$ around $\nu$ is primitive:
  \[
  i_\nu^*(u)\in \PA^k(\st_{\tallSigma}(\nu),i_\nu^*(\ell_1),i^*_\nu(\lef')).
  \]
  We may write $\ell_1=\sum_{\nu\in \Sigma(1)}c_\nu x_\nu$ where
  each coefficient $c_\nu>0$, since we can represent $\ell_1$ by a piecewise linear function
  which is strictly positive on each ray.  We have 
  \begin{align*}
    0 =\deg_\Sigma(\lef \, \cdot u\cdot u)
    =\deg_\Sigma(\sum_{\nu\in \Sigma(1)}c_\nu x_\nu \lef'\cdot u\cdot u)
    &=\sum_{\nu}c_\nu \deg_{\st_{\tallSigma}(\nu)}(i^*_\nu(\lef')\cdot
    i^*_\nu(u)\cdot i^*_\nu(u))\\
    &=(-1)^{k-1}\sum_{\nu\in \Sigma(1)}c_\nu
    \angl{i_\nu^*(u),i_\nu^*(u)}_{i^*_\nu(\lef')},
  \end{align*}
  where the last summands are the Hodge--Riemann forms for $i^*_\nu(\lef')$.
  Since the $c_\nu$'s are strictly positive, each summand is zero, and the
  mixed \HR\ property in $\st_{\tallSigma}(\nu)$ implies $i^*_\nu(u)=0$,
  for each $\nu$.  By Lemma~\ref{lem:raysum}, we have $u=0$, and $\lef \, \cdot$ is
  injective.
\end{proof}

As an application, we see that the mixed Lefschetz properties  in Definition \ref{def:mixed lefschetz} are
actually no stronger than the pure ones. See \cite{Cat08} for a discussion in a more general context.

\begin{theorem}\label{thm:pure->mixed}
  If $\Sigma$ is a Lefschetz fan, then it also has the mixed
  \HL\ and mixed \HR\ properties.
\end{theorem}

\begin{proof}
  We use induction
  on dimension.  If $\dim\Sigma=1$, the mixed and
  pure properties are identical, so let us suppose the claim is true for all
  Lefschetz fans of dimension less than $d$, for some $d>1$.
  Let $\Sigma$ be a Lefschetz fan of dimension $d$.  By induction,
  $\st_{\tallSigma}(\nu)$ satisfies mixed \HR\ for all rays
  $\nu\in\Sigma(1)$.  
  By Proposition~\ref{prop:HR=>HL+1}, then $\Sigma$ satisfies mixed \HL.

  Now we establish mixed \HR\ for $\Sigma$.  For any $\ell\in \Kone(\Sigma)$
  and $0\leq k\leq \frac{d}{2}$, the ``pure'' property $\HR^k(L')$ holds
  for $L'=\ell^{d-2k}$.  Corollary~\ref{cor:HR+1} states that
  mixed \HL\ and mixed $\HR^i$ for $i<k$ implies mixed $\HR^k$.  Setting
  $k=0$, we see $\Sigma$ has the mixed $\HR^0$ property.  Arguing by
  induction on $k$, we obtain mixed $\HR^k$ for all $k\leq \frac{d}{2}$.
\end{proof}

We now examine how the Hodge--Riemann forms fare under stellar subdivisions.
As before, we write $\bl\colon \TSigma\to\Sigma$ for the map of fans given by the edge subdivision under consideration, write
$x_\rho=a_1x_1+a_2x_2$ for some positive scalars $a_1$, $a_2$, and consider the diagram
\[
\begin{tikzcd}
  E\ar[d,"q"]\ar[hookrightarrow,r,"j"] &  
  \TSigma\ar[d,"\bl"] \\
  Z\ar[hookrightarrow,r,"i"] & \Sigma.
\end{tikzcd}
\]

\begin{lemma}\label{lem:push pi}
  We have $\bl_*(x_\rho)=0$ and $\bl_*(x_\rho^2)=-a_1a_2x_\sigma$.
\end{lemma}

\begin{proof}
  The first identity follows from the definition pushforward $\bl_*$ between the Chow groups.
  Now let $x_1$, $x_2$ be the Courant functions for the rays $\nu_1$, $\nu_2$
  of the cone $\sigma$, so $x_\sigma=x_1x_2$. For $i=1,2$, we have
\[
0 = \bl_*(x_\rho)x_i = \bl_*\big(x_\rho(x_i+a_ix_\rho)\big),
\]
  so $\bl_*(x_\rho x_i)=-a_i \bl_*(x_\rho^2)$.
  Since $\set{\nu_1,\nu_2}$ is not contained in a
  cone of $\TSigma$, we have 
\[
    x_\sigma = \bl_*\bl^*(x_1x_2)
    =\bl_*\big((x_1+a_1x_\rho)(x_2+a_2x_\rho)\big)
    =(0-2a_1a_2+a_1a_2)\bl_*(x_\rho^2). \qedhere
    \]
\end{proof}

\begin{lemma}\label{lem:codim2 ortho}
  Suppose that Poincar\'e duality holds for $\Sigma$. If $\TSigma$ is an ordinary edge subdivision of $\Sigma$,
 then, for
  all $0\leq k\leq \frac{d}{2}$ and all $\lef\in \Sym^{d-2k}\Kone(\Sigma)$, we have
the orthogonal direct sum
\[
  \hr^k(\TSigma,\bl^*\lef)\cong \hr^k(\Sigma,\lef)\oplus
  \hr^{k-1}(\st_{\tallSigma}(\sigma),i_\sigma^*(\lef)). 
\]
\end{lemma}

\begin{proof}
  We consider $\hr^k(\TSigma,\bl^*\lef)$ under the isomorphism
  $\phi\colon A^k(\Sigma)\oplus A^{k-1}(\st_\Sigma(\sigma)) \cong A^k(\TSigma)$ from
  Proposition \ref{thm:blowup chow}.  Recall that $\phi(u,v)=\bl^*(u)+j_*q^*(v)$.
  Let $\hat{v}$ be any preimage of $v$ through
  the surjective map $i_\sigma^*$: then $j_*q^*(v)=j_*j^*\bl^*(\hat{v})=
  x_\rho\bl^*(\hat{v})$.


  We use the notation $\angl{-,-}$ to pair elements under the various Hodge--Riemann forms, and check first that the two summands are indeed orthogonal.
  We calculate:
%
  \begin{align*}
    (-1)^k\angl{(u,0),(0,v)} &= \deg_{\TSigma}\big(\bl^*(\lef)\cdot
    \bl^*(u)\cdot x_\rho\bl_\sigma^*(\hat{v})\big)\\
    &=\deg_{\Sigma}\big(\lef \, \cdot u\hat{v}\cdot \bl_*(x_\rho)\big)=0,
  \end{align*}
  using the projection formula and the fact that $\bl_*(x_\rho)=0$.
  
  If $u,v\in A^k(\Sigma)$, the equality
  $\angl{(u,0),(v,0)}_{\bl^*(\lef)}=\angl{u,v}_{\lef}$ is
  straightforward.
  If $u,v\in A^{k-1}(\st_{\tallSigma}(\sigma))$, as before write
  $u=i_\sigma^*(\hat{u})$ and $v=i_\sigma^*(\hat{v})$ for some
  $\hat{u},\hat{v}\in  A^{k-1}(\Sigma)$. Then, calculating as above,
  \begin{align*}
    \angl{(0,u),(0,v)} &= (-1)^k\deg_{\TSigma}\big(\bl^*(\lef)\cdot
    \bl^*(\hat{u})\bl^*(\hat{v})\cdot x^2_\rho\big)\\
    &=(-1)^k\deg_{\Sigma}\big(\lef \, \cdot \hat{u}\hat{v}
    \cdot \bl_*(x^2_\rho)\big)\\
    &=-(-1)^ka_1a_2\deg_{\Sigma}\big(\lef \, \cdot \hat{u}\hat{v}
    \cdot x_\sigma\big)\quad\text{by Lemma~\ref{lem:push pi};}\\
    &=(-1)^{k-1}a_1a_2\deg_{\st_{\tallSigma}(\sigma)}\big(
    i^*_\sigma(\lef)\cdot i^*_\sigma(\hat{u})i^*_\sigma(\hat{v})
    \big)=a_1a_2\angl{u,v}_{i^*_\sigma(\lef)}.
  \end{align*}
  The conclusion follows, since $a_1,a_2>0$.
\end{proof}

Next we address star-shaped subdivisions.
Set $e \coloneq \dim \st_{\tallSigma}(\sigma)=\dim\Sigma-2$.  

\begin{lemma}\label{lem:schur}
  Suppose $P$ and $Q$ are $n\times n$ matrices with real entries
  and $Q=Q^T$.  Let
  \[
  M\coloneq \begin{pmatrix}
    0 & P\\ P^T & Q
  \end{pmatrix}.
  \]
  If $P$ is nonsingular, then $M$ has signature zero.
\end{lemma}

\begin{proof}
  Assume first that $Q$ is invertible, and let $S=-PQ^{-1}P^T$ (the Schur
  complement.)  Then it is easily seen that $M$ is congruent to a
  block-diagonal matrix:
  \[
  M = \begin{pmatrix} I_n & PQ^{-1} \\ 0 & I_n \end{pmatrix}
  \begin{pmatrix} S & 0 \\ 0 & Q \end{pmatrix}
  \begin{pmatrix} I_n & 0\\ Q^{-1}P^T & I_n \end{pmatrix},
  \]
  and the signature of $S$ is the negative of the signature of $Q$.  It
  follows that $M$ has signature zero.

  Now suppose $Q$ is singular.
  We replace $Q$ by $Q(\epsilon)$ to define
  $M(\epsilon)$ as above, for some real, invertible symmetric matrices
  $Q(\epsilon)$ with $\lim_{\epsilon\to0}Q(\epsilon)=Q$.
  Then $\det(M(\epsilon))=  (-1)^n\det(P)^2\neq0$, regardless of $\epsilon$,
  so the argument above
  shows $M(\epsilon)$ has $n$ positive eigenvalues and $n$ negative eigenvalues.
  By continuity, so does $M$.
\end{proof}

The last result in this section relates \HL\ and \HR\ along an edge
subdivision.

\begin{proposition}\label{prop:star-shape HR}
Suppose that at least one of $\st_{\tallSigma}(\sigma)$ and $\st_{\TSigma}(\rho)$
satisfies Poincar\'e duality, and that $\ell\in \Kone(\st_{\tallSigma}(\sigma))$ has
the hard Lefschetz property.  Then
\begin{enumerate}[(1)]\itemsep 5pt
\item  $\ell_\epsilon\coloneq \ell-\epsilon\cdot x_\rho\in
\Kone(\st_{\TSigma}(\rho))$ has the \HL\ property for sufficiently small
$\epsilon>0$, and
\item  for such $\epsilon$, the fan
$\st_{\TSigma}(\rho)$ satisfies $\HR(\ell_\epsilon)$ if 
$\st_{\tallSigma}(\sigma)$ satisfies $\HR(\ell)$.
\end{enumerate}
\end{proposition}

\begin{proof}
Let $\Delta=\st_{\tallSigma}(\sigma)$ and $\TDelta=\st_{\TSigma}(\rho)$.
By Proposition~\ref{prop:PD}, we may assume both $\Delta$ and $\TDelta$
have Poincar\'e duality.  By Proposition \ref{lem:blowup K}, we have
$\ell_\epsilon\in \Kone(\TDelta)$ for small enough positive $\epsilon$.

If $k<(e+1)/2$, we
use the \HR\ property of $\ell\in \Kone(\Delta)$ and Proposition \ref{thm:blowup chow}
to obtain a decomposition
\[
A^k(\TDelta)=\PA^k(\Delta,\ell)
\oplus \ell A^{k-1}(\Delta)\oplus x_\rho A^{k-1}(\Delta),
\]
with respect to which $\hr^k(\TDelta,\ell_\epsilon)$
is represented by a block matrix
\[
\hr^k(\Delta,\ell_\epsilon)=
\begin{pmatrix}
H_{11}(\epsilon) & H_{12}(\epsilon) & H_{13}(\epsilon)\\
H_{21}(\epsilon) & H_{22}(\epsilon) & H_{23}(\epsilon)\\
H_{31}(\epsilon) & H_{32}(\epsilon) & H_{33}(\epsilon)\\
\end{pmatrix}.
\]

For any $\epsilon>0$, the matrix above is congruent to the matrix
\[
\overline{\hr}^k(\epsilon) \coloneq
\begin{pmatrix}
\epsilon^{-1}H_{11}(\epsilon)& \epsilon^{-1}H_{12}(\epsilon)& H_{13}(\epsilon)\\
\epsilon^{-1}H_{21}(\epsilon)& \epsilon^{-1}H_{22}(\epsilon)& H_{23}(\epsilon)\\
H_{31}(\epsilon) & H_{32}(\epsilon) & \epsilon H_{33}(\epsilon)\\
\end{pmatrix},
\]
the entries of which we will see are polynomial in $\epsilon$.
For elements $p_1,p_2\in \PA^k(\Delta,\ell)$, we have
\begin{align*}
\angl{p_1,p_2}_{\ell_\epsilon}&=(-1)^k \deg_{\TDelta}\big((\ell-\epsilon
x_\rho)^{e+1-2k}p_1p_2\big)\\
&=-(-1)^k\cdot\epsilon\cdot\deg_{\TDelta}
\big(\ell^{e-2k}(e+1-2k)p_1p_2x_\rho\big)+O(\epsilon^2)\\
&=(-1)^k\epsilon(e+1-2k)\deg_{\Delta}\big(\ell^{e-2k}p_1p_2\big)+O(\epsilon^2)\\
&=(e+1-2k)\epsilon\cdot\angl{p_1,p_2}_{\ell}+O(\epsilon^2).
\end{align*}
so the block $H_{11}(\epsilon)$ represents a positive multiple of
the pairing $\hr^k(\Delta,\ell)$, modulo $\epsilon^2$.

Similar computations show that
the block $H_{22}(\epsilon)$ is the matrix of the pairing
$(e+1-2k)\epsilon\cdot\hr^{k-1}(\Delta,\ell)$, modulo $\epsilon^2$,
and $H_{23}(\epsilon)=H_{32}(\epsilon)=-\hr^{k-1}(\Delta,\ell)$ modulo
$\epsilon$.  Along the same lines, we see $H_{12}(\epsilon)=H_{21}(\epsilon)$
are divisible by $\epsilon^2$, and $H_{13}(\epsilon)=H_{31}(\epsilon)$ is
divisible by $\epsilon$.  Returning to the matrix for $\overline{\hr}^k(\epsilon)$, we have
\[
\overline{\hr}^k(\epsilon) =
\begin{pmatrix}
(d+1-2k)\hr^k(\Delta,\ell)\mid_{\PA^k} & 0 & 0\\
0 & -(d+1-2k)\hr^{k-1}(\Delta,\ell) & -\hr^{k-1}(\Delta,\ell)\\
0 & -\hr^{k-1}(\Delta,\ell) & 0
\end{pmatrix}+O(\epsilon).
\]

Given our assumption that $k<(e+1)/2$,
the matrix $\overline{\hr}^k(0)$ is invertible, because each nonzero block
is nondegenerate (since $\ell$ has the \HL\ property).  It follows that
$\ell_\epsilon$ has the $\HL^k$ property for all $0\leq k<(e+1)/2$, for
some sufficiently small $\epsilon>0$.  Using Lemma~\ref{lem:schur}, we
see the signature of $\overline{\hr}^k(\epsilon)$ agrees with that of
the top-left block.  By hypothesis, $\hr^k(\Delta,\ell)$ is positive definite
on $\PA^k(\Delta,\ell)$.  Now $\dim\PA^k(\Delta,\ell)=
b_k(\Delta)-b_{k-1}(\Delta)$, which  by Propositions  \ref{thm:blowup chow} and \ref{thm:HR as signature}
is the expected signature for $\overline{\hr}^k(
\epsilon)$; that is, $\HR^k(\ell_\epsilon)$ holds for sufficiently small
$\epsilon$.

It remains to consider the case where $e$ is odd and $k=(e+1)/2$.  In this case we have
$A^k(\TDelta)=A^{k-1}(\Delta)\oplus x_\rho A^{k-1}(\Delta)$, and, up to a
sign, the pairing is equal to the Poincar\'e pairing $M^k(\TDelta)$.
In
the middle dimension, $M^k(\Delta)=M^{k-1}(\Delta)$, so we have a block
decomposition 
\[
M^k(\TDelta)=
\begin{pmatrix}
0 & -M^k(\Delta)\\
-M^k(\Delta) & Q
\end{pmatrix}
\]
for some square matrix $Q$.  The matrix $M^k(\Delta)$ is nonsingular,
by $\HL^k$, so $M^k(\TDelta)$ has signature zero by
Lemma~\ref{lem:schur}, which shows $\ell_\epsilon$ has $\HR^k$ for
any $\epsilon$ by by Propositions  \ref{thm:blowup chow} and \ref{thm:HR as signature}.
\end{proof}

\subsection{Proofs of the main results}\label{ss:blowup theorems}
We are now ready to prove the main result of this section.  We will treat the
star-shaped and ordinary cases separately, beginning with the former.
As before, let $\TSigma$ be a subdivision of a
$d$-dimensional simplicial fan $\Sigma$ by a ray $\rho$ contained in
a two-dimensional cone $\sigma$, and set
$\TDelta=\st_{\TSigma}(\rho)$ and $\Delta=\st_{\tallSigma}(\sigma)$.


\begin{theorem}\label{thm:lefschetz I}
The fan $\Delta$ is Lefschetz if and only if the fan $\TDelta$ is Lefschetz.
\end{theorem}

\begin{proof}
  First, suppose that $\Delta$ is Lefschetz, and let
  $\nu_1,\nu_2$ denote the two extreme rays of $\sigma$.  First, we check
  that the star of each cone $\tau\in\TDelta$ is Lefschetz.
  This is easy if $\tau$ does not contain $\nu_1$ or $\nu_2$, since
  then $\tau$ is a cone of $\Delta$.  Otherwise, $\tau$
  contains (exactly) one such ray, say $\nu_1$.  The remaining rays of
  $\tau$ span a cone $\tau'$ of $\Delta$, and by inspection,
  $\st_{\TDelta}(\tau)=  \st_{\Delta}(\tau')$, which is again Lefschetz by
  hypothesis.

Poincar\'e duality for $\TDelta$ follows from Proposition~\ref{prop:PD}.
  To establish \HL, we use Proposition~\ref{prop:HR=>HL+1}.  For this,
  we need to know that the star of each ray satisfies mixed \HR, but
  the star of a ray of $\TDelta$ is also a star in $\Delta$, so \HL\ for
  $\TDelta$ follows.
  Finally, we use Proposition~\ref{prop:star-shape HR}:
  for any $\ell\in \Kone(\Delta)$, there exists some $\ell_\epsilon\in
  \Kone(\TDelta)$
  with the \HR\ property.  By Corollary~\ref{cor:HR+1}, $\TDelta$ has \HR.

Conversely, if $\TDelta$ is Lefschetz, then
  $\st_{\TDelta}(\nu_1)=\Delta$, so $\Delta$ is Lefschetz too.
 \end{proof}
 
We note that $\Kone(\Delta)$ is nonempty if and only if
$\Kone(\TDelta)$ is nonempty.  The forward implication follows immediately from
Proposition \ref{lem:blowup K}.  The converse holds
because  $\Delta$ is a star in $\TDelta$.
However,  $\Kone(\TSigma)$ can be nonempty  while  $\Kone(\Sigma)$ is empty.

\begin{theorem}\label{thm:lefschetz II}
  If $\Sigma$
  is a Lefschetz fan with nonempty $\Kone(\Sigma)$,
  then $\TSigma$ is a Lefschetz fan.
  Conversely, if $\TSigma$ is a Lefschetz fan, then $\Sigma$ is a Lefschetz fan.
\end{theorem}

\begin{proof}
  We prove the first statement by induction on the dimension $d$.  The statement is
  vacuously true if $d=1$, so let us assume it holds for all Lefschetz
  fans of dimension less than $d$.

  First we check that the star of every cone $\tau\in\TSigma$ is Lefschetz,
  for which we  consider two cases.  First suppose
  $\tau\in \Sigma$.  If $\sigma\not\in\st_{\tallSigma}(\tau)$, then
  $\st_{\TSigma}(\tau)=\st_{\tallSigma}(\tau)$, which is Lefschetz.  If, on the
  other hand, $\sigma\in\st_{\tallSigma}(\tau)$, then
  $\st_{\TSigma}(\tau)=\stellar_\rho(\st_{\tallSigma}(\tau))$, which is
  a star-shaped subdivision.  Since $\st_{\tallSigma}(\tau)$ is Lefschetz, so is
  $\st_{\TSigma}(\tau)$, by Theorem~\ref{thm:lefschetz I}.

  Now suppose $\tau\not\in\Sigma$.  Then $\rho\in\tau$, so
  $\st_{\TSigma}(\tau)\subseteq \st_{\TSigma}(\rho)$: in fact,
  $\st_{\TSigma}(\tau)=\st_{\Sigma'}(\tau)$, where
  $\Sigma'=\st_{\TSigma}(\rho)$.  Since $\Sigma'=
  \stellar_\rho(\st_{\tallSigma}(\sigma))$, a star-shaped subdivision, $\Sigma'$
  is Lefschetz by Theorem~\ref{thm:lefschetz I},
  and it follows that $\st_{\TSigma}(\tau)$ is Lefschetz too.

  By Propositions~\ref{prop:PD} and \ref{prop:HR=>HL+1}, respectively,
  the fan $\TSigma$ satisfies \PD\ and \HL. It remains to check that
  $\TSigma$ satisfies \HR\ as well.

  Consider any $0\leq k\leq d/2$ and $\ell\in \Kone(\Sigma)$.
  By Lemma~\ref{lem:codim2 ortho},
  we have $\hr^k(\TSigma,\bl^*\ell)=\hr^k(\Sigma,\ell)\oplus
  \hr^{k-1}(\st_{\tallSigma}(\sigma),i_\sigma^*(\ell))$.  The summands are
  nondegenerate, because $\Sigma$ and $\st_{\tallSigma}(\sigma)$ satisfy
  $\HL(\ell)$ and $\HL(i^*_\sigma \ell)$, respectively, so
  $\hr^k(\TSigma,\bl^*\ell)$ is nondegenerate as well.

  By the \HR\ signature test (Proposition \ref{thm:HR as signature})
  the signature of $\hr^k(\TSigma,\bl^*\ell)$ equals
\begin{align*}
  &\sum_{i=0}^k (-1)^{k-i}\big(
  b_i(\Sigma)-b_{i-1}(\Sigma)\big) +
  \sum_{i=0}^{k-1} (-1)^{k-(i-1)}\big(
  b_{i-1}(\st_{\tallSigma}(\sigma))-b_{i-2}(\st_{\tallSigma}(\sigma))\big)\notag \\
  =& \sum_{i=0}^k (-1)^{k-i}\big(b_i(\Sigma)+b_i(\st_{\tallSigma}(\sigma))-
  b_{i-1}(\Sigma)-b_{i-1}(\st_{\tallSigma}(\sigma))\big)\notag \\
  =& \sum_{i=0}^k (-1)^{k-i}\big(b_i(\TSigma)-b_{i-1}(\TSigma)\big). 
  \end{align*}
  Proposition \ref{lem:blowup K} states $\bl^*\ell$ is in the closure of  $\Kone(\TSigma)$.
  Then there exists an open ball $U\subseteq A^1(\TSigma)$
  containing $\bl^*\ell$ on which $\hr^k(\TSigma,-)$ is nondegenerate.
  Choosing any $\ell'\in U\cap \Kone(\TSigma)$, we can use
  Corollary~\ref{cor:HR+1} to conclude that $\TSigma$ satisfies $\HR^k$.

The converse is similar in spirit. 
  Again, we argue by induction on the dimension $d$.  The base case being trivial,
  we assume that, if $\TSigma$ is Lefschetz and has dimension less than $d$,
  then $\Sigma$ is Lefschetz as well.  Now assume $\TSigma$ is a Lefschetz
  fan of dimension $d$, and we show $\Sigma$ is as well.

  \PD\ for $\Sigma$ follows from Proposition~\ref{prop:PD}.  Next, consider
  a ray $\nu\in\Sigma(1)$.  If $\nu\not\in\clst_{\Sigma}(\sigma)(1)$, then
  $\st_{\tallSigma}(\nu)=\st_{\TSigma}(\nu)$, which is Lefschetz.
  If, on the other hand,
  $\nu\in\clst_{\Sigma}(\sigma)(1)$, then $\sigma\in
  \clst_{\Sigma}(\nu)(2)$, and $\clst_{\TSigma}(\nu)=
  \stellar_\rho(\clst_{\tallSigma}(\nu))$.  Since $\st_{\TSigma}(\nu)$
  is Lefschetz, so is $\st_{\tallSigma}(\nu)$, by Theorem~\ref{thm:lefschetz I}.
  Either way,
  $\st_{\tallSigma}(\nu)$ has the \HR\ property for each ray $\nu$,
  so $\Sigma$ has the \HL\ property (by Proposition~\ref{prop:HR=>HL+1}).

  A similar argument shows that $\st_{\tallSigma}(\tau)$ is Lefschetz for all
  cones $\tau$ of $\Sigma$: if the star remains a star in $\TSigma$, it is
  Lefschetz by hypothesis.  Otherwise, a subdivision of it is a star in
  $\TSigma$.  If $\tau=\sigma$, the subdivided edge, we invoke
  Theorem~\ref{thm:lefschetz I}.  Otherwise, we note the dimension is
  less than $d$, so $\st_{\tallSigma}(\tau)$ is Lefschetz by induction.

  It remains to establish $\HR^k$ for $\Sigma$,
  for $0\leq k\leq d/2$.  The condition is vacuous if $\Kone(\Sigma)=\varnothing$.
  Otherwise, choose any $\ell\in \Kone(\Sigma)$.
  By Lemma~\ref{lem:codim2 ortho},
  \begin{equation*}
  \hr^k(\TSigma,\bl^*\ell)=\hr^k(\Sigma,\ell)\oplus
  \hr^{k-1}(\st_{\tallSigma}(\sigma),i^*_\sigma(\ell)).\label{eq:hr sum}
  \end{equation*}
  Since the second factor is the blowdown of
  $\st_{\TSigma}(\rho)$, it is Lefschetz by Theorem~\ref{thm:lefschetz I},
  and the first factor
  is Lefschetz by the argument above.  So both summands are nondegenerate,
  and so is $\hr^k(\TSigma,\bl^*\ell)$.

  By \HR, the bilinear form $\hr^k(\TSigma,\widetilde{\ell})$ has the expected
  signature for all $\widetilde{\ell}\in \Kone(\TSigma)$.  It follows by
  Proposition~\ref{prop:fixed sig} that
  $\hr^k(\TSigma,\bl^*\ell)$ also has that signature, since it is nondegenerate
  and $\bl^*\ell$ lies in the boundary of $\Kone(\TSigma)$.

  The \HR\ property for $\st_{\tallSigma}(\sigma)$ determines the signature of
  $\hr^{k-1}(\st_{\tallSigma}(\sigma),i_\sigma^*(\ell))$, and we obtain the
  signature of $\hr^k(\Sigma,\ell)$ by subtraction.  
  By the HR signature test again, we find that it equals $\sum_{i=0}^k(-1)^{k-i}
  \big(b_i(\Sigma)-b_{i-1}(\Sigma)\big)$, and we conclude $\Sigma$ has the
  $\HR^k$ property.
\end{proof}
Putting the pieces together gives a proof that the Lefschetz property is an
invariant of the support of a fan.

\begin{reptheorem}{thm:lefschetz}
  Let $\Sigma_1$ and $\Sigma_2$ be simplicial fans that have the same support
  $\abs{\Sigma_1}=\abs{\Sigma_2}$.
If $\Kone(\Sigma_1)$ and $\Kone(\Sigma_2)$
  are nonempty, then $\Sigma_1$ is Lefschetz if and only if $\Sigma_2$ is
  Lefschetz.
\end{reptheorem}

\begin{proof}[Proof of Theorem~\ref{thm:lefschetz}]
  Suppose $\abs{\Sigma}=\abs{\Sigma'}$.  According to Theorem~\ref{thm:fan MFT},
  there is a sequence of fans $(\Sigma_0,\Sigma_1,\cdots,
  \Sigma_N)$ with $\Sigma=\Sigma_0$, $\Sigma_N=\Sigma'$, and
  for which either $\Sigma_i\to\Sigma_{i+1}$
  or $\Sigma_{i+1}\to\Sigma_i$ is an edge subdivision, for each $i$.
  Furthermore, there is some $i_0$ for which $\Sigma_i\to \Sigma$ is a
  projective map of fans for each $i\leq i_0$, and $\Sigma_i\to\Sigma'$
  is a projective map of fans for each $i\geq i_0$.  By
  Proposition~\ref{prop:proj morphism}, we see that the cone $\Kone(\Sigma_i)$
  is nonempty for each $i$.  By Theorem~\ref{thm:lefschetz II},
  if any one of these fans is Lefschetz, then they all are.
\end{proof}



In our terminology, the main result of \cite{AHK15} says that the Bergman fan of ${\M}$ is Lefschetz.
We use the result to show that the conormal fan of $\M$ is Lefschetz.

\begin{lemma}\label{lem:lef products}
If $\Sigma_1$ and $\Sigma_2$ are Lefschetz fans, then so is $\Sigma_1\times
\Sigma_2$.
\end{lemma}
\begin{proof}
  It was shown in \cite[Section 7.2]{AHK15} that, if $\Sigma_1$ and $\Sigma_2$
  have \PD, \HL, and \HR, then so does $\Sigma_1\times\Sigma_2$.  Since
  stars of cones in a product are products of stars in the factors,
  we conclude that $\Sigma_1\times\Sigma_2$ is a Lefschetz fan, by induction
  on dimension.
\end{proof}



\begin{theorem}\label{thm:conormal Lefschetz}
For any matroid $\M$, the conormal fan $\Sigma_{\M,\M^\perp}$   is Lefschetz.

\end{theorem}

\begin{proof}
We may assume that $\M$ is loopless and coloopless. 
Since the Bergman fan is Lefschetz, from Lemma~\ref{lem:lef products} we
see the fan $\Sigma_{\vphantom{\M^\perp}\smash{\M}}\times\Sigma_{\M^\perp}$ is Lefschetz.  Moreover,
its support is equal to that of $\Sigma_{\M,\M^\perp}$.  Bergman fans are
quasiprojective, since they are subfans of the permutohedral fan, so
$\Kone(\Sigma_{\M}\times\Sigma_{\M^\perp})$ is nonempty.
We saw that the bipermutohedral fan $\Sigma_{E,\E}$ is the
normal fan of the bipermutohedron, so the conormal fan is also
quasiprojective, and $\Kone(\Sigma_{\M,\M^\perp})$ is nonempty as well.  By
Theorem~\ref{thm:lefschetz}, then, $\Sigma_{\M,\M^\perp}$ is Lefschetz.
\end{proof}


The extra structure present in the Chow rings of Lefschetz fans leads
easily to an Aleksandrov--Fenchel-type inequality.
\begin{theorem}\label{thm:AFinequality}
  Let $\Sigma$ be a Lefschetz fan of dimension $d$, and
  $\ell_2,\ell_3,\ldots,\ell_d$ elements in the closure of $\Kone(\Sigma)$.  Then
  for any $\ell_1\in A^1(\Sigma)$,
\[
  \deg(\ell_1\ell_2\cdots \ell_d)^2\geq
  \deg(\ell_1\ell_1\ell_3\cdots \ell_d)\cdot \deg(\ell_2\ell_2\ell_3\cdots
  \ell_d).
  \]
\end{theorem}
\begin{proof}
  We first verify the inequality when $\ell_i\in \Kone(\Sigma)$ for each
  $2\leq i\leq d$.  For this,
  let $\lef=\ell_3\cdots\ell_d$, a Lefschetz element, and consider
  $\angl{-,-}\coloneq \angl{-,-}_{\lef}$ on $A^1(\Sigma)$.

  If $\angl{\ell_2,\ell_2}\neq 0$, let $\ell'_1=\ell_1-
  \frac{\angl{\ell_1,\ell_2}}{\angl{\ell_2,\ell_2}}\ell_2$,
  so that $\angl{\ell'_1,\ell_2}=0$.  This means $\ell'_1\in
  \PA^1(\Sigma,\ell_2)$, so by \HR,
\[
0  \leq \angl{\ell'_1,\ell'_1}
 = \angl{\ell_1,\ell'_1}
 = \angl{\ell_1,\ell_1}-
\frac{\angl{\ell_1,\ell_2}}{\angl{\ell_2,\ell_2}}\angl{\ell_1,\ell_2}.
\]
By the signature test, $\angl{-,-}$ is negative-definite on the
orthogonal complement of $\ell_1'$. Therefore $\angl{\ell_2,\ell_2}<0$,
and we see that
  \[
  \angl{\ell_1,\ell_2}^2\geq \angl{\ell_1,\ell_1}
  \cdot \angl{\ell_2,\ell_2}.
  \]
  If, on the other hand,
$\angl{\ell_2,\ell_2}=0$, then the displayed inequality is obvious.

If we relax the hypothesis to consider
$\ell_2,\ldots,\ell_d$ in the closure of $\Kone(\Sigma)$, then the desired inequality
continues to hold by continuity, as in
\cite[Theorem 8.8]{AHK15}.
\end{proof}

\begin{reptheorem}{thm:conjs}
For any matroid $\M$, the $h$-vector of the broken circuit complex of $\M$ is log-concave.
\end{reptheorem}

\begin{proof}
It suffices to assume that $\M$ is loopless and coloopless.
  The
  classes $\gamma=\gamma_i$
  and $\delta=\delta_i$ are pullbacks of the nef classes $\alpha = \alpha_i \in A^1(\Sigma_{\M})$ and $\alpha=\alpha_i \in A^1(\Delta_E)$,
  along the two maps  $\pi\colon \Sigma_{\M,\M^\perp}\to \Sigma_{\M}$ and $\mu\colon \Sigma_{\M,\M^\perp}\to \Delta_E$, respectively.
  The pullback of a convex function on a fan is convex, so both $\gamma$
  and $\delta$ are represented by convex functions on the conormal fan.  Since
  $\Kone(\Sigma_{\M,\M^\perp})$ is nonempty by Proposition \ref{thm:ardilahedron}, we see that
  $\gamma$ and $\delta$ are in the closure of $\Kone(\Sigma_{\M,\M^\perp})$, following the discussion
  at the end of Section \ref{ss:K}.
  By Theorem~\ref{degtheorem}, we have
\[
h_{r-k}(\bc(\M)) = \deg_{\Sigma_{\M,\M^\perp}}(\gamma^k \delta^{n-k-1})
= \angl{\gamma,\delta}_{\lef},
\]
where $\lef=\gamma^{k-1}\delta^{n-k-2}$.  Since
$\Sigma_{\M,\M^\perp}$ is Lefschetz by Theorem~\ref{thm:conormal Lefschetz}, 
the log-concave inequalities follow from Theorem~\ref{thm:AFinequality}.
\end{proof}

\begin{remark}\label{rem:shortcut}
In the above proof of Theorem \ref{thm:conjs}, our use of the existence of the bipermutohedron (Proposition \ref{thm:ardilahedron}) can be avoided.
The toric Chow lemma \cite[Theorem 6.1.18]{CLS} guarantees that the conormal fan has a refinement that is the normal fan of a polytope, and we may apply the same argument to that refinement.
\end{remark}

\bibliographystyle{amsalpha}
\bibliography{refs}
\end{document}